\documentclass[10pt]{ijnam}
\def\currentvolume{}
\def\currentissue{}
\def\currentyear{2026}
\def\month{May}
\pagespan{1}{38}
\copyrightinfo{2026}{} % copyright year

\usepackage{subfig}
\usepackage{lipsum}
\usepackage{amsfonts}
\usepackage{graphicx}
\usepackage{algorithmic}
\usepackage{xcolor}
\usepackage{multirow}
\usepackage{amsthm}
\usepackage{algorithm}
\newtheorem{proposition}{Proposition}
\usepackage{amsthm}
\newtheorem{theorem}{Theorem}
\usepackage{url}
\begin{document}
%%%%%%%%%%%%%%%%

\title[Two-stage Multi-Grade Deep Learning for Nonlinear PDEs]
{Two-stage Multi-Grade Deep Learning for Nonlinear Partial Differential Equations with Applications to the Burgers Equation}

% author and address
\author[Y. Xu]{Yuesheng Xu}
\address{
  Department of Mathematics and Statistics, 
  Old Dominion University, 
  Norfolk, 
  Virginia 23529, 
  USA 
}
\email{y1xu@odu.edu}
%\urladdr{http://www.ucs.lousiana.edu/$\sim$hxy6615/}

% two authors have the same address
\author[T. Zeng]{Taishan Zeng*}
\address{
  School of Mathematical Science, 
  South China Normal University, 
  Guangzhou 510631, 
  China
}
\email{zengtsh@m.scnu.edu.cn}
%\urladdr{http://www.ucs.lousiana.edu/$\sim$hxy6615/ \and
%  http://www.ucs.lousiana.edu/$\sim$hxy6615/}

% dedication
%\dedicatory{This paper is dedicated to our authors}

% communicated
%\commby{Hilbert}

% Received by the editors ?
\date{}

% it is suggested to put it in Acknowledgments
\thanks{Y. Xu is supported in part by the US National Science Foundation under grant DMS-2208386.}
\thanks{*Corresponding author.}


\abstract{Deep neural networks (DNNs) show great promise for solving partial differential equations (PDEs), but their deep architectures introduce complex, large-scale, non-convex optimization challenges. Nonlinear PDEs, like the viscous Burgers' equation, compound these difficulties due to steep gradients and shock-like solutions. To address this, we propose a two-stage multi-grade deep learning (TS-MGDL) method. In the first stage, shallow networks are trained progressively grade by grade to fit the target function from low- to high-frequency components; previously learned grades are frozen, and each new residual block is trained solely to minimize the remaining approximation error. The second stage unfreezes and retrains selected layers using the first-stage network as initialization, achieving an interpretable, stable hierarchical refinement while mitigating optimization complexity. Furthermore, we theoretically prove that each grade and stage in TS-MGDL monotonically reduces the loss function under an appropriate optimization strategy. Numerical experiments on 1D, 2D, and 3D viscous Burgers' equations demonstrate that TS-MGDL significantly outperforms single-grade learning (SGL), reducing predictive errors by up to a factor of 60.}

%\abstract{Deep neural networks (DNNs) have attracted growing attention for their promising capabilities in solving partial differential equations (PDEs). However, deep networks introduce the challenge of solving large-scale, non-convex, and nonlinear optimization problems. Nonlinear PDEs, such as the viscous Burgers equation, pose additional challenges due to steep gradients and shock-like features in their solutions. To address these difficulties, we propose a two-stage multi-grade deep learning (TS-MGDL) method. The core idea is to progressively train the neural network from shallow to deep layers, gradually fitting the target function from low- to high-frequency components.  In the first stage, the DNN is trained incrementally, grade by grade, with each grade learning from the residuals of the previous one. To further minimize prediction error, the second stage unfreezes and retrains selected layers leveraging the network obtained from stage one as the initiation. This two-stage approach mitigates the complexity of optimizing a network with a large number of parameters while effectively capturing residual components left unresolved in earlier grades. Furthermore, we prove that each grade and stage in TS-MGDL reduces the loss function when using an appropriate optimization strategy. Numerical experiments on 1D, 2D, and 3D viscous Burgers equations show that TS-MGDL significantly outperforms single-grade learning (SGL), achieving predictive errors up to 60 times smaller.}


\maketitle

\section{Introduction}
Over the past decade, deep neural networks (DNNs) have achieved remarkable success across various scientific and technological fields. In domains such as natural language processing, computer vision, and robotics, deep learning has outperformed traditional artificial intelligence approaches \cite{Goodfellow2016, vaswani2017attention}. By leveraging multiple layers, DNNs can effectively extract multi-scale features and representations from input data.
Beyond traditional machine learning applications, DNNs have found widespread use in diverse areas, including physics \cite{diefenthaler2022deeply, raissi2018deep}, scientific computing \cite{cuomo2022scientific}, and finance \cite{deng2021crude}. The success of DNNs is largely attributed to their exceptional expressiveness in function representation. Consequently, understanding the mathematical foundations of DNNs has emerged as a key area of research in applied mathematics \cite{daubechies2022nonlinear, poggio2017why, shen2020deep, xu2022convergence, Xu_Zhang2024IEEE, ZHOU2020787, ZouBalanSingh2020}.

Partial differential equations (PDEs) play a crucial role in modern science and engineering, with applications spanning numerous disciplines. However, obtaining accurate solutions for nonlinear PDEs remains a challenging and active area of research. Numerical methods for solving PDEs can generally be divided into mesh-based and meshless approaches.
Mesh-based methods, including finite difference \cite{cheng2017dispersion, qiu2003finite}, finite element \cite{brenner2008mathematical}, finite volume schemes \cite{chen2012higher, chen2015construction, Li2000Generalized}, and spectral methods \cite{shen2011spectral}, rely on the creation of grids or meshes to discretize the problem domain. While these techniques have achieved substantial success, they often struggle with complex geometries and irregular domains, necessitating elaborate mesh generation. Moreover, mesh-based methods are prone to the curse of dimensionality, significantly limiting their effectiveness for high-dimensional problems.
Meshless methods, including techniques such as radial basis functions \cite{babuška_banerjee_osborn_2003} and, more recently, DNN-based approaches, eliminate the need for mesh generation. Unlike traditional meshless methods that depend on manually selected basis functions, DNNs automatically learn critical features from data, making them highly adaptable to a wide range of problems without requiring hand-designed functions.

DNN-based methods have gained significant attention for solving PDEs numerically. These approaches represent PDEs as constraints or loss functions that the neural network is trained to minimize. DNNs serve as universal solvers capable of handling diverse problems and are particularly effective in high-dimensional spaces and complex domains.
Recent advances in DNN-based methods include physics-informed neural networks (PINNs) \cite{PINN_2019}, the deep Ritz method \cite{E2018TheDR}, the deep Galerkin method \cite{sirignano2018dgm}, convolutional neural networks \cite{gao2021phygeonet}, generative adversarial networks \cite{bao2020numerical}, multi-scale DNNs \cite{wang2020multi, wang2021eigenvector}, and sparse DNNs \cite{XuZeng2023}. These methods have been successfully applied in fields such as fluid dynamics, structural mechanics, and molecular dynamics.

Despite their remarkable success in solving PDEs, DNNs still face several significant challenges. Typically, a DNN solves a single optimization problem to train all its parameters. As the network depth increases, solving the corresponding minimization problem becomes progressively more difficult, demanding larger datasets and greater computational resources.
One direct consequence of these optimization challenges is spectral bias, also known as the frequency principle. This inherent bias causes DNNs to learn low-frequency components of functions more quickly and easily than high-frequency components or steep gradients. This preferential learning significantly slows the convergence of high-frequency features, thereby limiting the network's ability to accurately capture fine details and rapid variations within the solution. Such limitations are particularly problematic for scientific and engineering problems that inherently involve high-frequency phenomena or rapid variations.

Recently, there has been considerable interest in improving the training of DNNs \cite{bengio2007greedy, haber2017stable}. Layer-wise training has emerged as a commonly used technique to address these challenges by training only a subset of layers at a time. One specific approach involves training a few initial layers and gradually adding new layers while keeping the previously trained layers fixed \cite{bengio2007greedy}. However, this method is susceptible to getting trapped in local minima.
The incorporation of residual connections into network architectures has been a highly effective strategy for improving the training of deep networks \cite{he2016deep, wang2024piratenets}. Nevertheless, as networks grow deeper, solving the resulting large-scale optimization problems remains a formidable challenge. Moreover, for multi-scale problems, relying on a single neural network is often insufficient to capture information across diverse scales.

Addressing spectral bias in deep learning involves various methods aimed at enhancing the model’s capacity to approximate high-frequency components effectively. One approach is Fourier feature mappings, which project the input onto sine and cosine functions at different frequencies, thus expanding the model’s ability to capture high-frequency information \cite{tancik2020fourier}. Multiscale DNNs \cite{wang2020multi} apply transformations across multiple scales to encompass a wider range of frequencies. Adaptive activation functions \cite{Jagtap2020} facilitate this process by enabling the model to adjust dynamically to different frequency scales during training. Phase DNNs \cite{Cai2020} utilize frequency-shifting techniques, simplifying the learning of high-frequency components. The multi-grade deep learning (MGDL) method, proposed in \cite{xu2023multi, xu2023successive}, builds upon the principle of residual learning. In MGDL, each ``grade" involves adding new layers to the existing network while keeping previous parameters fixed as known features. This allows the model to learn the residual of the current approximation, thereby enhancing its ability to capture complex features, including high-frequency components of the residual.
Despite these advancements, the challenge of addressing the spectral bias problems persists, underscoring the need for more robust and generalizable solutions.

In this paper, we introduce a two-stage multi-grade deep learning (TS-MGDL) method to address the limitations of DNNs, particularly when solving nonlinear PDEs. Because the original multi-grade deep learning (MGDL) approach relies on greedy training, while it proved effective for linear problems where the residual could be precisely determined, it may encounter difficulties in achieving high-accuracy numerical solutions when dealing with the complexities of nonlinearity. To prevent the MGDL method from getting trapped in an undesirable local minimizer, we propose a two-stage training strategy. In the first stage, the multi-grade deep network is trained incrementally, grade by grade, with each new grade learning additional features. In the second stage, to achieve higher prediction accuracy, we unfreeze the parameters of certain layers from the last several grades and retrain them over a broader range. Furthermore, we prove that each grade and stage in TS-MGDL reduces the loss function when choosing an appropriate optimization strategy.

%The MGDL method was previously proposed in \cite{xu2023multi, xu2023successive} in a general setting. The main idea of the MGDL is based on residual learning, where each grade stacks additional layers on top of the previous grade (with fixed parameters) to learn the residual of the approximation. 

%Recent literature has explored similar strategies of using successive neural networks to learn from the residuals of previous networks. For instance, approaches such as stacked networks \cite{howard2023stacked} and multifidelity deep operator networks \cite{howard2023multifidelity} have been proposed to improve physics-informed training and problem-solving capabilities. Multi-stage and multi-level neural networks, as discussed in \cite{wang2024multi} and \cite{aldirany2024multi}, are used to enhance accuracy in solving boundary-value problems and function approximation. Additionally, precision machine learning techniques relevant to function approximation have been explored \cite{michaud2023precision}. In contrast, our proposed TS-MGDL method introduces a novel two-stage training strategy that not only integrates multi-grade deep learning with residual learning but also systematically addresses the issue of local minima through incremental training and selective parameter retraining, thus offering a more robust and effective approach for solving PDEs.

The Burgers equation is a fundamental equation in nonlinear dynamics, combining both nonlinear propagation effects and diffusive phenomena \cite{whitham2011}. It is well-known for its potential to develop a jump discontinuity, or shock wave, even from a smooth initial condition. Due to its inherent nonlinearity and sharp gradients, solving the Burgers equation presents significant challenges, making it a benchmark problem for testing PDE solvers \cite{LuLiuShu2021}. This study focuses on the Burgers equation due to its complexity and relevance in various applications such as nonlinear wave propagation, turbulence, and shock waves.

%%%%%%%%%%%%%%%%%
 MGDL improves DNN training through a unique architectural and optimization strategy that distinguishes it from several established frameworks:
\begin{itemize}
    \item {\bf Residual Networks (ResNet):} While ResNets \cite{he2016deep} use skip connections to stabilize gradient flow in a single, jointly-optimized model, MGDL decomposes the problem. It builds depth incrementally by training shallow subnetworks on the residuals of previous grades, reducing optimization complexity by focusing on one grade at a time.

\item {\bf Relay Backpropagation:} Similar to MGDL, relay backpropagation \cite{shen2016relay}  partitions networks into segments with local losses. However, it optimizes all segments simultaneously. MGDL adopts a sequential freezing mechanism, where earlier subnetworks remain fixed, preventing the gradient interference common in joint multi-loss optimization.

\item {\bf Curriculum Learning (CL):} CL \cite{bengio2009curriculum, graves2017automated} stages training by increasing data difficulty over a fixed architecture. MGDL instead stages the architecture itself, decomposing the model into a sequence of shallow subnetworks. This simplifies optimization at the structural level rather than the data level.

\item {\bf Gradient Boosting (GB):} MGDL shares the stage-wise, residual-driven logic of GB \cite{chen2016xgboost, friedman2001greedy, ke2017lightgbm}. However, while GB typically creates a flat additive ensemble, MGDL constructs a hierarchical, compositional model. Each grade transforms features from the previous stage, allowing the system to learn deep, nested representations rather than just additive corrections.
\end{itemize}
In essence, MGDL's grade-wise sequential training of shallow subnetworks provides a stable alternative to end-to-end learning, facilitating the progressive construction of complex hierarchies.

%%%%%%%%%%%%%%
The organization of this paper is as follows: Section 2 reviews traditional DNN and PINN for solving PDEs. Section 3 introduces the proposed MGDL method for solving PDEs and provides theoretical justification. In Section 4, we discuss the two-stage training strategy for the MGDL method, demonstrating how the second stage can reduce the residue error of the approximate solution obtained from the first stage. Section 5 details the implementation of the TS-MGDL method to solve the Burgers equations in 1D, 2D, and 3D, validates the theoretical results, and compares them with those produced by SGL. Section 6 presents the numerical findings and discusses their implications.
Finally, we conclude our work in Section 7.

\section{Background}
In this section, we first describe the PDE setting and then review the conventional DNN, followed by a description of the PINN for solving PDEs.

%Although this paper focuses on the Burgers equation, 
We first describe the general initial-boundary value problem of nonlinear PDEs to be studied in this paper.
For an integer $d\geq 1$, we let $\Omega \subset \mathbb{R}^d$ be an open domain with boundary $\Gamma$. We denote by $\mathcal{F}: C^l(0,T]\times C^k(\Omega)\to C((0,T]\times \Omega)$ the nonlinear differential operator, by $\mathcal{I}: C(\Omega)\to C(\Omega)$ and $\mathcal{B}: C((0,T]\times \Omega)\to  C((0,T]\times \Omega)$ the (nonlinear) operators for the initial and boundary conditions, respectively.
Assuming that $u\in C^l(0, T]\times C^k(\Omega)$ is the unknown solution to be learned, we consider the initial-boundary value problem of nonlinear PDEs as follows
\begin{align}\label{eq:PDE}
	\mathcal{F} (u(t, x)) &= 0, \ \ x \in \Omega, \ t \in (0, T],  \\
	 \label{eq:PDE:initial}
	 \mathcal{I}(u(0,x)) &=0, \ \ x \in \Omega, \ t=0, \\
	\label{eq:PDE:boundary}
	 \mathcal{B}(u(t, x))  &= 0, \ \ x \in \Gamma, \ t \in (0, T], 
\end{align} 
where $T>0$, the initial state $u(0,x)$ and the data of $u$ on $\Gamma$ are given. The formulation \eqref{eq:PDE}-\eqref{eq:PDE:boundary} offers a comprehensive framework for many problems, ranging from the wave equation, Maxwell's equations, and the Burgers equation. For example, the 1D Burgers equation can be written as 
$
\mathcal{F}(u) := \frac{\partial}{\partial t} u + u\frac{\partial}{\partial x} u - \nu \frac{\partial^2}{\partial x^2} u,
$
where $\nu$ is the viscosity coefficient.

%We further assume that the initial-boundary value problem \eqref{eq:PDE}-\eqref{eq:PDE:boundary} satisfies the following assumption:

We recall the definition of a feed-forward neural network (FNN). Let $s,t$ be two positive integers.
An FNN with depth $D$,  a function that maps an input vector of $s$ dimensions to an output vector of $t$ dimensions, is
a neural network consisting of an input layer, $D-1$ hidden layers, and an output layer. 
Let $d_i$ denote the number of neurons in the $i$-th hidden layer, and let $W_i \in \mathbb{R}^{d_i \times d_{i-1}}$ and $b_i \in \mathbb{R}^{d_i}$ represent the weight matrix and the bias vector, respectively, for the $i$-th layer.  Let $\sigma : \mathbb{R} \rightarrow \mathbb{R}$ denote an activation function and  $\mathbf{x} :=[x_1, x_2, \dots, x_s]^T \in \mathbb{R}^s$ be the input vector.
The output of the first hidden layer, denoted as $\mathcal{H}_1(\mathbf{x})$, is defined by applying the activation function to an affine map defined by $W_1 \in \mathbb{R}^{d_1 \times s}$ and $b_1 \in \mathbb{R}^{d_1}$. That is,
$
\mathcal{H}_1(\mathbf{x}) := \sigma (W_1 \mathbf{x}+b_1)$, $\mathbf{x}\in \mathbb{R}^s$.
For a neural network with depth $D \ge 3$, the output of the ($i$+1)-th hidden layer can be identified as a recursion of the output of the $i$-th hidden layer, defined as  
$\mathcal{H}_{i+1}(\mathbf{x}) := \sigma (W_{i+1} \mathcal{H}_{i}(\mathbf{x}) + b_{i+1})$, for $i = 1, 2, \dots, D-2$.
Finally, the output of the neural network with depth $D$ is defined by
\begin{equation}
    \label{def:neural}
    \mathcal{N}_{D} (\mathbf{x}) := W_{D} \mathcal{H}_{D-1}(\mathbf{x}) + b_{D}.
\end{equation}
We denote by $\Theta:= \{W_i, b_i\}_{i=1}^{D}$ the set of trainable network parameters, consisting of all weight matrices and bias vectors for all layers, and specifically write $  \mathcal{N}_{D} (\Theta:\mathbf{x}):=\mathcal{N}_{D} (\mathbf{x})$, for $\mathbf{x}\in \mathbb{R}^s$.
Typically, the activation function is a nonlinear function, making the neural network capable of learning complex patterns and representations.

The physics-informed neural network (PINN) model \cite{PINN_2019} utilizes a deep neural network to learn the solution of a PDE with initial/boundary value conditions. %by incorporating governing PDE in the loss function. 
The loss function of PINN consists of three components: loss of the PDE, loss of the initial condition, and loss of the boundary condition. Assuming that $\mathcal{N}_{D}(\Theta; \mathbf{x})$, for $\mathbf{x}:=(t,x)\in (0,T]\times \Omega$, is a deep neural network to be learned,
the three components of the loss function of PINN are defined as follows:
\begin{itemize}
\item Loss of the PDE:
$
Loss_{PDE}(\mathcal{N}_{D}(\Theta; \bullet)) :=\frac{1}{N_f}\sum_{i=1}^{N_f}|\mathcal{F}(\mathcal{N}_{D}(t_f^i, x_f^i))|^2$, where $(t_f^i, x_f^i)\in (0, T] \times \Omega$ are collocation points obtained by using the Hammersley sampling method \cite{Sampling1997}.
\item
Loss of the initial condition: 
$
Loss_{I}(\mathcal{N}_{D}(\Theta; \bullet)) := \frac{1}{N_0}\sum_{i=1}^{N_0} |\mathcal{I}(\mathcal{N}_{D}(0, x_0^i))|^2$, where $x_0^i \in \Omega$ are randomly sampled for the initial condition \eqref{eq:PDE:initial}.
\item
Loss of boundary conditions:
$Loss_{B}(\mathcal{N}_{D}(\Theta; \bullet)) := \frac{1}{N_{b}} \sum_{i=1}^{N_{b}} \left|\mathcal{B}(\mathcal{N}_{D}(t_b^i, x_b^i) )\right|^2$, where $(t_b^i, x_b^i) \in (0, T] \times \Gamma$ are randomly sampled for the boundary condition \eqref{eq:PDE:boundary}.
\end{itemize}
We define the loss function (the residue error) by
\begin{equation}
\label{loss:PINN}
    Loss(\mathcal{N}_{D}(\Theta; \bullet)) := Loss_{PDE}(\mathcal{N}_{D}(\Theta; \bullet)) + Loss_{I}(\mathcal{N}_{D}(\Theta; \bullet)) + Loss_{B}(\mathcal{N}_{D}(\Theta; \bullet))
\end{equation}
and solve
\begin{equation}
    \label{eq:PINN}
    \min_{\Theta} Loss(\mathcal{N}_{D}(\Theta;\bullet))
\end{equation}
for an optimal $\Theta^*$. Then, the function $\mathcal{N}_{D}(\Theta^*;\bullet)$ provide us an approximate solution of the initial-boundary value problem \eqref{eq:PDE}-\eqref{eq:PDE:boundary}. Throughout this paper without further mentioning, we assume that minimization problem \eqref{eq:PINN} has a solution.
Selecting an appropriate optimization algorithm for \eqref{eq:PINN} is crucial for obtaining an accurate solution of PDE.
There are many optimization algorithms for solving deep learning problems, such as stochastic gradient descent \cite{Bottou2012_sgd},  Adam \cite{KingmaB14_adam}, L-BFGS \cite{Liu1989lbfgs}, and others.

To further generalize the optimization objective and address potential gradient imbalances between different terms, the loss function can be extended to a weighted form:
\begin{equation}
\label{loss:PINN_weighted}
\begin{split}
Loss_{w}(\mathcal{N}_{D}(\Theta; \bullet)) := & Loss_{PDE}(\mathcal{N}_{D}(\Theta; \bullet)) \\
& + \lambda_I Loss_{I}(\mathcal{N}_{D}(\Theta; \bullet)) + \lambda_B Loss_{B}(\mathcal{N}_{D}(\Theta; \bullet)),
\end{split}
\end{equation}
where $\lambda_I$ and $\lambda_B$ are positive penalty weights assigned to the initial and boundary constraints, respectively. When $\lambda_I = \lambda_B = 1$, this formulation reduces to the standard unweighted loss function as defined in \eqref{loss:PINN}. In this paper, we primarily employ the unweighted loss configuration ($\lambda_I = \lambda_B = 1$) as the default setting. The primary motivation for this choice is to demonstrate the intrinsic efficacy and robustness of the proposed multi-grade network architecture itself, rather than relying on intensive hyperparameter tuning of penalty weights. It is important to note, however, that the theoretical theorems and algorithms developed in the subsequent sections can be naturally and directly extended to the weighted loss scenarios.

It is important to point out that the traditional fully connected neural network, which involves training a neural network with a fixed architecture and optimizing all the parameters at once,  has limitations. Although it can work well for certain problems, it can become challenging to train DNNs with many layers due to issues such as the vanishing gradient problem. We address this issue in the next two sections.

\section{Multi-Grade Deep Learning}
In this section, we develop a multi-grade deep learning (MGDL) method, originally introduced in \cite{xu2023multi} for function approximation, for the numerical solution of the initial-boundary value problem \eqref{eq:PDE}-\eqref{eq:PDE:boundary}.   

%During the training process, we learn a neural network grade by grade, in a way that mimics the learning process of humans. 

The standard PINN  solves the optimization problem \eqref{eq:PINN} for $\Theta^*$ and uses the function $\mathcal{N}_{D}(\Theta^*;\bullet)$ as a numerical solution of the initial-boundary value problem \eqref{eq:PDE}-\eqref{eq:PDE:boundary}. Due to large number $D$ of layers, the optimization problem \eqref{eq:PINN} is highly non-convex and has large number of parameters. Solving it is a challenging task: Gradient-based methods for solving the optimization problem often get stuck in its poor solutions \cite{bengio2007greedy}. Note that gradient-based methods are likely to converge when solving problem \eqref{eq:PINN} with a small number of hidden layers \cite{Wu2005}. To address computational issues involved in problem \eqref{eq:PINN} with a large number of hidden layers,  instead of learning all parameters together, we solve several related optimization problems with less layers, grade by grade.

We start with %the multi-grade deep learning method originally proposed in \cite{xu2023multi}. 
a neural network $\mathcal{N}_{k_1}$ of depth $k_1$  defined in \eqref{def:neural}, with parameters $\Theta_1:= \{W_i^1, b_i^1\}_{i=1}^{k_1}$. We choose 
$
u_1 = u_1(\Theta_1;\bullet):=\mathcal{N}_{k_1}(\Theta_1;\bullet)
$ 
as the neural network of grade 1. To learn the parameters of grade 1, we solve the minimization problem %equation \eqref{eq:min:grade_1},
\begin{equation}
    \label{eq:min:grade_1}
    \min_{\Theta_1} Loss(u_1(\Theta_1;\bullet)),
\end{equation}
where $Loss$ is the loss function defined in \eqref{loss:PINN}. Upon solving this problem, we obtain an approximate solution of grade 1, denoted by $u_1^{*}= \mathcal{N}_{k_1}^{*}:= \mathcal{N}_{k_1} (\Theta_1^*;\bullet)$, where $\Theta_1^{*}:= \{W_i^{1*}, b_i^{1*}\}_{i=1}^{k_1}$ are the learned parameters.  Using the definition \eqref{def:neural},  the approximate solution $u_1$ can be represented as
$u_1^{*}(\mathbf{x}) = W_{k_1}^{1*} \mathcal{H}_{k_1 -1}^{1*} (\mathbf{x}) + b_{k_1}^{1*}$, 
where $\mathcal{H}_{k_1 -1}^{1*}$ is the network $\mathcal{N}_{k_1}^*$ without the output layer, $W_{k_1}^{1*}$ denotes the weight matrix connecting the last hidden layer with the output layer, and $b_{k_1}^{1*}$ represents the corresponding bias vector. Note that the neural network $\mathcal{H}_{k_1 -1}^{1*}$ has the learned parameters $\{W_i^{1*}, b_i^{1*}\}_{i=1}^{k_1-1}$. This formulation highlights the role of the hidden layers in transforming the input $\mathbf{x}$ into a feature representation captured by $\mathcal{H}_{k_1 -1}^{1*}$. 

\begin{figure}[ht]
  \centering
  \includegraphics[width=0.8\textwidth]{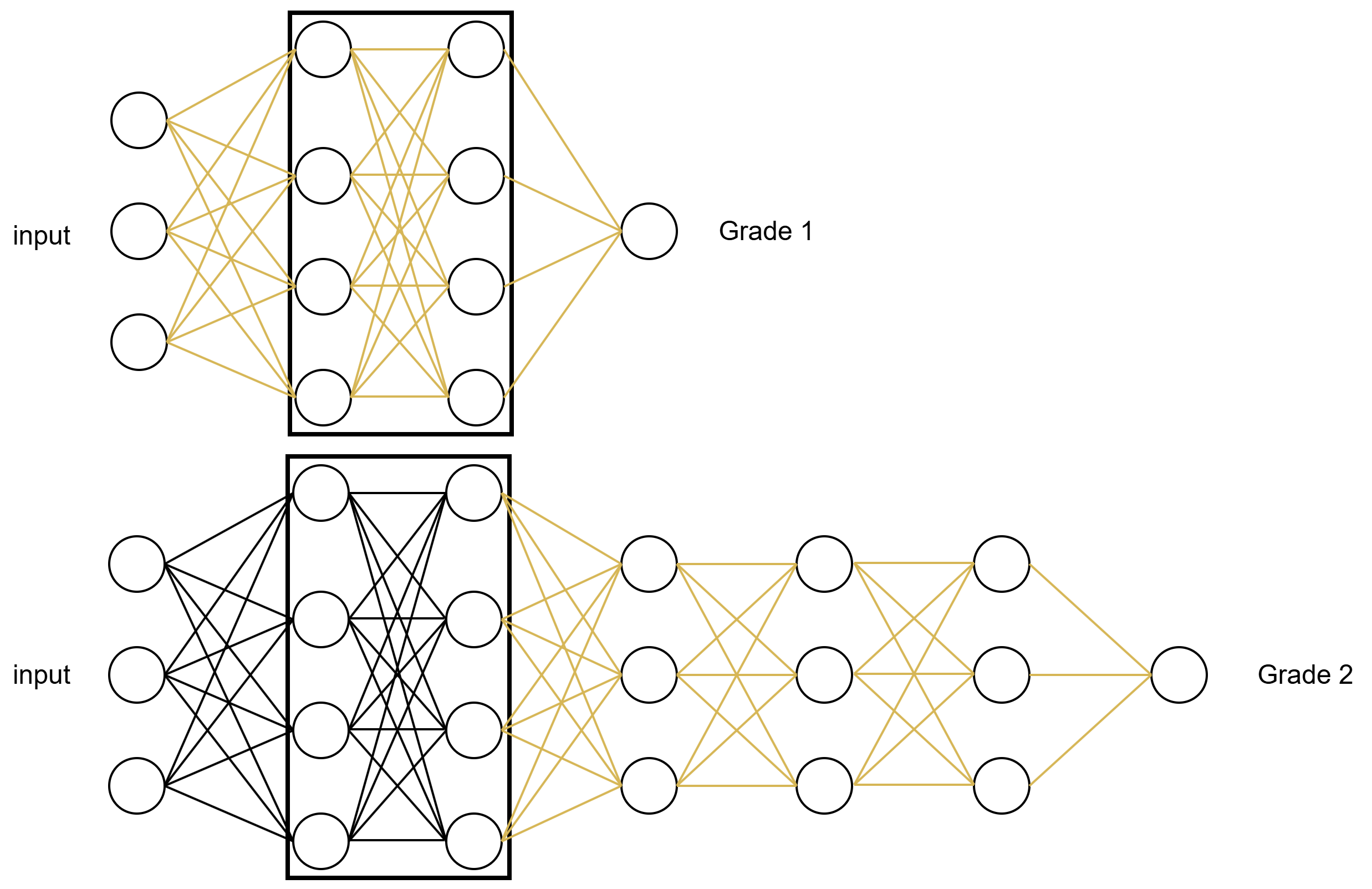}
  \caption{Multi-grade network with 2 grades}
  \label{fig:multi_grade_network}
\end{figure}

Next, we construct the neural network of grade 2, denoted by $u_2$, which is built on top of the neural network of grade 1, using a neural network  $\mathcal{N}_{k_2}:=\mathcal{N}_{k_2}(\Theta_2; \bullet)$ with parameters $\Theta_2=\{W_j^2, b_j^2\}_{j=1}^{k_2}$. 
The composition of two functions (or neural networks) $f: X \to Y$ and $g: Y \to Z$ is defined as $(g \circ f)(x) = g(f(x))$, meaning the output of $f$ becomes the input of $g$. 
To define the grade 2 network, we remove the output layer of grade 1 and stack $\mathcal{N}_{k_2}$ on its last hidden layer:
\begin{equation}\label{Def:u-2}
    u_2= u_2(\Theta_2;\bullet) := \mathcal{N}_{k_2}(\Theta_2; \bullet) \circ \mathcal{H}_{k_1 -1}^{1*},
\end{equation}
Equation \eqref{Def:u-2} indicates that the output of the last hidden layer of the grade 1 network, $\mathcal{H}_{k_1 - 1}^{1*}$, serves as the input to the grade 2 network $\mathcal{N}_{k_2}$. In other words, in grade 2, we use a new shallow network $\mathcal{N}_{k_2}$ composed with the shallow network $\mathcal{H}_{k_1 - 1}^{1*}$ trained in grade 1 as an ``adaptive basis'' to learn higher-scale information left in the residue of grade 1. In this way, on one hand, we only train a shallow network and on the other hand, the new shallow network composed with the old shallow network learned in grade enhances its expressiveness. 
We illustrate in Figure \ref{fig:multi_grade_network} the network structure of grades 1 and 2, with the trainable weight parameters highlighted in yellow. From the construction process of a multi-grade network, it can be seen that traditional fully connected neural networks can be viewed as single-grade learning (SGL).

We learn the optimal parameters $\Theta_2^*:=\{W_j^{2*}, b_j^{2*}\}_{j=1}^{k_2}$  of grade 2 by solving the minimization problem 
\begin{equation}
    \label{eq:min:grade_2}
    \min_{\Theta_2} Loss(u_1^{*}(\bullet)+u_2(\Theta_2;\bullet))
\end{equation} 
and as \eqref{Def:u-2} define
$
u_2^{*}:= u_2(\Theta_2^*;\bullet)= \mathcal{N}_{k_2}(\Theta_2^{*}; \bullet ) \circ \mathcal{H}_{k_1 -1}^{1*}.
$ 
Note that the parameters of $\mathcal{H}_{k_1 -1}^{1*}$ are fixed during the minimization process \eqref{eq:min:grade_2}. We can see from problem \eqref{eq:min:grade_2} that $u_2^{*}$ learns the residual of the solution $u_1^{*}$ of grade 1 to better approximate the solution of the PDE.

We repeat this process to construct a neural network of $\ell+1$ grades, assuming that the neural networks $u_i^*$ of grade $i$, $1 \le i \le \ell$, have been learned. Namely, we will learn a neural network $\mathcal{N}_{k_{\ell+1}}(\Theta_{\ell+1}; \bullet)$ with trainable parameters $\Theta_{\ell+1}=\{W_j^{\ell+1}, b_j^{\ell+1}\}_{j=1}^{k_{\ell+1}}$ to define the neural network $u_{\ell+1}$ of grade $\ell+1$. That is,
\begin{equation}
    \label{eq:grade:Lplus1}
    u_{\ell+1}(\Theta_{\ell+1};\mathbf{x}) := (\mathcal{N}_{k_{\ell+1}}(\Theta_{\ell+1}; \bullet) \circ \mathcal{H}_{k_{\ell} -1}^{\ell*} \circ \dots \circ \mathcal{H}_{k_2-1}^{2*} \circ \mathcal{H}_{k_1-1}^{1*})(\mathbf{x}),
\end{equation}
where  $\mathcal{H}_{k_i -1}^{i*}$ denotes the neural network  $\mathcal{N}_{k_{i}}^*:=\mathcal{N}_{k_{i}}(\Theta_i^{*}; \bullet )$ without the output layer, with the learned parameters $\{W_j^{i*}, b_j^{i*}\}_{j=1}^{k_i-1}$, $i = 1, 2, \dots, \ell$. 
The optimal parameters $\Theta_{\ell+1}^*=\{W_j^{\ell+1*}, b_j^{\ell+1*}\}_{j=1}^{k_{\ell+1}}$ of grade $\ell+1$ is obtained by solving the problem
\begin{equation}
    \label{eq:min:grade_L}
    \min_{\Theta_{\ell+1}} Loss\left( \sum_{i=1}^{\ell} u_i^{*}(\bullet)+u_{\ell+1}(\Theta_{\ell+1};\bullet)\right).
\end{equation}
It should be noted that the parameters of $\mathcal{H}_{k_i -1}^{i*}$, for $i=1, 2, \dots, \ell$, are all fixed during the minimization process \eqref{eq:min:grade_L}.
We then let $u_{\ell+1}^*:=u_{\ell+1}(\Theta_{\ell+1}^*;\bullet)$.  In training of grade $\ell+1$, the network $\mathcal{H}_{k_{\ell} -1}^{\ell*} \circ \dots \circ \mathcal{H}_{k_2-1}^{2*} \circ \mathcal{H}_{k_1-1}^{1*}$ serves as an adaptive basis to enhance the expressiveness of $u_{\ell+1}^*$.
%We learn the parameters of grade $\ell+1$ by solving the minimization problem given by equation \eqref{eq:min:grade_L} to obtain the optimal $u_{\ell+1}^{*}$.
Finally, by summing the approximations of all $\ell+1$ grades,  the resulting neural network 
\begin{equation}\label{Multiscale-Rep}
    \bar{u}_{\ell+1}^{*}: = \sum_{i=1}^{\ell+1} u_i^{*}
\end{equation}
serves as an approximate solution of the initial-boundary value problem \eqref{eq:PDE}-\eqref{eq:PDE:boundary}.

This strategy replaces a single large end-to-end optimization with a sequence of smaller minimization problems involving shallow networks. By successively learning the residuals of previous grades, MGDL progressively refines the solution and achieves higher accuracy with improved optimization stability. The effectiveness of this approach is supported by the theoretical result stated below. Moreover, the representation \eqref{Multiscale-Rep} naturally provides a multiscale decomposition of the solution. Such a form is particularly well suited for problems exhibiting sharp gradients or localized features, as the MGDL procedure extracts and captures the intrinsic multiscale structure of the PDE solution across successive grades.

Although in this work MGDL is implemented using fully connected networks, the same framework can be integrated with other architectures, such as residual networks, to further enhance their performance.

In the following, we use $\Theta_{\ell+1}^0$ to denote the parameters of the neural network $ u_{\ell+1} $ when the output layer is set to zero, defined as
$\Theta_{\ell+1}^0 := \{W_j^{\ell+1}, b_j^{\ell+1}\}_{j=1}^{k_{\ell+1}}$,
with
$W_{k_{\ell+1}}^{\ell+1} = 0$ and $b_{k_{\ell+1}}^{\ell+1} = 0$.

\begin{proposition}\label{thm:loss}
If $\bar{u}^*_{\ell}$ and $\bar{u}^*_{\ell+1}$ are the neural network solutions of the initial-boundary value problem \eqref{eq:PDE}-\eqref{eq:PDE:boundary} obtained by the MGDL method after $\ell$ and $\ell+1$ grades, respectively, then 
\begin{equation}\label{Justification}
    Loss\left( \bar{u}_{\ell+1}^{*}(\bullet)\right) \le Loss\left(\bar{u}_{\ell}^{*}(\bullet)\right).
\end{equation}  
\end{proposition}

\begin{proof}
It can be seen directly from definition \eqref{def:neural} of neural networks that if the weight matrix and bias vector of the output layer are both set to 0, then the output of a fully connected neural network will be 0. In other words, the resulting neural network is a zero function. Following this observation, if we initialize the weight matrix and bias vector of the output layer to 0, that is, by selecting $\Theta_{\ell+1}^0$,
then we observe that 
$u_{\ell+1}(\Theta_{\ell+1}^0;\mathbf{x}) = 0$, for all $\mathbf{x}\in\mathbb{R}^s$. 
Hence, the loss function will satisfy the equation
\begin{equation}\label{Justification0}
    Loss\left( \bar{u}_{\ell}^{*}(\bullet)+u_{\ell+1}(\Theta_{\ell+1}^0;\bullet)\right) = Loss\left(\bar{u}_{\ell}^{*}(\bullet)\right).
\end{equation} 

Let $\Theta_{\ell+1}^*$ be a solution of the minimization problem \eqref{eq:min:grade_L}. Then, for the initialized parameters $\Theta_{\ell+1}^0$ chosen as above, there holds that 
\begin{align*}
Loss\left(\bar{u}_{\ell+1}^{*}(\bullet)\right)&=
    Loss\left(\bar{u}_{\ell}^{*}(\bullet)+u_{\ell+1}(\Theta_{\ell+1}^*;\bullet)\right)\\
    &=\min_{\Theta_{\ell+1}}Loss\left(\bar{u}_{\ell}^{*}(\bullet)+u_{\ell+1}(\Theta_{\ell+1};\bullet)\right)\\
    &\leq  Loss\left(\bar{u}_{\ell}^{*}(\bullet)+u_{\ell+1}(\Theta_{\ell+1}^0;\bullet)\right)\\
    &= Loss\left(\bar{u}_{\ell}^{*}(\bullet)\right).
\end{align*} 
We have used equation \eqref{Justification0} in the last equality of the above inequality. This yields inequality \eqref{Justification}. 
\end{proof}

%{\bf TS: The next result requires additional hypothesis on the smoothness of the operators F, B, I. We need to carefully formulate it.}
%{\color{red}

Proposition \ref{thm:loss} guarantees that adding a new layer does not increase the loss. Under additional assumptions, we can further establish a strict decrease of the loss in grades. We employ the gradient descent algorithm (GDA) to demonstrate this strict decrease. The GDA iteratively updates parameters by moving in the opposite direction of the gradient, aiming to reduce the function’s value. Specifically, let $\Delta_{\ell+1}(\Theta_{\ell+1};\bullet):=\nabla_{\Theta_{\ell+1}} Loss\left(\bar{u}_{\ell}^{*}(\bullet)+u_{\ell+1}(\Theta_{\ell+1};\bullet)\right)$ denote the gradient of the loss with respect to the parameter ${\Theta_{\ell+1}}$.
According to the differentiability of loss function \( Loss\left(\bar{u}_{\ell}^{*}(\bullet)+u_{\ell+1}(\Theta_{\ell+1};\bullet)\right) \) with respect to parameters \( \Theta \), the update rule of the GDA at each step is as follows
\[
\Theta^{\text{new}} = {\Theta_{\ell+1}^0} - \eta \Delta_{\ell+1}(\Theta_{\ell+1};\bullet)|_{\Theta_{\ell+1}=\Theta_{\ell+1}^0},
\]
where \( \eta > 0 \) is the learning rate, controlling the step size in the parameter space.

\begin{proposition} \label{thm:strict_loss_decrease}
Suppose that the functions \(\mathcal{F}\), \(\mathcal{B}\), and \(\mathcal{I}\) are continuously differentiable with respect to \(u(\cdot)\), and that the neural network employs continuously differentiable activation functions. Let \(\Theta_{\ell+1}^0\) denote the initial parameter vector at layer \(\ell+1\). If the gradient of the loss  with respect to \(\Theta_{\ell+1}\) at \(\Theta_{\ell+1}^0\) is non-zero, i.e., \(\Delta_{\ell+1}(\Theta_{\ell+1};\bullet)|_{\Theta_{\ell+1}=\Theta_{\ell+1}^0} \neq 0\),
then the loss decreases strictly:
\begin{equation}\label{eq:strict_decrease}
    Loss\left( \bar{u}_{\ell+1}^{*}(\bullet) \right) < Loss\left( \bar{u}_{\ell}^{*}(\bullet) \right).
\end{equation}
\end{proposition}

\begin{proof} 
Since \(\mathcal{F}\), \(\mathcal{B}\), and \(\mathcal{I}\) are continuously differentiable with respect to \(u(\cdot)\), and the neural network’s activation functions are also continuously differentiable, it follows that the loss $Loss$ is continuously differentiable with respect to \(\Theta_{\ell+1}\).
%
%To demonstrate the strict decrease in the loss function, we consider the concept of directional derivatives. 
By hypothesis that \(\Delta_{\ell+1}(\Theta_{\ell+1};\bullet)|_{\Theta_{\ell+1}=\Theta_{\ell+1}^0} \neq 0\),
the directional derivative of the loss  at \(\Theta_{\ell+1}^0\) in the direction of the negative gradient is given by:
\begin{align*}
&\Delta_{\ell+1}(\Theta_{\ell+1};\bullet)|_{\Theta_{\ell+1}=\Theta_{\ell+1}^0}\cdot (-\Delta_{\ell+1}(\Theta_{\ell+1};\bullet)|_{\Theta_{\ell+1}=\Theta_{\ell+1}^0}) \\
&= -\left\| \Delta_{\ell+1}(\Theta_{\ell+1};\bullet)|_{\Theta_{\ell+1}=\Theta_{\ell+1}^0}\right\|^2 < 0.
\end{align*}
Since this directional derivative is negative, it indicates that iteration moving in the negative gradient direction will reduce the loss.

Using the GDA to update \(\Theta_{\ell+1}\), we define the updated parameter as
\[
\Theta_{\ell+1}^{\text{new}} = \Theta_{\ell+1}^0 - \eta \Delta_{\ell+1}(\Theta_{\ell+1};\bullet)|_{\Theta_{\ell+1}=\Theta_{\ell+1}^0},
\]
where \(\eta > 0\) is the learning rate. Given the negative directional derivative and the continuity and differentiability of the loss, there exists a sufficiently small \(\epsilon > 0\) such that for all learning rates \(\eta \in (0, \epsilon]\), the loss strictly decreases after one update step, that is,
\begin{equation}\label{Directional_Updates}
     Loss\left( \bar{u}_{\ell}^{*}(\bullet) + u_{\ell+1}(\Theta_{\ell+1}^{\text{new}}; \bullet)\right)
< Loss\left( \bar{u}_{\ell}^{*}(\bullet) + u_{\ell+1}(\Theta_{\ell+1}^{0}; \bullet)\right).
\end{equation}
This inequality demonstrates that the loss function decreases following the parameter update. By using \eqref{Directional_Updates}, we conclude that
\begin{align*}
Loss\left(\bar{u}_{\ell+1}^{*}(\bullet)\right)
    &=\min_{\Theta_{\ell+1}}Loss\left(\bar{u}_{\ell}^{*}(\bullet)+u_{\ell+1}(\Theta_{\ell+1};\bullet)\right)\\
    &\leq  Loss\left(\bar{u}_{\ell}^{*}(\bullet)+u_{\ell+1}(\Theta_{\ell+1}^{\text{new}};\bullet)\right)\\
    &< Loss\left( \bar{u}_{\ell}^{*}(\bullet) + u_{\ell+1}(\Theta_{\ell+1}^{0}; \bullet)\right)\\
    &= Loss\left(\bar{u}_{\ell}^{*}(\bullet)\right),
\end{align*} 
which confirms the desired inequality \eqref{eq:strict_decrease} and completes the proof.
\end{proof}

Proposition \ref{thm:strict_loss_decrease} establishes that the value of the loss function decreases strictly as the algorithm proceeds from one grade to the next. This result extends Theorem 4.5 of \cite{xu2023multi}, originally developed for function approximation, to the setting of initial–boundary value problems of PDEs. Moreover, when each grade consists of a single hidden ReLU layer, it was recently proved in \cite{zhang2026multigrade} that the MGDL network achieves vanishing approximation error for any continuous target function as the number of grades tends to infinity. In light of this theoretical guarantee, MGDL places deep neural networks on a footing comparable to classical approximation tools such as Fourier or wavelet bases, while retaining the flexibility of data-driven learning. These findings provide strong theoretical support for the expressive power and accuracy of the MGDL framework.

From an optimization perspective, it was shown in \cite{fang2025computational} that the enhanced stability of MGDL stems from its ability to keep the eigenvalues of the iteration matrices within the contractive interval $(-1,1)$, thereby effectively mitigating the oscillatory dynamics and learning-rate sensitivity commonly encountered in single-grade learning. As a result, MGDL permits the use of substantially larger learning rates while still ensuring optimization stability and convergence.

We emphasize that Proposition \ref{thm:strict_loss_decrease} pertains to the optimization of the empirical training loss. Nevertheless, the discrete loss serves as a numerical approximation of the continuous $L_2$ residual norm; as the number of collocation points increases, the training loss converges to the test loss. Consequently, optimizing over a progressively enriched set of collocation points leads to a genuine reduction in the global residual error.

Compared with SGL, the MGDL approach offers several advantages. First, by decomposing the original problem into a sequence of smaller subproblems, MGDL accelerates the convergence of deep network training by stabilizing the optimization process. Second, its hierarchical learning strategy mitigates overfitting and alleviates spectral bias by enabling each grade to capture relevant features before further refinement at subsequent grades. Third, the resulting representation is inherently multiscale, making it well suited for PDE solutions with complex structures. Fourth, this strategy leads to improved overall approximation accuracy. These advantages will be further demonstrated through numerical experiments.

\section{A Two-Stage Multi-Grade Model}

This section is devoted to an introduction of a two-stage multi-grade deep learning (TS-MGDL) model for training a DNN solution for PDEs.

%In the MGDL method described in the last section, 
%the first stage focuses solely on training newly added layers, which may limit the network's ability to fully capture complex patterns in the solution. 
%a present grade is to learn from the residual of the previous grade. As the number of grades increases, the frequency of the residual oscillations becomes higher, making them more challenging to capture. However, due to a relatively small number of layers in the neural network of each grade (typically less than 6 layers in our experiments to be presented later), 
%the network may struggle to capture more complex patterns contained in the underlying residual oscillations.
%Moreover, the MGDL model may be trapped in a local minimizer and miss a global minimizer. Addressing these issues requires us to enlarge the search region for a minimizer. To this end, we {\it unfreeze} some layers of the network that have been trained in the last grade and some previous grades, and retrain them all together to improve the accuracy of the resulting DNN solution. 

The MGDL method described in the previous section primarily focuses on training newly added layers. However, this approach may limit the network's ability to capture complex patterns in the solution. To address this limitation, we expand the search region for a minimizer. Specifically, we unfreeze the layers trained in the last grade and some layers in the previous grades and retrain them all together to improve the accuracy of the resulting DNN solution. This step is known as the `second stage' of training. This process is analogous to human learning, where we periodically revisit and review material previously learned from several courses or even several grades to reinforce our understanding and identify gaps in our knowledge. 

The second stage of training can be described below.
Assume that after the first stage of training, we have constructed the neural network $\bar{u}_{L}^{*}:= \sum_{i=1}^{L} u_i^{*}$, with $L$ grades, which approximates the solution of PDE. 
Note that the function $u_L^*$ has the following expression
$u_{L}^{*}:= \mathcal{N}_{k_{L}}^{*} \circ \mathcal{H}_{k_{L}-1}^{L-1*} \circ \dots \circ  \mathcal{H}_{k_1-1}^{1*},$
where $\mathcal{N}_{k_L}^*$ is trained in grade $L$ with learned parameters $\Theta_L^*$. 
In the second stage of training, we choose an integer $k>k_L$ and set the last $k$ layers of $u_L^{*}$ as trainable layers. We denote by
$
\Theta_{L,k}:=\{{W}_j^{L,k}, {b}_j^{L,k}\}_{j=1}^k
$
the parameters of the $k$ trainable layers in the second stage of training. It is important to note that $\Theta_{L,k}$ includes not only the parameters of $\mathcal{N}_{k_L}$, which are explicitly targeted for training at grade $L$, but also incorporates $k - k_L$ layers from networks previously learned in earlier grades. This indicates that $k - k_L$ layers originally frozen during prior stages are now unfrozen and included in the current training process.

We use the notation $\widetilde{u}_L:=\widetilde{u}_L(\Theta_{L, k}; \bullet)$ to denote the neural network obtained from $u_L^*$  with only the last $k$ layers made trainable and all preceding layers kept frozen. We solve the minimization problem 
\begin{equation}    \label{eq:min:grade_L:retrain}
    \min_{\Theta_{L, k}}  Loss(\bar{u}_{L-1}^{*}(\bullet)+\widetilde{u}_L(\Theta_{L, k};\bullet))
\end{equation}
using the previously trained parameters of the $k$ layers of $u_L^*$ as an initial guess, obtain the new parameters $\Theta_{L, k}^*$, and define the function $\widetilde{u}_{L}^{*}:=\widetilde{u}_L(\Theta_{L, k}^*; \bullet)$. As such, the second stage generates the approximate solution $
\hat{u}^*_L := \bar{u}_{L-1}^{*} + \widetilde{u}_{L}^{*}
$ of the PDE.
We summarize in Algorithm \ref{alg:MGL} the TS-MGDL algorithm described above for solving nonlinear PDEs.

\begin{algorithm}[htb]
\caption{Two-stage multi-grade deep learning (TS-MGDL)}
\label{alg:MGL}
\begin{algorithmic}
\REQUIRE  Number $L$ of total grades and number $k_i$ of layers for grade $i=1,2, \dots, L$. Nonlinear differential operator $\mathcal{F}$, initial condition $\mathcal{I}$, and boundary condition $\mathcal{B}$.
\ENSURE Approximate solution of the PDE.
\\ \textbf{Stage 1: Learn the network grade by grade.}
\STATE Initialize: $\mathcal{N}_{k_1}$ with parameters $\Theta_1 := \{W_i^1, b_i^1\}_{i=1}^{k_1}$.
\STATE Obtain $u_1^{*}:=\mathcal{N}_{k_1}^{*}$ by solving minimization problem \eqref{eq:min:grade_1}.
\STATE Initialize $\ell := 1$.
\WHILE{$\ell < L$}
\STATE Initialize $\mathcal{N}_{k_{\ell+1}}$ with parameters $\Theta_{\ell+1} := \{W_j^{\ell+1}, b_j^{\ell+1}\}_{j=1}^{k_{\ell+1}}$.
\STATE Set grade $\ell+1$ network $u_{\ell+1} := \mathcal{N}_{k_{\ell+1}} \circ \mathcal{H}_{k_{\ell}-1}^{\ell*} \circ \dots \circ \mathcal{H}_{k_2-1}^{2*} \circ \mathcal{H}_{k_1-1}^{1*}$.
\STATE Obtain $u_{\ell+1}^{*}$ by solving the minimization problem \eqref{eq:min:grade_L}.
\STATE $\ell \leftarrow \ell + 1$.
\ENDWHILE
\\ \textbf{Stage 2: Learn the network by layer unfreezing.}
\STATE Let $\widetilde{u}_L := u_L^{*}$ and set the $k$ last layers trainable in the network $\widetilde{u}_L$.
\STATE Obtain $\widetilde{u}_{L}^{*}$ by solving the minimization problem \eqref{eq:min:grade_L:retrain} with the initial guess $u_L^*$.
\STATE Let $\hat u_L^* := \bar{u}_{L-1}^{*} + \widetilde{u}_{L}^{*}$.
\STATE Output: Approximate solution  $\hat u_L^*$ of the PDE.
\end{algorithmic}
\end{algorithm}

Recall that we have shown in Proposition \ref{thm:loss} that the value of the loss function is reduced as we move from grade $\ell$ to grade $\ell + 1$ in the first stage.
The next proposition shows that the second stage of the TS-MGDL method can further reduce the value of the loss function.

\begin{proposition}\label{thm:loss-TS}
If $\hat{u}_{L}^{*}$ is the neural network solution of the initial-boundary value problem \eqref{eq:PDE}-\eqref{eq:PDE:boundary} obtained by the TS-MGDL method (Algorithm \ref{alg:MGL}) with $k>k_L$, 
then 
\begin{equation}\label{Justification-TS}
    Loss\left( \hat{u}_{L}^{*}(\bullet)\right) \le Loss\left(\bar{u}_{L}^{*}(\bullet)\right).
\end{equation}
\end{proposition}
\begin{proof} Without loss of generality, we assume that the parameters of the last $k$ layers of $u_L^*$ are unfrozen and retrained. We denote by 
\begin{equation}\label{Parameters'}
    \Theta_{L,k}':=\{W_j^{L,k*}, b_j^{L,k*}\}_{j=1}^k
\end{equation}
the parameters of the last $k$ layers of $u_L^*$, which have been trained in the first stage of the MGDL method of $L$ grades.
We then have that 
$
\widetilde{u}_{L}(\Theta_{L,k}';\bullet)=u_L^*(\bullet).
$
This together with the relation $\bar{u}_{L}^{*}=\bar{u}_{L-1}^{*}+u_L^*$ yields that
\begin{equation}\label{Justification1}
    Loss\left(\bar{u}_{L-1}^{*}(\bullet)+\widetilde{u}_{L}(\Theta_{L,k}';\bullet)\right) = Loss\left(\bar{u}_{L}^{*}(\bullet)\right).
\end{equation} 
Let $\Theta_{L,k}^*$ be a solution of the minimization problem \eqref{eq:min:grade_L:retrain}.  Then, for the parameters $\Theta_{L,k}'$ chosen as \eqref{Parameters'}, by employing equation \eqref{Justification1}, we obtain that 
\begin{align*}
 Loss\left(\hat{u}_{L}^{*}(\bullet)\right)&=
    Loss\left(\bar{u}_{L-1}^{*}(\bullet)+\widetilde{u}_{L}(\Theta_{L,k}^*;\bullet)\right)\\
    &=\min_{\Theta_{L,k}}Loss\left(\bar{u}_{L-1}^{*}(\bullet)+\widetilde{u}_{L}(\Theta_{L,k};\bullet)\right)\\
    &\leq  Loss\left(\bar{u}_{L-1}^{*}(\bullet)+\widetilde{u}_{L}(\Theta_{L,k}';\bullet)\right)\\
    &= Loss\left(\bar{u}_{L}^{*}(\bullet)\right),
\end{align*} 
which leads to inequality \eqref{Justification-TS}. 
\end{proof}

Inequality \eqref{Justification-TS} of Proposition \ref{thm:loss-TS} may be made strict for TS-MGDL under an additional condition on the gradient of the loss function, as Proposition \ref{thm:strict_loss_decrease} for MGDL. Specifically, we let \(\Theta_{L,k} := \{W_j^{L,k}, b_j^{L,k}\}_{j=1}^k\) denote the parameters associated with the last \(k\) layers of \(u_L\) to be retrained in the second stage. 
Let \(\Theta_{L,k}' := \{W_j^{L,k*}, b_j^{L,k*}\}_{j=1}^k\) represent the parameters of the last \(k\) layers of \(u_L^*\), which have been trained in the first stage of the MGDL method of $L$ grades. We denote by $\widetilde{\Delta}$ 
the gradient of the loss function $Loss\left( \bar{u}_{L}^{*}(\bullet) + u_{L}(\Theta_{L,k};\bullet) \right)$ with respect to $\Theta_{L,k}$. We now state the result without proof as it is similar to that for Proposition \ref{thm:strict_loss_decrease}.
% at $\Theta_{L,k}'$ is denoted as $\nabla_{\Theta_{L,k}} Loss\left( \bar{u}_{L}^{*}(\bullet) + u_{L}(\Theta_{L,k};\bullet) \right) |_{\Theta_{L,k}'}$.

%{\color{red} Analogous to Proposition \ref{thm:strict_loss_decrease}, it can be demonstrated that during the second stage of training, the loss function will strictly decrease if the parameters introduced or updated in this stage are not at a local minimum within the loss landscape. 
%We define \(\Theta_{L,k} = \{W_j^{L,k}, b_j^{L,k}\}_{j=1}^k\) as the parameters associated with the last \(k\) layers of \(u_L\). Let \(\Theta_{L,k}' = \{W_j^{L,k*}, b_j^{L,k*}\}_{j=1}^k\) represent the parameters of the last \(k\) layers of \(u_L^*\), which have been trained in the first stage of the MGDL method of $L$ grades. The gradient of the loss function $Loss\left( \bar{u}_{L}^{*}(\bullet) + u_{L}(\Theta_{L,k};\bullet) \right)$ with respect to $\Theta_{L,k}$ at $\Theta_{L,k}'$ is denoted as $\nabla_{\Theta_{L,k}} Loss\left( \bar{u}_{L}^{*}(\bullet) + u_{L}(\Theta_{L,k};\bullet) \right) |_{\Theta_{L,k}'}$.
 
\begin{proposition} \label{thm:second:strict_loss_decrease}
Assume that $\bar{u}_{L}^{*}$ and $\hat{u}_{L}^{*}$ represent neural network solutions of the initial-boundary value problem \eqref{eq:PDE}-\eqref{eq:PDE:boundary}, obtained by stages 1 and 2 of the TS-MGDL method, respectively.
Suppose that the functions \(\mathcal{F}\), \(\mathcal{B}\), and \(\mathcal{I}\) are continuously differentiable with respect to \(u(\cdot)\), and that the neural network employs continuously differentiable activation functions. If the gradient of the loss  with respect to the trainable parameter $\Theta_{L,k}$ at \(\Theta_{L,k}'\) is non-zero, that is, \(\widetilde{\Delta}|_{\Theta_{L,k}=\Theta_{L,k}'}\neq 0\),
then 
\begin{equation}\label{eq:strict:decrease:2stage}
    Loss\left( \hat{u}_{L}^{*}(\bullet)\right) < Loss\left(\bar{u}_{L}^{*}(\bullet)\right).
\end{equation}
\end{proposition}
Proposition \ref{thm:loss-TS} confirms that the value of the loss function is reduced as the TS-MGDL method moves from stage 1 to stage 2. Propositions \ref{thm:loss}, \ref{thm:strict_loss_decrease}, \ref{thm:loss-TS}, and \ref{thm:second:strict_loss_decrease} together reveal that every grade and stage can reduce the value of the loss function. This fact will be further validated by numerical experiments in the next section.

We can extend the second stage training for every grade after the first grade. We describe below the {\it generalized} TS-MGDL method. Starting grade 2, we choose an integer $\widetilde{k}_2>k_2$ and unfreeze $\widetilde{k}_2-k_2$ layers of grade 1 and retrain the parameters $\Theta_{2,\widetilde{k}_2}:=\{{W}_j^{2,\widetilde{k}_2}, {b}_j^{2,\widetilde{k}_2}\}_{j=1}^{\widetilde{k}_2}$ to obtain the optimal parameters $\Theta_{2,\widetilde{k}_2}^*$, and let $\widetilde{u}_2^*:=\widetilde{u}_2(\Theta_{2, \widetilde{k}_2}^*; \bullet)$. Repeating this process for each grade, we obtain updated solution components $\widetilde{u}_j^*$ for $j=3,4,\dots, L$. We then obtain a sequence of approximate solutions  
\begin{equation}\label{Multiscale-pre}
    \hat{u}_\ell^*:=u_1^*+\sum_{j=2}^\ell \widetilde{u}_j^*, \ \ \ell=2, 3, \dots, L, 
\end{equation}
of initial-boundary value problem  \eqref{eq:PDE}-\eqref{eq:PDE:boundary}. We have the following theorem for the generalized TS-MGDL method. The proof of this theorem is excluded for conciseness.

\begin{theorem}
If $\hat{u}^*_{\ell}$ and $\hat{u}^*_{\ell+1}$ are the neural network solutions of the initial-boundary value problem \eqref{eq:PDE}-\eqref{eq:PDE:boundary} obtained by the generalized TS-MGDL method after $\ell$ and $\ell+1$ grades, respectively, then 
\begin{equation}\label{Justification-hat}
    Loss\left( \hat{u}_{\ell+1}^{*}(\bullet)\right) \le Loss\left(\hat{u}_{\ell}^{*}(\bullet)\right)
\end{equation} 
and
\begin{equation}\label{Justification-hat2}
    Loss\left( \hat{u}_{\ell+1}^{*}(\bullet)\right) \le Loss\left(\hat{u}_{\ell}^{*}(\bullet)+u_{\ell+1}^*(\bullet)\right),
\end{equation} 
where $u_{\ell+1}^*$ denotes the solution component of grade $\ell+1$ without the second stage retraining.
\end{theorem}

Let $\widehat{\Delta}$ be
the gradient of the loss function $Loss\left( \hat{u}_{\ell}^{*}(\bullet) + u_{\ell+1}(\Theta_{{\ell+1},\widetilde{k}_{\ell+1}};\bullet) \right)$ with respect to $\Theta_{{\ell+1},\widetilde{k}_{\ell+1}}$. 
Let \(\Theta_{{\ell+1},\widetilde{k}_{\ell+1}}' := 
\{W_j^{\ell+1,\widetilde{k}_{\ell+1}*}, b_j^{{\ell+1},\widetilde{k}_{\ell+1}*}\}_{j=1}^{\widetilde{k}_{\ell+1}}\) represent the parameters of the last \(\widetilde{k}_{\ell+1}\) layers of 
\(u_{\ell+1}^*\), which have been trained in the first stage of the MGDL method of $\ell+1$ grades. The following theorem shows the exact decrease of the loss function using generalized TS-MGDL. 
The proof is omitted here due to its similarity to the reasoning presented in Proposition \ref{thm:strict_loss_decrease}.

\begin{theorem}\label{thm:GL:TSMGDL:strict_loss_decrease}
Suppose that $\hat{u}^*_{\ell}$ and $\hat{u}^*_{\ell+1}$ are the neural network solutions of the initial-boundary value problem \eqref{eq:PDE}-\eqref{eq:PDE:boundary} obtained by the generalized TS-MGDL method after $\ell$ and $\ell+1$ grades, respectively. If the gradient of the loss  with respect to the trainable parameter $\Theta_{\ell+1,\widetilde{k}_{\ell+1}}$ evaluated at \(\Theta_{\ell+1,\widetilde{k}_{\ell+1}}'\) is non-zero, i.e., $$\widehat{\Delta}|_{\Theta_{\ell+1,\widetilde{k}_{\ell+1}}=\Theta_{\ell+1,\widetilde{k}_{\ell+1}}'}\neq 0,$$
then it follows that
\begin{equation}\label{Strict_Generalized_TS_MGDL}
    Loss\left( \hat{u}_{\ell+1}^{*}(\bullet)\right) < Loss\left(\hat{u}_{\ell}^{*}(\bullet)+u_{\ell+1}^*(\bullet)\right),
\end{equation} 
where $u_{\ell+1}^*$ denotes the solution component of grade $\ell+1$ without the second stage retraining.
\end{theorem}

We devote the remainder of this section to several key implementation issues.

The number 
$k$ of trainable layers in TS-MGDL strongly influences performance. Too few trainable layers limit expressive capacity, whereas too many increase computational cost and may destabilize training. To balance expressiveness and stability, we set 
$k:=8$, selecting the layers closest to the network output 
$u_L^*$
for second-stage optimization. In our implementation, TS-MGDL employs three grades: Grade 2 typically contains 2–3 hidden layers and Grade 3 contains 3–4. With 
$k=8$, all layers from Grades 2–3 and the final layers of Grade 1 are unfrozen, enabling cross-grade adjustment of high-frequency residuals while preserving previously learned low-frequency components. This choice provides a practical rule: 
$k$ should be large enough to cover the most recent grades, yet small enough to avoid the optimization difficulties of full end-to-end training.

Careful selection of learning rates across grades and stages is essential. MGDL allows grade-specific learning rates, enabling early grades to capture coarse features quickly while later grades progressively refine fine-scale structures. 
Specifically, larger learning rates are assigned to early grades to rapidly capture coarse-scale global features, while progressively smaller rates are utilized in later grades and stages to refine fine-scale structures with high stability. For instance, in our implementation for the Burgers equation, the learning rates were set at 1e-3, 3e-4, and 2e-4 for Grades 1, 2, and 3, respectively, and further reduced to 1e-4 during the second-stage fine-tuning. This decaying schedule ensures that the global structure is learned first, followed by stable fine-tuning without risk of overshooting.

The number of epochs per grade significantly affects results. Excessive epochs in early grades risk convergence to poor local minima, while later grades require more iterations as training becomes harder. We therefore use relatively few epochs in the first stage to exploit the initial rapid loss decay, and more epochs in the second stage for accurate refinement. An adaptive strategy can further improve efficiency: early stopping based on validation loss (e.g., no improvement for 300 epochs) prevents overfitting, and dynamic epoch allocation across grades can respond to problem complexity.

Finally, TS-MGDL is compatible with advanced refinements. A multi-stage progressive unfreezing scheme can shift the trainable window toward earlier layers to refine foundational representations. Layer-wise learning rate decay, adaptive activation functions, self-adaptive loss weighting, and residual-based collocation sampling can also be incorporated to focus capacity on high-error regions. These extensions enable TS-MGDL to address highly nonlinear, multiscale PDEs while retaining the stability of greedy training and the accuracy of localized global fine-tuning.

\section{Applications to the Burgers Equation and Numerical Results}

%To evaluate the effectiveness of the proposed TS-MGDL method, we conducted numerical experiments on the 1D, 2D, and 3D Burgers equations. Numerical results demonstrate that the MGDL approach can achieve higher accuracy and faster convergence compared to the single-grade learning method. We will present our numerical findings in the next section.

In this section, we evaluate the performance of the proposed TS-MGDL method by applying it to numerical solutions of the Burgers equations of 1D, 2D, and 3D. We validate the theoretical results established in the last two sections and compare the numerical performance of the TS-MGDL method and the SGL method.

All the experiments were carried out on a 64-bit Linux server equipped with 125GB of physical memory and an Intel(R) Xeon(R) Silver 4116 CPU with a clock speed of 2.10GHz. The Nvidia Tesla V100 graphics card is utilized to accelerate training. Both the TS-MGDL method and SGL method were implemented using the DeepXDE 1.10.0 \cite{lu2021deepxde} and Tensorflow 2.8 frameworks. All computations were performed using single-precision floating-point format (float32). Furthermore, we employed the XLA (Accelerated Linear Algebra) compiler to optimize the computational graph.

In the experiments that follow, we employ $\tanh$ as the activation function for all the neural networks and solve all the optimization problems by the Adam optimizer. It is important to note that the optimization process employs full-batch training, meaning that all available collocation points are utilized in each training iteration. Unless explicitly stated otherwise (as in the robustness analysis section), a fixed random seed of 12300 is utilized across all numerical examples to ensure consistency and reproducibility.
Furthermore, regarding the Adam optimizer configuration in the numerical examples, we implemented an inverse time learning rate decay schedule defined as:
\begin{equation}
    \eta_k = \frac{\eta_0}{1 + \gamma k},
    \label{eq:inverse:time:decay}
\end{equation}
where $\eta_k$ denotes the learning rate at epoch $k$, $\eta_0$ is the initial learning rate, and $\gamma$ represents the decay rate.
The relative $L_2$ error between the approximate solution $u^{\ast}$ and the true solution $u$ is defined by
$ 
\| u - u^{\ast} \|_2/\| u \|_2,
$
where $\|v\|_2: = (\sum_{i=1}^N v(\mathbf{x_i})^2)^{1/2}$, with $\mathbf{x_i}$ being test points.

\subsection{The 1D Burgers equation}
In this example, we consider the 1D Burgers equation, which takes the form
\begin{equation}\label{eq:Burgers:1D}
	u_t (t,x)+  u(t,x) u_x(t,x) - (0.01/\pi)  u_{xx} (t,x)= 0, \quad t \in (0,1], \ x\in (-1, 1),
\end{equation}
where $u(t,x)$ represents the velocity of the fluid at time $t$ and position $x$, $u_x(t,x)$ and $u_{xx}(t,x)$ denote its first and second spatial derivatives, respectively. The initial condition for this equation is given by $u(0,x) = -\sin(\pi x)$,
which represents a sinusoidal velocity profile at time $t=0$. The boundary conditions are given by
$u(t,-1) = u(t, 1) = 0$, which specify that the velocity of the fluid is zero on the boundaries of the domain, $x=-1$ and $x=1$, for all times $t$.
 
\begin{figure}[htbp]
  \centering
  \includegraphics[width=0.4\textwidth]{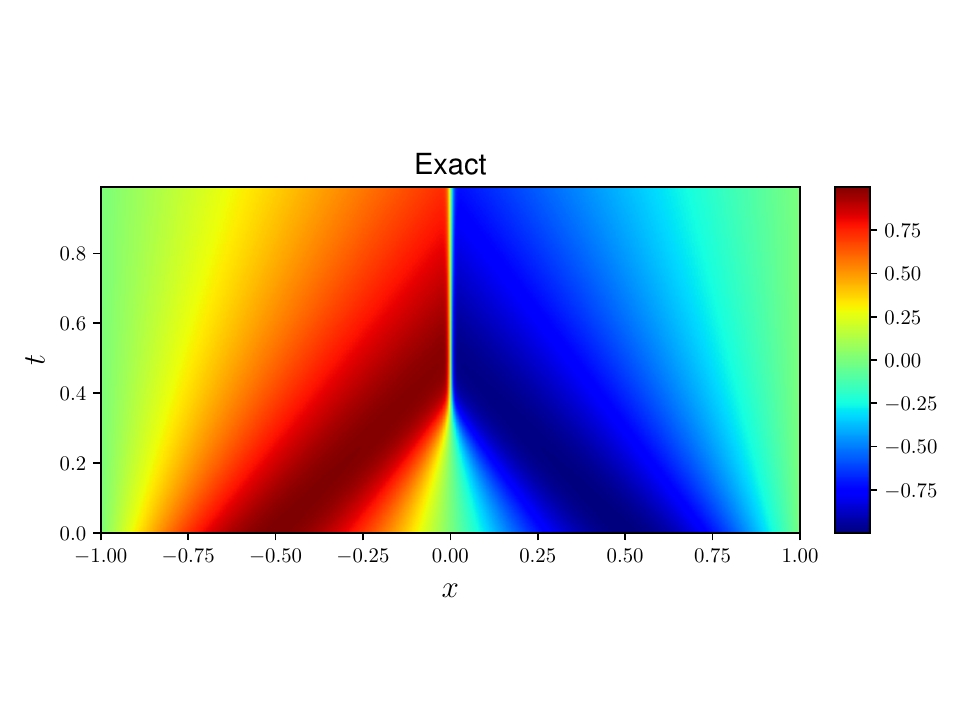}
  \caption{Exact solution of 1D Burgers' equation.}
  \label{fig:1d:exact}
\end{figure}

According to \cite{BASDEVANT1986}, the analytic solution of the initial-boundary value problem of the 1D Burgers equation takes the form  
$$
u(t, x) = -\frac{\int_{-\infty}^{+\infty} \sin \pi(x-\eta) h(x-\eta) \exp \left(-\eta^{2} / 4 \nu t\right) \mathrm{d} \eta }
{\int_{-\infty}^{+\infty} h(x-\eta) \exp \left(-\eta^{2} / 4 \nu t\right) \mathrm{d} \eta},\ \ t \in [0,1], \ x\in [-1, 1],
$$
where $\nu := 0.01/\pi$ and $h(y) := \exp (-\cos \pi y / 2 \pi \nu)$.
The analytic solution, shown in Figure \ref{fig:1d:exact}, which exhibits a big change at $x=1$ when $t>0.4$, will be used for comparison of the TS-MGDL and the SGL methods.

We employed the Hammersley sampling method to randomly generate training sample points. Along the boundaries, a total of $N_b := 80$ random points were generated. Furthermore, $N_0 := 120$ random points were generated based on the initial conditions. In the interior region $(0,1] \times (-1, 1)$, there are $ N_f := 10,000$ random points generated. The training points are resampled every 5,000 epochs during the training process. The test set consists of evenly spaced grid points obtained by dividing the spatiotemporal region $[0,1] \times [-1, 1]$. The total number of test points is $100 \times 256$.

Table \ref{table:1D:Burgers:SGL:structure} presents three different network structures (SGL-1, SGL-2, and SGL-3) of the SGL method for solving the 1D Burgers equation. They all consist of 2 input neurons and 1 output neuron. Each of the network structures is represented by a list of integers, which indicate the number of neurons in the corresponding layers of the network.

Table \ref{table:1D:Burgers:MGL:structure} displays the network structures  (Grade 1, Grade 2, and Grade 3) of the TS-MGDL model for solving the 1D Burgers equation. Each neural network has 2 input neurons and 1 output neuron. To facilitate comparison, the network structures of SGL-1 and Grade 1 of the TS-MGDL model are set to be the same. Likewise, the network structures of SGL-2 and Grade 2 of the TS-MGDL model are set to be the same, and the network structures of SGL-3 and Grade 3 of the TS-MGDL model are also set to be the same. This ensures that any differences in performance between these methods can be attributed to the different learning strategies employed, rather than differences in network architecture.
An asterisk (*) next to a number in the list indicates that the corresponding weight parameters are fixed to be those trained in the previous grades during training in the present grade of the first stage. Only those layers without an asterisk will be trained in the current grade. In the second stage of training, the last 8 hidden layers of the Grade 3 network will be unfrozen, which are located near the output.  This means that the weights of these layers will be allowed to change during training. 

\begin{table}[htbp]
		\centering
\begin{tabular}{l|l}
 \hline
   Methods & Network structure  \\
   \hline
   SGL-1 &   [2, 128, 128, 128, 128, 128, 128, 1 ]     \\
   SGL-2 &  [2, 128, 128, 128, 128, 128, 128, 256, 256, 1]     \\
   SGL-3 &   [2, 128, 128, 128, 128, 128, 128, 256, 256, 256, 256, 256, 128, 1]  \\
   \hline
\end{tabular}
  \caption{Network structure of single-grade learning for the 1D Burgers equation.}
 \label{table:1D:Burgers:SGL:structure}
\end{table}

\begin{table}[htpb]
		\centering
\begin{tabular}{c|l}
 \hline
   Grade & Network structure  \\
   \hline
    1 &   [2, 128, 128, 128, 128, 128, 128, 1 ]     \\
    2 &  [2, 128*, 128*, 128*, 128*, 128*, 128*, 256, 256, 1]     \\
    3 &   [2, 128*, 128*, 128*, 128*, 128*, 128*, 256*, 256*, 256, 256, 256, 128, 1]  \\
   \hline
\end{tabular}
  \caption{Network structure of multi-grade learning for the 1D Burgers equation.}
 \label{table:1D:Burgers:MGL:structure}
\end{table}

\subsubsection{Performance comparison under standard optimization}

In this section, we compare the performance of the proposed TS-MGDL method with standard single-grade learning (SGL) models using a purely gradient-based optimization approach. Specifically, both the TS-MGDL and the SGL baseline models (SGL-1, SGL-2, and SGL-3) are optimized using the Adam optimizer for all epochs. To provide a detailed analysis of the training dynamics and error profiles, the results presented in this subsection are obtained from a representative experiment initialized with a fixed random seed of 12300. It is important to note that the data reported here correspond to a single experimental trial; a comprehensive robustness analysis aggregating metrics across multiple independent runs is presented in Section \ref{sec:1d:robustness}. Furthermore, unless explicitly stated otherwise, the learning rate schedule for the Adam optimizer follows the inverse time decay scheme defined in Equation \eqref{eq:inverse:time:decay}.

Table \ref{table:1D:Burgers:MG} presents the numerical results of the TS-MGDL model for solving the 1D Burgers equation. The training process is divided into two stages. 
In the first stage of training, we use a learning rate of 1e-3 for Grade 1. As we move to Grades 2 and 3, the learning rate is reduced while the number of epochs is increased to 100,000. This helps balance the trade-off between learning complex representations and preventing overfitting. In the second stage of training, we use a learning rate of 3e-4 and 300,000 epochs to fine-tune the network's learned representations.

As we can see from Table \ref{table:1D:Burgers:MG}, the relative $L_2$ error decreases as we move from lower to higher grades within the first stage of training. Specifically, in Stage 1, Grade 1, we obtain a relative $L_2$ error of 7.45e-04, which decreases to 3.30e-04 in Grade 2 and further to 1.15e-04 in Grade 3. This demonstrates that increasing the number of grades can improve performance. In Stage 2, the relative $L_2$ error is further reduced to 9.65e-06 through fine-tuning of the Grade 3 network. 
 
\begin{table}[htpb]
		\centering
\begin{tabular}{c|c|c|c|c}
 \hline
   Stage \& grade &  Learning rate & Decay rate & Epochs & Relative $L_2$ error  \\
   \hline
   Stage 1, Grade 1 & 1e-3  & 1e-4   & 50,000  &  7.45e-04\\
   Stage 1, Grade 2 & 3e-4  & 1e-4   & 100,000   & 3.30e-04\\   
   Stage 1, Grade 3 & 2e-4 & 1e-4  & 100,000  & 1.15e-04\\
   Stage 2          & 3e-4 & 1e-4  & 300,000  & \textbf{9.65e-06}\\
   \hline
\end{tabular}
  \caption{Numerical results of the TS-MGDL method for the 1D Burgers equation.}
 \label{table:1D:Burgers:MG}
\end{table}

\begin{figure}[htbp]
  \centering
  \includegraphics[width=0.5\textwidth]{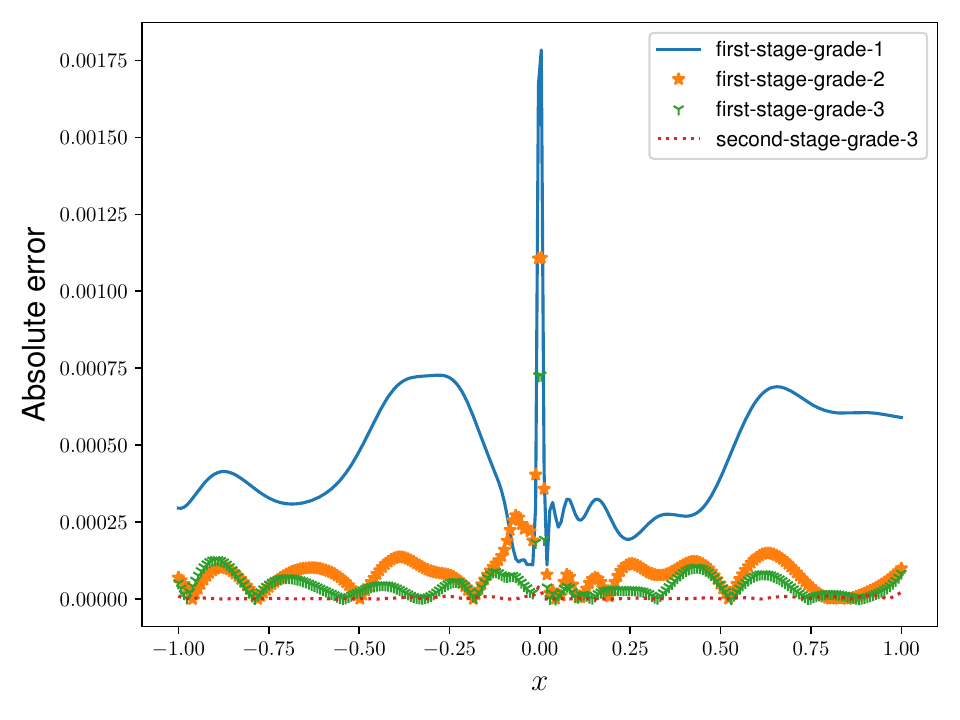}
  \caption{Absolute errors of TS-MGDL for the 1D Burgers equation at $t=0.99$.}
  \label{fig:1d:ts-mgdl-error}
\end{figure}

Table \ref{table:1D:Burgers} compares the performance of single-grade learning models (SGL-1, SGL-2, and SGL-3) and the TS-MGDL model for solving the 1D Burgers equation. The hyperparameters used for training the SGL neural networks include a learning rate of 1e-3, a decay rate of 1e-4, and 550,000 epochs for each model. The relative $L_2$ error for SGL-1 is 5.75e-04, for SGL-2 is 3.02e-04, and for SGL-3 is 2.53e-04. 
In comparing the numerical results of the TS-MGDL and SGL methods, we can see that TS-MGDL achieves lower relative $L_2$ errors than SGL. The TS-MGDL method achieves a relative $L_2$ error of 9.65e-06, which is significantly lower than the lowest error achieved by any of the SGL methods, which is 2.53e-04. Moreover, Grades 2 and 3 of Stage 1 generate approximate solutions having accuracy comparable to those generated by SGL-2 and SGL-3, but require significantly fewer epochs.
All these facts suggest that TS-MGDL is more effective than the SGL method for problems with complex solutions, as it approximates the solution of the PDE by gradually increasing the complexity of the neural network.

\begin{table}[htpb]
		\centering
\begin{tabular}{c|c|c|c|c}
 \hline
   Methods &  Learning rate &  Decay rate & Epochs & Relative $L_2$ error  \\
   \hline
   SGL-1 & 1e-3  & 1e-4   & 550,000    & 5.75e-04  \\
   SGL-2 & 1e-3  & 1e-4     & 550,000    & 3.02e-04   \\
   SGL-3 & 1e-3  & 1e-4    & 550,000   & 2.53e-04  \\
   TS-MGDL  & -- & -- &550,000   & \textbf{9.65e-06}
   \\
   \hline
\end{tabular}
  \caption{Numerical results for the 1D Burgers equation.}
 \label{table:1D:Burgers}
\end{table}

%
%an improved approximation of the residual compared to $u_3^*$. 
%Moreover, Figure \ref{fig:ex1-stage2} presents the predicted solution of stage 2, which is %obtained by combining $u_1^*$, $u_2^*$, and $\widetilde{u}_{3}^{*}$.
Figure \ref{fig:1d:ts-mgdl-error} provides a comparison of the absolute errors between the predicted and true solutions at $t=0.99$, for Grades 1, 2, 3, and Stage 2.
The results indicate that the TS-MGDL method approximates the solution well, with the error gradually decreasing from grade 1 to grade 3. After the training in the second stage, the error in the solution becomes significantly smaller than that of the first stage, indicating the effectiveness of the TS-MGDL method in improving the accuracy of the solution.

\begin{figure}[htbp]
  \centering
  \subfloat[SGL-1]{
    \includegraphics[width=0.45\textwidth]{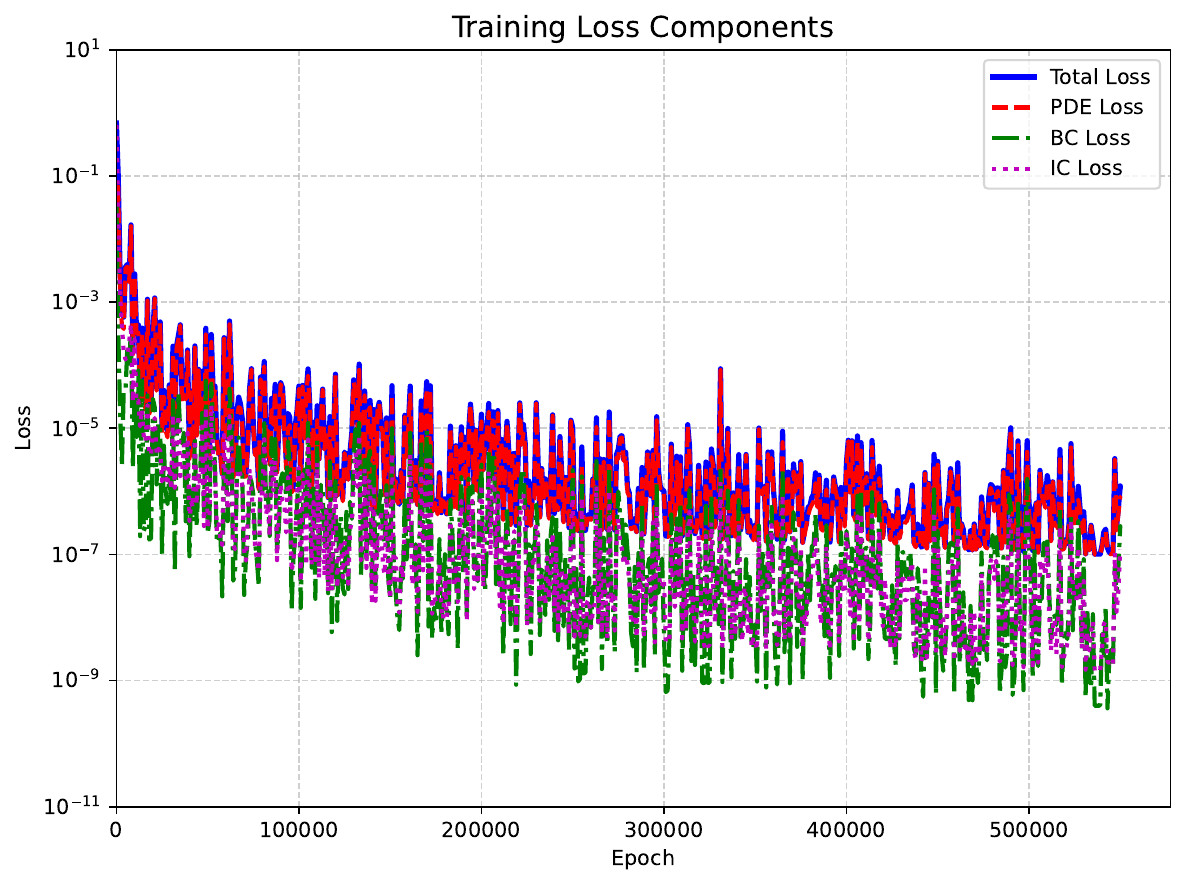}
    \label{fig:loss:sgl-1:1d}
  }
  \subfloat[SGL-2]{
    \includegraphics[width=0.45\textwidth]{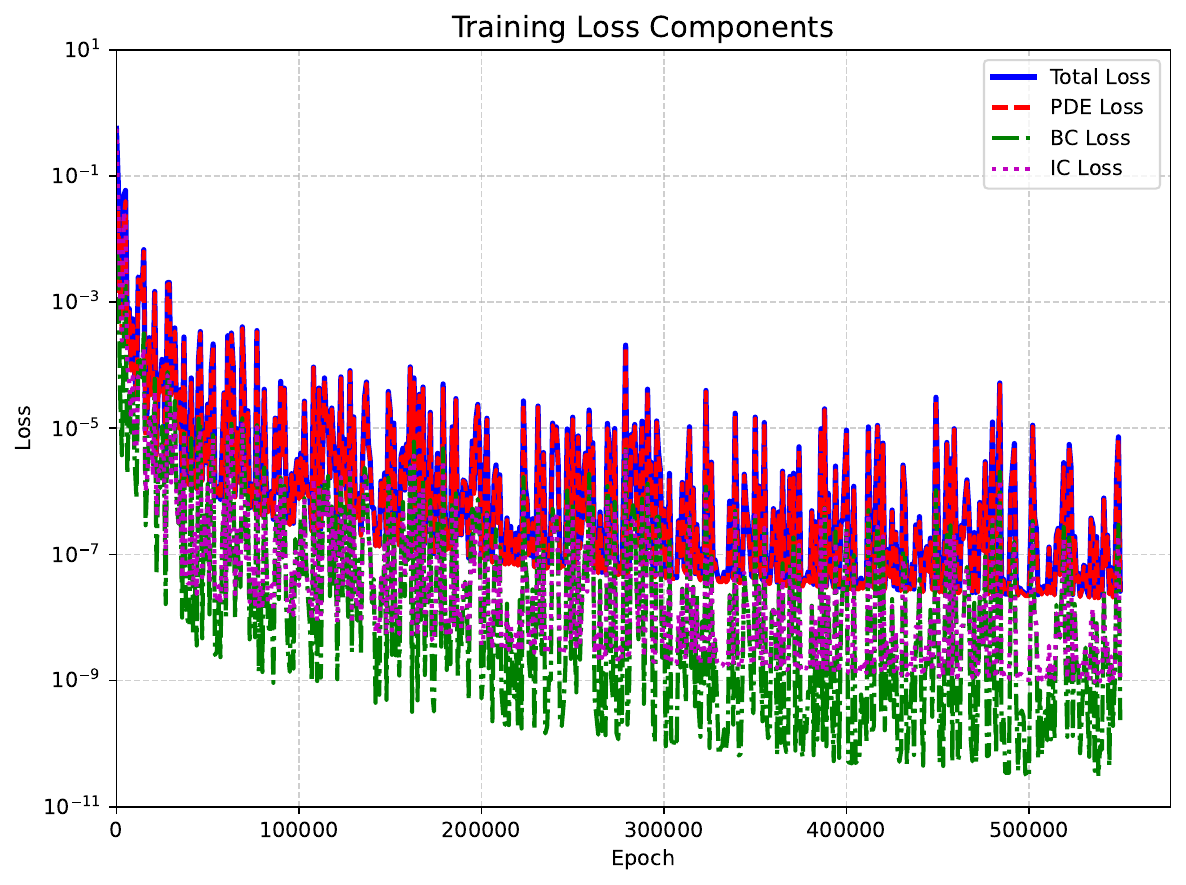}
    \label{fig:loss:sgl-2:1d}
  }\\
    \subfloat[SGL-3]{
    \includegraphics[width=0.45\textwidth]{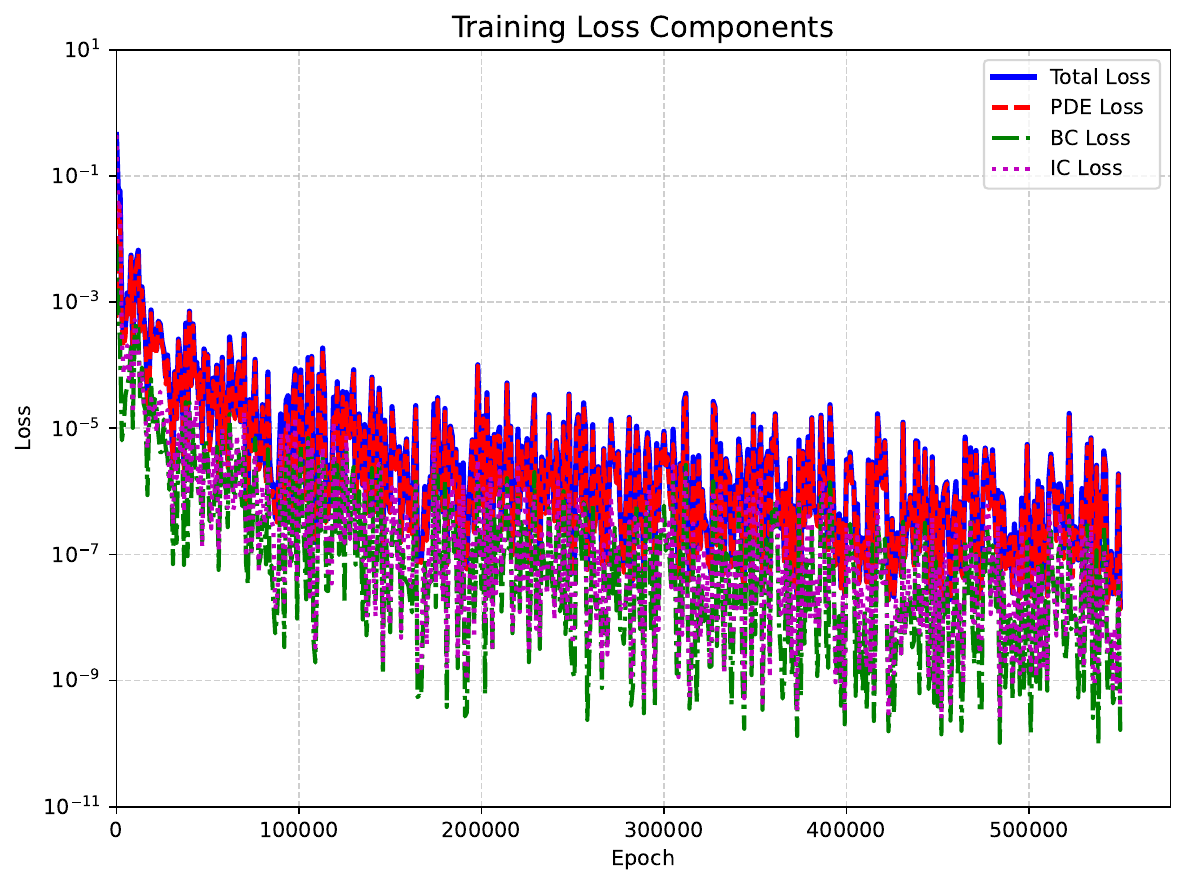}
    \label{fig:loss:sgl-3:1d}
    }
   \subfloat[TS-MGDL]{
    \includegraphics[width=0.45\textwidth]{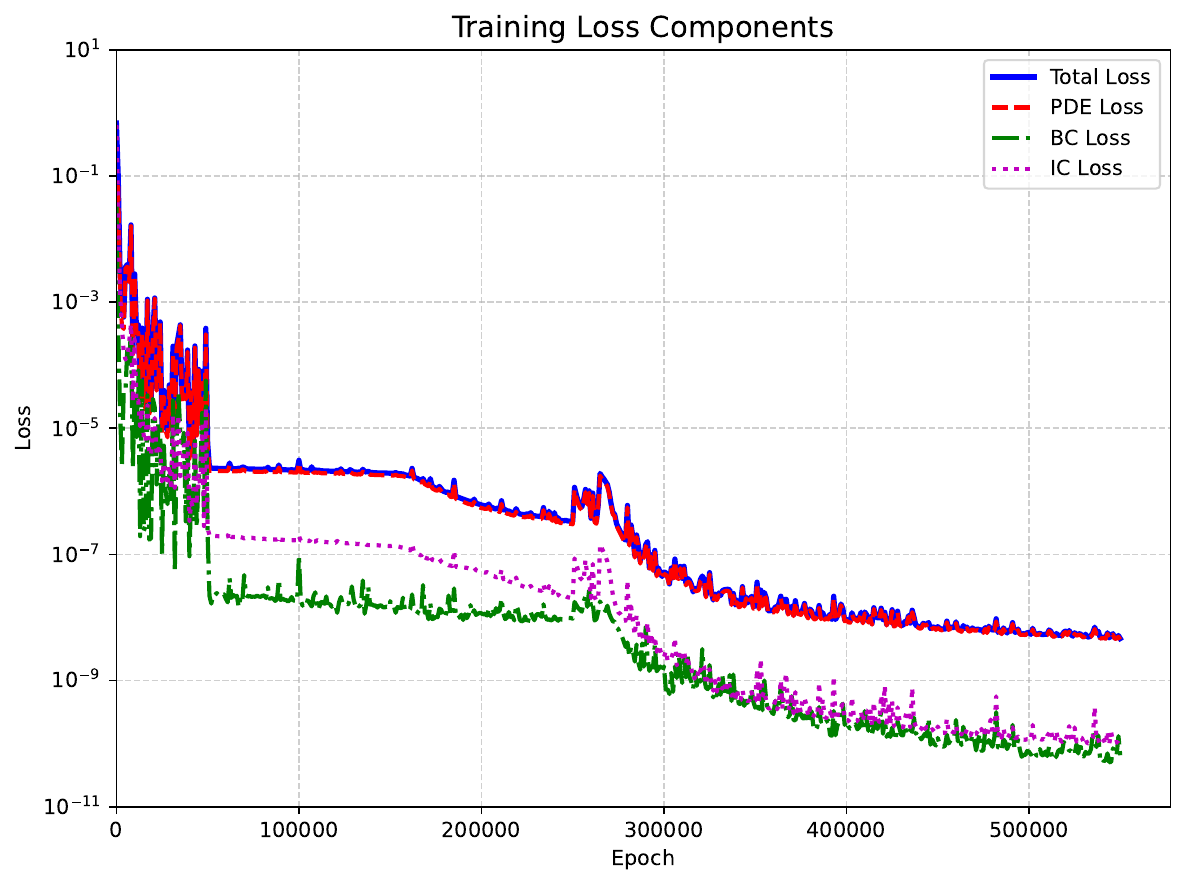}
    \label{fig:loss:mgl:1d}
  }
  \caption{Training loss for SGL-1, SGL-2, SGL-3, and TS-MGDL for the 1D Burgers equation.}
  \label{fig:loss:1d}
\end{figure}

\begin{figure}[htbp]
  \centering
  \subfloat[SGL-1]{
    \includegraphics[width=0.45\textwidth]{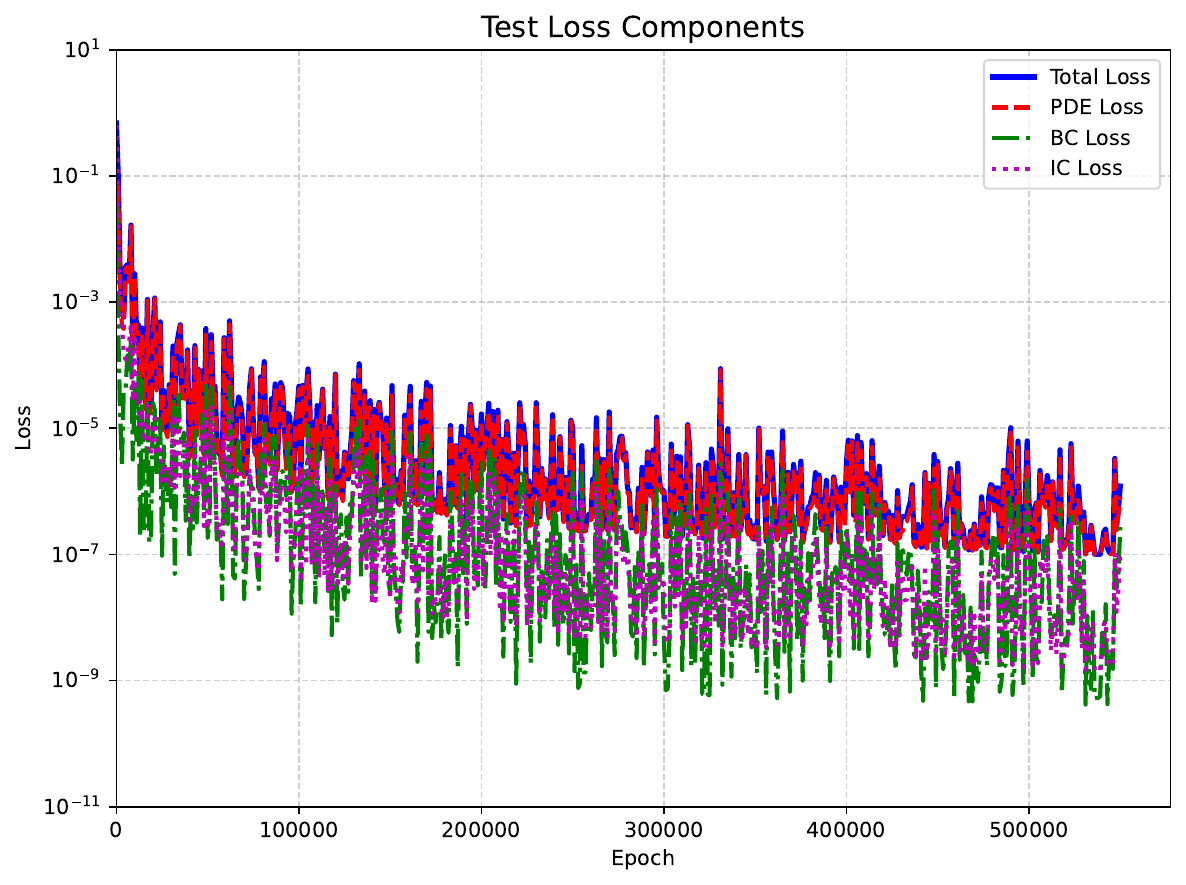}
    \label{fig:testloss:sgl-1:1d}
  }
  \subfloat[SGL-2]{
    \includegraphics[width=0.45\textwidth]{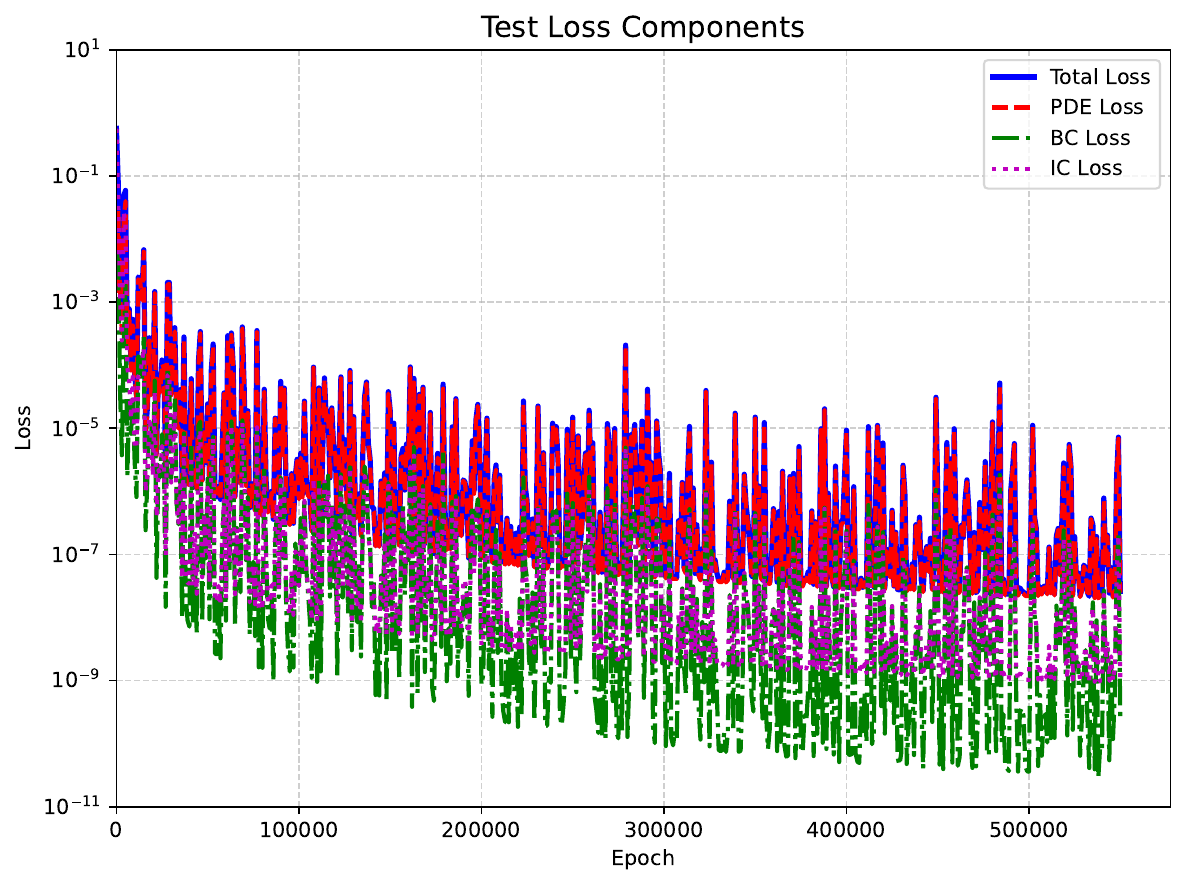}
    \label{fig:testloss:sgl-2:1d}
  }\\
    \subfloat[SGL-3]{
    \includegraphics[width=0.45\textwidth]{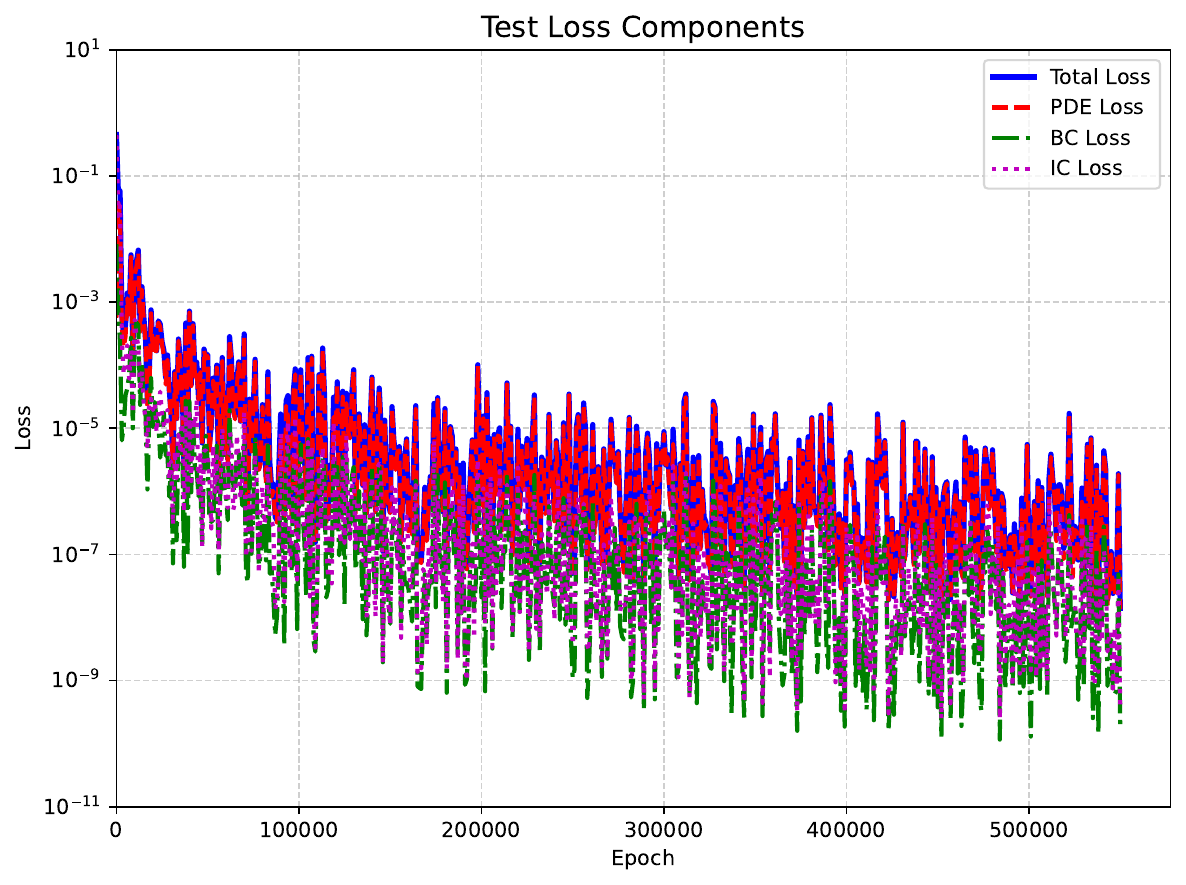}
    \label{fig:testloss:sgl-3:1d}
    }
   \subfloat[TS-MGDL]{
    \includegraphics[width=0.45\textwidth]{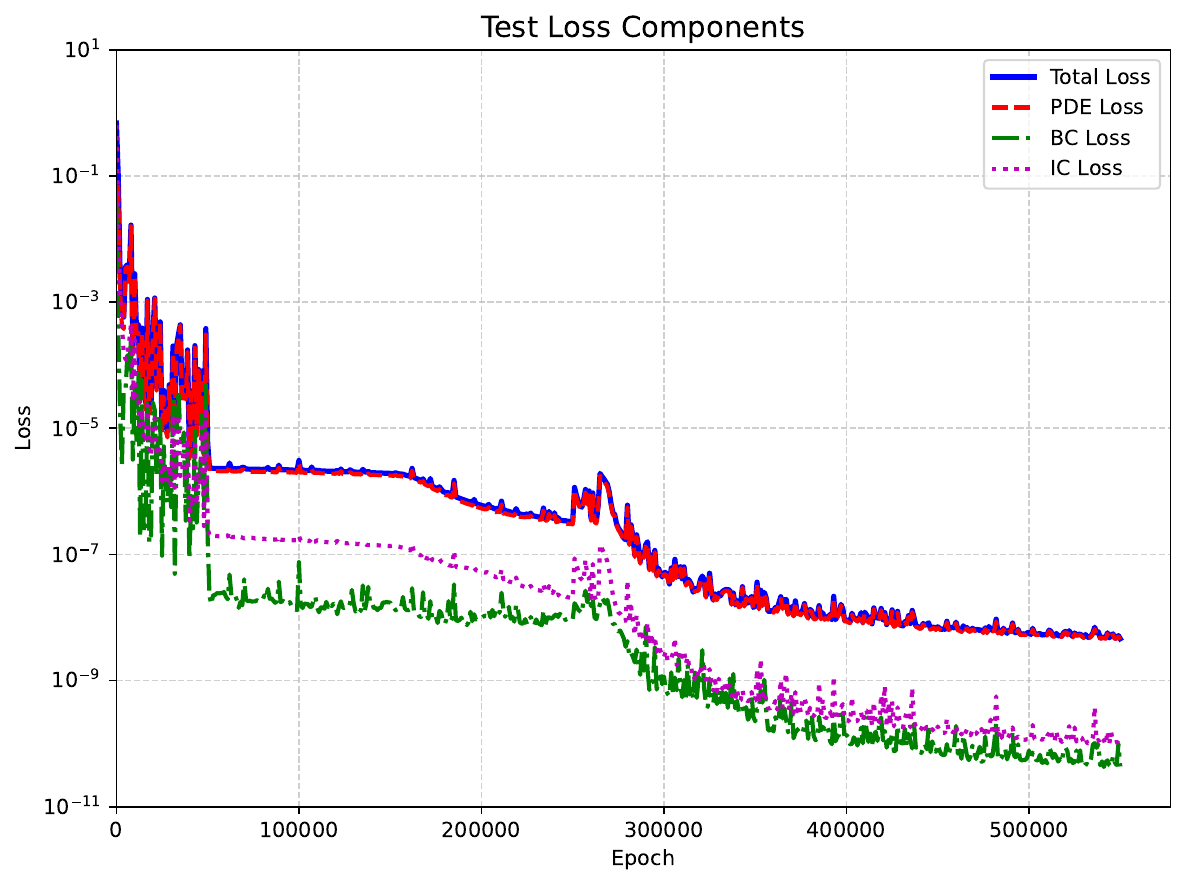}
    \label{fig:testloss:mgl:1d}
  }
  \caption{Test loss for SGL-1, SGL-2, SGL-3, and TS-MGDL for the 1D Burgers equation.}
  \label{fig:testloss:1d}
\end{figure}

Figure \ref{fig:loss:1d} shows the training loss for SGL-1, SGL-2, SGL-3, and TS-MGDL. To further dissect the training dynamics, we examine the individual components of the total loss: the PDE loss ($Loss_{PDE}$), which enforces the governing equation; the BC loss $Loss_{B}$, which imposes boundary conditions; and the IC loss ($Loss_{I}$), representing initial condition constraints. It can be observed that the SGL methods exhibit significant oscillations in the loss curve during the later stages of training, leading to slow convergence and a delayed reduction in loss. The loss curve for the TS-MGDL method exhibits a smooth descent without pronounced oscillations, indicating stable convergence. Moreover, the TS-MGDL method demonstrates a faster rate of loss reduction. In particular, Figure \ref{fig:loss:1d}\subref{fig:loss:mgl:1d} confirms the theoretical results stated in both Propositions \ref{thm:loss} and \ref{thm:loss-TS} for the 1D Burgers equation.
These findings highlight the superiority of the TS-MGDL approach in terms of both accuracy and efficiency when approximating solutions for the 1D Burgers equation. 
Figure \ref{fig:testloss:1d} illustrates the test loss for SGL-1, SGL-2, SGL-3, and TS-MGDL, which exhibits trends similar to the training loss profiles. This consistency indicates that given a sufficiently large number of training points, the training loss serves as a reliable indicator of the model's generalization performance, resulting in synchronized behavior between training and testing metrics.

\begin{figure}[htpb]
  \centering
  
  %\subfloat[Predicted solution of SGL-1]{
  %  \includegraphics[width=0.45\textwidth]{1d-sgl-1-predict.pdf}
  %  \label{fig:ex1-sgl-1}
  %}
  \subfloat[Absolute error of SGL-1]{
    \includegraphics[width=0.3\textwidth]{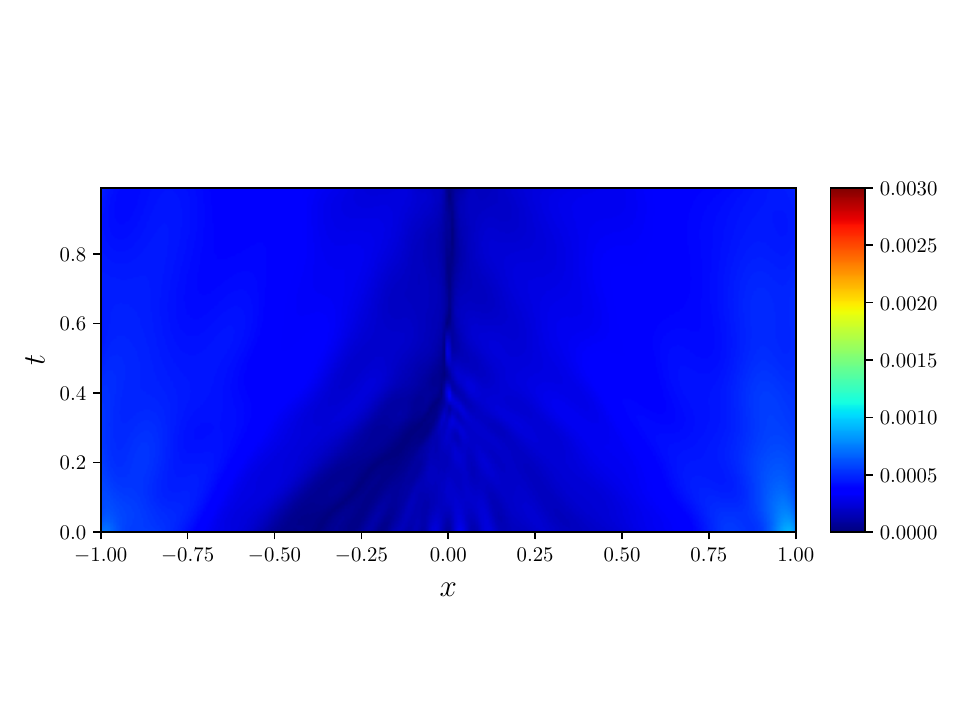}
    \label{fig:ex1-sgl-1:error}
  }
  %  \subfloat[Predicted solution of SGL-2]{
  %  \includegraphics[width=0.45\textwidth]{1d-sgl-2-predict.pdf}
  %  \label{fig:ex1-sgl-2}
  %}\\
  \subfloat[Absolute error of SGL-2]{
    \includegraphics[width=0.3\textwidth]{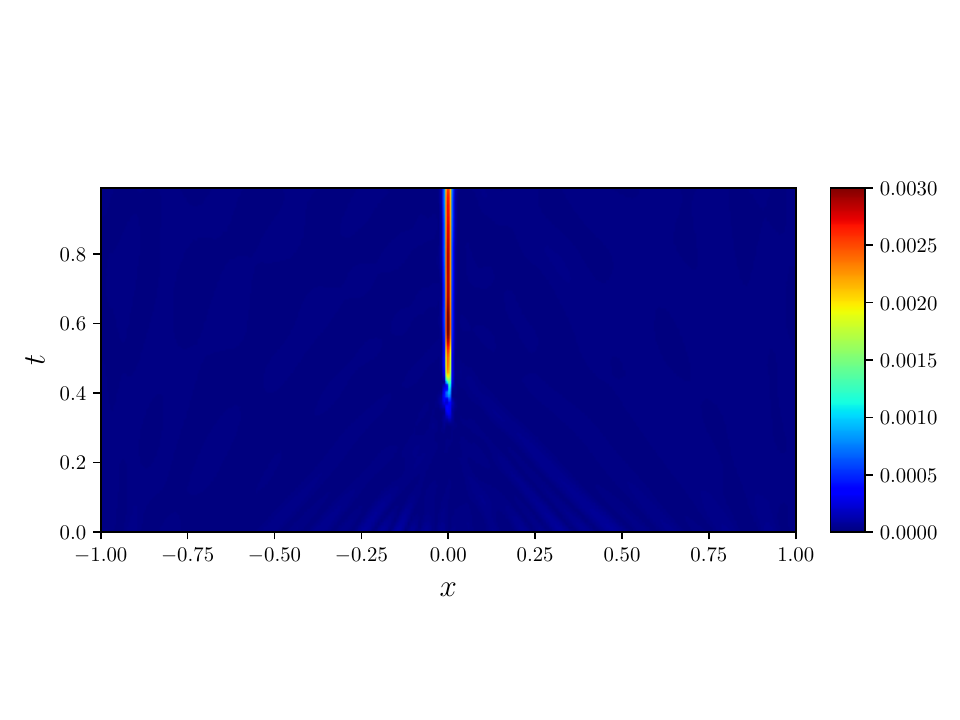}
    \label{fig:ex1-sgl-2:error}
  }\\
  %  \subfloat[Predicted solution of SGL-3]{
  %  \includegraphics[width=0.45\textwidth]{1d-sgl-3-predict.pdf}
  %  \label{fig:ex1-sgl-3}
  %}
  \subfloat[Absolute error of SGL-3]{
    \includegraphics[width=0.3\textwidth]{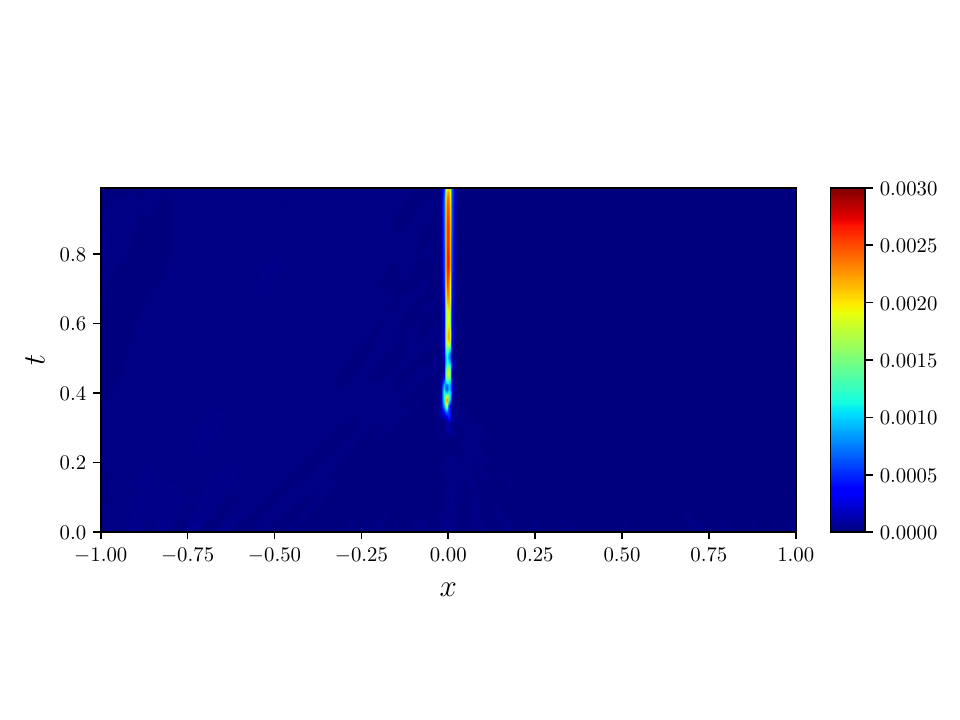}
    \label{fig:ex1-sgl-3:error}
  }
  %  \subfloat[Predicted solution of TS-MGDL]{
  %  \includegraphics[width=0.45\textwidth]{mgl-1d-predict.pdf}
  %  \label{fig:ex1-mgl}
  %}
  \subfloat[Absolute error of TS-MGDL]{
    \includegraphics[width=0.3\textwidth]{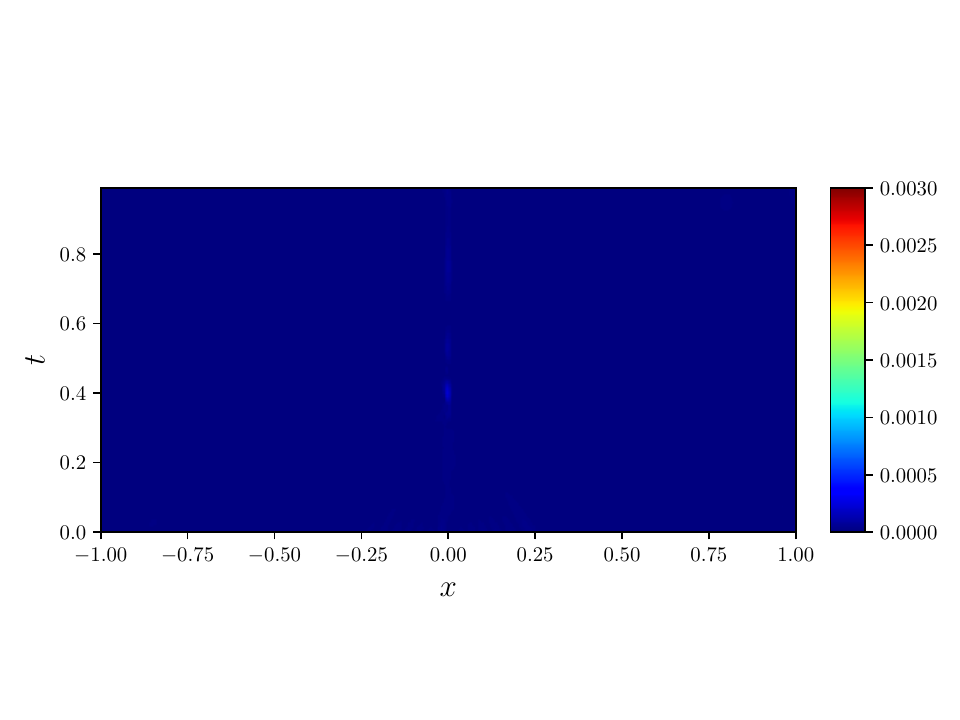}
    \label{fig:ex1-mgl:error}
  }

  \caption{Comparison of SGL1, SGL2, SGL3 with TS-MGDL  for the 1D Burgers equation.}
  \label{fig:1d}
\end{figure}

\begin{figure}[htbp]
    \centering
    \includegraphics[width=0.5\textwidth]{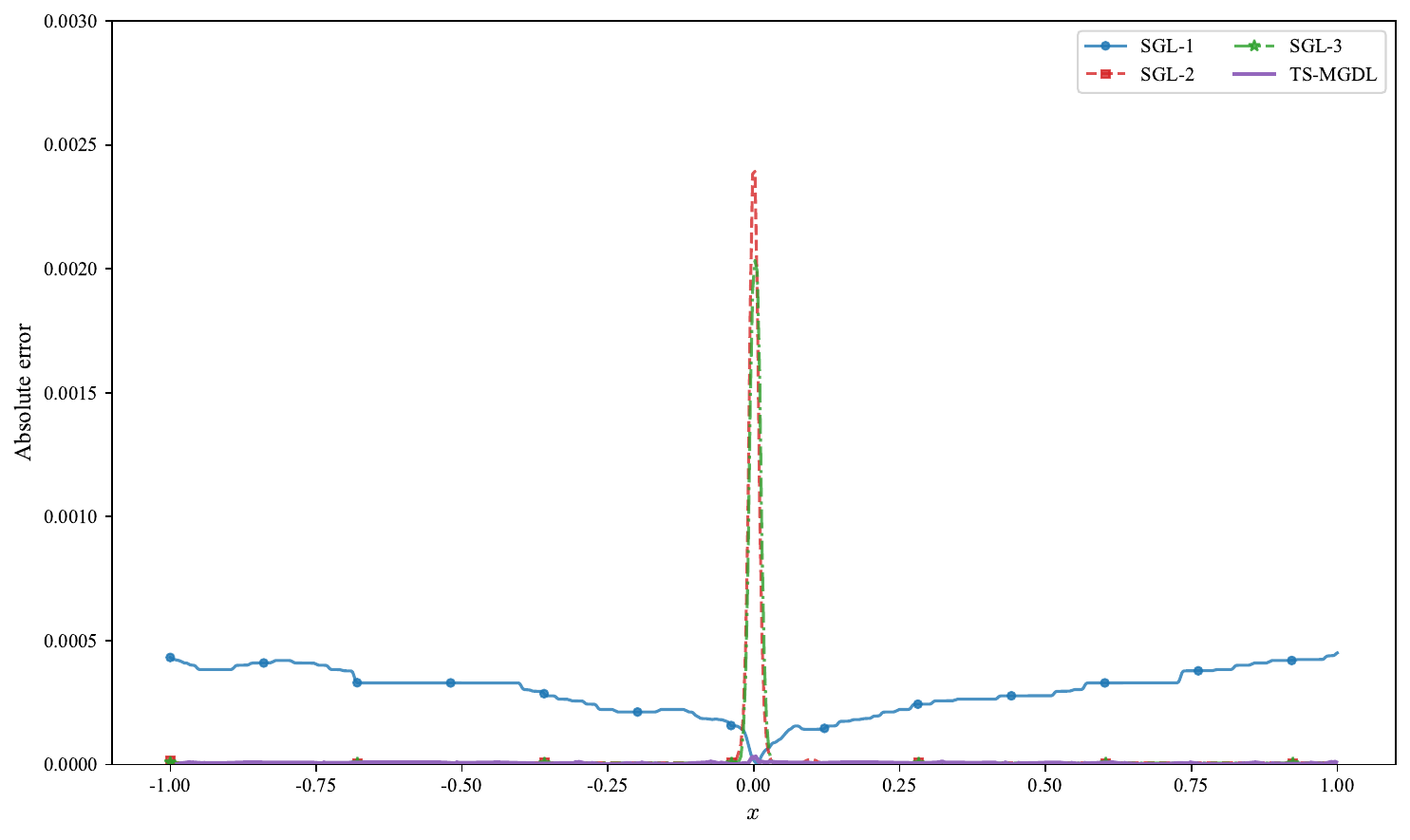} 
    \caption{Distribution of absolute errors for the 1D Burgers equation at $t=0.99$ using different numerical schemes: SGL-1, SGL-2, SGL-3, and the proposed TS-MGDL method.}
    \label{fig:error_comparison_t099}
\end{figure}

\begin{figure}[htbp]
    \centering
    \includegraphics[width=0.8\textwidth]{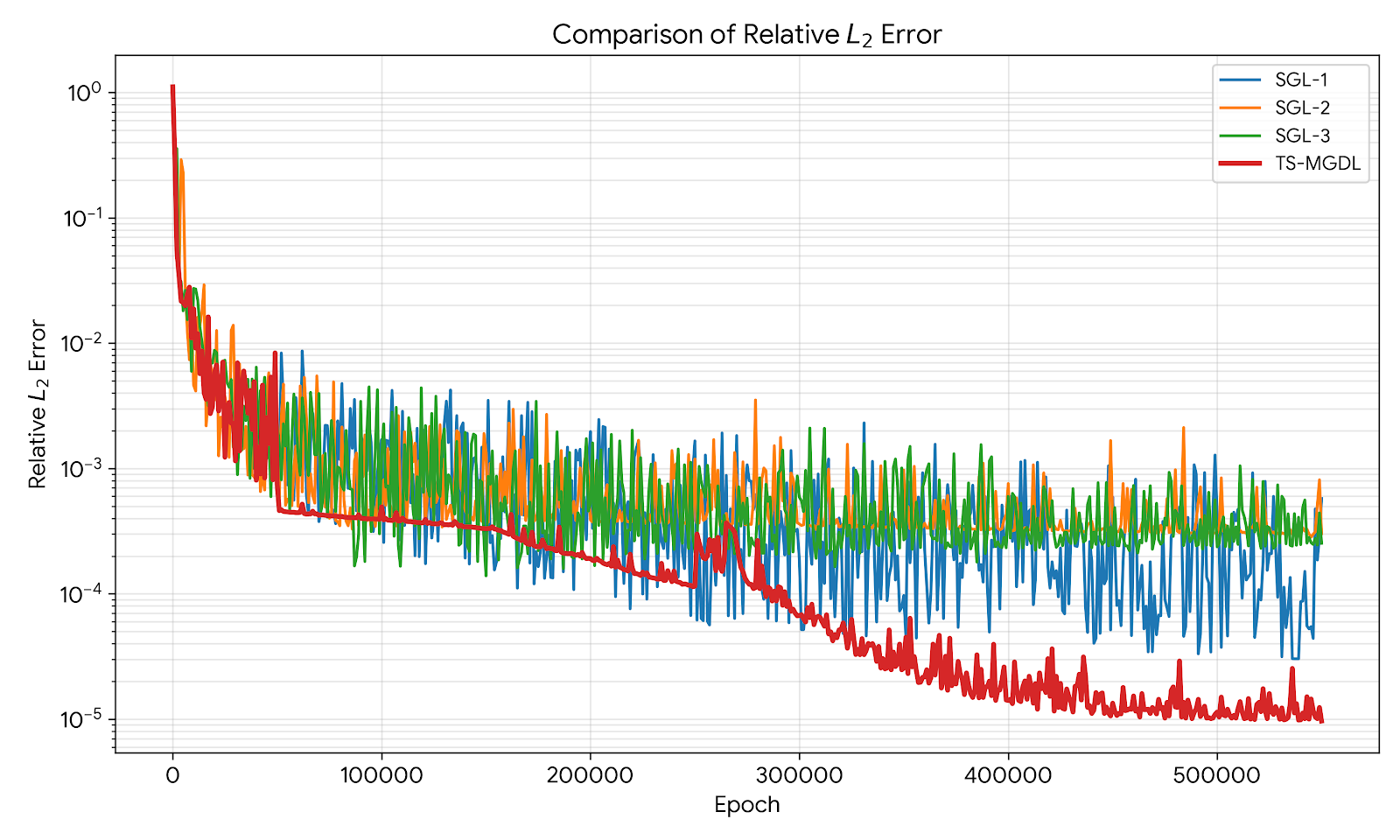}
    \caption{Relative $L_2$ error comparison between TS-MGDL and SGL methods for the 1D Burgers equation.}
    \label{fig:1d-burgers}
\end{figure}

The figures presented in Figure \ref{fig:1d} depict the predicted solutions and absolute errors for SGL-1, SGL-2, SGL-3, and TS-MGDL, for solving the 1D Burgers equation. Compared to the SGL methods, it is evident that the TS-MGDL method achieves a significantly smaller absolute error. Specifically, at the time instant $t=0.99$, Figure \ref{fig:error_comparison_t099} reveals that the absolute errors of the SGL-type schemes remain considerably high, on the order of $10^{-3}$. It clearly shows that the proposed TS-MGDL method achieves an error magnitude approximately two orders of magnitude lower than the SGL configurations, further validating its exceptional numerical precision. Further more, Figure \ref{fig:1d-burgers} compares the convergence of the relative $L_2$ error between TS-MGDL and SGL methods. While the SGL baselines fluctuate and plateau between $10^{-3}$ and $10^{-4}$, TS-MGDL demonstrates a superior training trajectory, ultimately achieving a high-precision error of 9.65e-06. This order-of-magnitude improvement signifies that the multi-grade learning strategy effectively overcomes the optimization bottlenecks encountered by standard single-grade networks.

To provide a more comprehensive evaluation, we further investigated the performance of single-grade learning models under various learning rate scheduling strategies, denoted as SGL-Piecewise and SGL-Exp. Specifically, SGL-Piecewise implements a custom multi-stage scheduler meticulously designed to mimic the learning rate evolution inherently present in the TS-MGDL framework, while SGL-Exp employs a conventional smooth optimization trajectory with a staircase exponential decay.

Specifically, the SGL-Piecewise scheduler was meticulously designed to mimic the learning rate evolution inherently present in the TS-MGDL framework. By segmenting the training of a single-grade network into discrete phases that correspond to the grade transitions in our multi-grade approach, we can isolate whether the advantages of TS-MGDL stem from the scheduling itself or from its hierarchical architectural refinement. The piecewise inverse time decay scheduler is defined as:
\begin{equation}
\eta(t) = \frac{\eta_{base}(t)}{1 + \gamma \cdot (t - t_{start}(t))},
\label{eq:piecewise_lr}
\end{equation}
where $t$ is the current training step, $\gamma = 10^{-4}$ is the local decay rate, and the base learning rate $\eta_{base}(t)$ and stage starting step $t_{start}(t)$ are piecewise constant functions:
\begin{equation}
(\eta_{base}(t), t_{start}(t)) = 
\begin{cases} 
(0.001, 0), & t < 50,000, \\
(0.0003, 50,000), & 50,000 \le t < 150,000, \\
(0.0002, 150,000), & 150,000 \le t < 250,000, \\
(0.0003, 250,000), & t \ge 250,000.
\end{cases}
\end{equation}

For the SGL-Exp models, the learning rate was initialized at $0.001$ and adjusted every $1,000$ steps using a staircase decay with a rate of $0.96$, representing a standard exponential decay optimization trajectory. These two variants ensure that the baseline SGL models are optimized under both conventional smooth decay and the multi-stage decay conditions that mirror our proposed method TS-MGDL. This comprehensive comparison, summarized in Table \ref{table:1D:Burgers:Schedules}, allows for a rigorous assessment of the performance gains attributable to the TS-MGDL architectural design rather than mere scheduling effects.

The results show that advanced scheduling strategies provide improvements primarily for shallower networks, with SGL-2-Exp and SGL-1-Piecewise achieving errors of 7.42e-05 and 1.23e-04, respectively. In contrast, the performance of deeper single-grade networks remains suboptimal or even degrades, as exemplified by the 3.58e-04 error of SGL-3-Piecewise. This finding shows that simply replicating the learning rate schedule of TS-MGDL is insufficient to address the optimization bottlenecks inherent in deep, static architectures.

Crucially, the TS-MGDL method significantly outperforms all SGL variants across all tested schedules. Its relative $L_2$ error (9.65e-06) is at least one order of magnitude lower than the best results achieved by any SGL configuration. This comparison further confirms that the superiority of the TS-MGDL framework is rooted in its progressive architectural refinement and multi-grade residual learning strategy, rather than mere hyperparameter configurations.

\begin{table}[htpb]
    \centering
    \caption{Comparison of SGL models under different learning rate schedules for the 1D Burgers equation. SGL-Exp represents smooth exponential decay, while SGL-Piecewise mimics the staged learning rate transitions of TS-MGDL.}
    \label{table:1D:Burgers:Schedules}
    \begin{tabular}{l|l|c|c}
    \hline
    Optimization Strategy & Model & Total Epochs & Relative $L_2$ Error \\
    \hline
    \multirow{3}{*}{SGL-Exp} & SGL-1-Exp & 550,000 & 8.12e-05 \\
                             & SGL-2-Exp & 550,000 & 7.42e-05 \\
                             & SGL-3-Exp & 550,000 & 9.00e-04 \\
    \hline
    \multirow{3}{*}{SGL-Piecewise} & SGL-1-Piecewise & 550,000 & 1.23e-04 \\
                                   & SGL-2-Piecewise & 550,000 & 1.70e-04 \\
                                   & SGL-3-Piecewise & 550,000 & 3.58e-04 \\
    \hline
    \textbf{TS-MGDL (Ours)} & \textbf{Two-Stage} & 550,000 & \textbf{9.65e-06} \\
    \hline
    \end{tabular}
\end{table}

\subsubsection{Robustness analysis against hybrid optimization strategies}
\label{sec:1d:robustness}

A widely adopted practice for achieving high-precision solutions in neural network-based PDE solvers is to couple the first-order Adam optimizer with the  Limited-memory Broyden–Fletcher–Goldfarb–Shanno (L-BFGS) method. Since the L-BFGS optimizer exhibits a faster convergence rate near local minima, it enables a more refined approximation of the target solution. Therefore, to evaluate the upper performance limit of the baseline methods for a fair comparison with TS-MGDL and ensure that any performance gap is due to the model architecture rather than the optimization routine, we employ a hybrid optimization approach for single grade learning method.

To evaluate the robustness of the TS-MGDL method against initialization stochasticity, we conducted repeated experiments across three independent random seeds: 12300, 12301, and 12302. To ensure a rigorous comparison regarding computational expenditure and optimization efficacy, we introduced a hybrid training protocol for the single-grade learning (SGL) baselines, denoted as SGL-1-H, SGL-2-H, and SGL-3-H. These models inherit the architectures of SGL-1, SGL-2, and SGL-3, respectively, as detailed in Table \ref{table:1D:Burgers:SGL:structure}, but undergo a dual-phase optimization process: an extensive Adam-based training phase of 550,000 epochs, followed by a fine-tuning stage of up to 1,000 L-BFGS iterations. For these hybrid baselines, an optimal hyperparameter selection strategy was implemented, where the initial learning rate for the Adam optimizer was carefully chosen from the candidate set $\{0.001, 0.002, 0.003\}$ to elicit peak performance. In parallel, the proposed TS-MGDL method employs the three-grade network configuration consistent with Table \ref{table:1D:Burgers:MGL:structure}. To maintain a consistent computational budget, the cumulative training length for TS-MGDL was also maintained at 550,000 epochs. The specific allocation of epochs and the learning rate schedules across the three grades and the final fine-tuning stage were specifically adjusted according to the different random seeds to optimize convergence, as summarized in Table \ref{tab:mgdl_hyperparams}.

\begin{table}[htbp]
  \centering
  \caption{Specific hyperparameter configurations for TS-MGDL across different random seeds. }
  \label{tab:mgdl_hyperparams}
  \resizebox{\columnwidth}{!}{ 
    \begin{tabular}{llcccc}
    \hline
    \textbf{Random Seed} & \textbf{Hyperparameter} & \textbf{Grade 1} & \textbf{Grade 2} & \textbf{Grade 3} & \textbf{Stage 2} \\
    \hline
    \multirow{2}{*}{12300} & Epochs & 50,000 & 100,000 & 100,000 & 300,000 \\
                           & Learning Rate & 1e-3 & 3e-4 & 2e-4 & 3e-4 \\
    \hline
    \multirow{2}{*}{12301} & Epochs & 100,000 & 100,000 & 100,000 & 250,000 \\
                           & Learning Rate & 1e-3 & 3e-4 & 2e-4 & 1e-4 \\
    \hline
    \multirow{2}{*}{12302} & Epochs & 80,000 & 80,000 & 80,000 & 310,000 \\
                           & Learning Rate & 1e-3 & 3e-4 & 2e-4 & 1e-4 \\
    \hline
    \end{tabular}%
  }
\end{table}

Table \ref{tab:error_comparison} summarizes the performance metrics across the three random seeds. The TS-MGDL method demonstrates superior stability and accuracy, achieving an average relative $L_2$ error of $\mathbf{2.78 \text{e-}05}$. This result is the lowest among all tested methods, outperforming even the best Single-Grade variant, SGL-1-H, which yields an average relative $L_2$ error of $4.92 \text{e-}05$.

Crucially, when comparing TS-MGDL to SGL-3-H, which shares a comparable network depth, TS-MGDL yields a relative error approximately one order of magnitude lower than that of SGL-3-H ($2.78 \text{e-}05$ vs. $2.40 \text{e-}04$). Furthermore, regarding computational efficiency, TS-MGDL requires less training time (10,981s) than the SGL-3-H model (11,264s). These results confirm that the TS-MGDL framework effectively leverages transfer learning to achieve higher accuracy with reduced computational cost compared to training a deep network from scratch, even when the latter is augmented with second-order optimization.

A comparison between the standard SGL results in Table \ref{table:1D:Burgers} and the hybrid results in Table \ref{tab:error_comparison} highlights the distinct advantages of the multi-grade strategy over purely optimizer-based enhancements. Specifically, using the results corresponding to the random seed 12300, we observe that while introducing L-BFGS optimization to the single-grade learning models (Table \ref{tab:error_comparison}) improved the accuracy of the shallower SGL-1 significantly compared to its Adam-only counterpart in Table \ref{table:1D:Burgers} (reducing the error from $5.75 \text{e-}04$ to $3.46 \text{e-}05$), this benefit diminished for the deeper SGL-3 architecture ($2.53 \text{e-}04$ vs. $2.47 \text{e-}04$). This stagnation indicates that simply adding second-order optimization is insufficient to overcome the training difficulties of deep networks initialized from scratch. In contrast, TS-MGDL effectively navigates this optimization landscape, achieving superior accuracy to even the best-performing hybrid SGL models. Finally, it is worth mentioning that the proposed TS-MGDL method is also compatible with the L-BFGS optimizer; the hybrid strategy of initializing with Adam and fine-tuning with L-BFGS can be seamlessly integrated into the training process of each grade to potentially further enhance precision.

\begin{table}[htbp]
  \centering
  \caption{Comparison of numerical errors, computational time, and memory usage between TS-MGDL and hybrid single-grade methods (SGL-*-H) across three random seeds.}
  \label{tab:error_comparison}
  \resizebox{\columnwidth}{!}{ 
    \renewcommand{\arraystretch}{1.3} 
    \begin{tabular}{lcccccc}
    \hline
    \multirow{2}{*}{\textbf{Method}} & \multicolumn{3}{c}{\textbf{Random Seed (Relative $L_2$ Error)}} & \multirow{2}{*}{\textbf{\shortstack{Relative $L_2$ Error\\(Mean $\pm$ Std. Dev.)}}} & \multirow{2}{*}{\textbf{\shortstack{Average\\Time (second)}}} & \multirow{2}{*}{\textbf{\shortstack{Average\\Memory (GB)}}} \\
    \cline{2-4}
     & 12300 & 12301 & 12302 &  &  &  \\
    \hline
    TS-MGDL & \textbf{9.65e-06} & \textbf{3.87e-05} & \textbf{3.50e-05} & \textbf{2.78e-05 $\pm$ 1.58e-05} & 10981.43 & 2.38 \\
    SGL-1-H & 3.46e-05 & 7.49e-05 & 3.82e-05 & 4.92e-05 $\pm$ 2.23e-05 &  4251.75 & 1.35 \\
    SGL-2-H & 1.99e-04 & 7.29e-05 & 4.20e-05 & 1.05e-04 $\pm$ 8.32e-05 & 6357.64 & 1.60 \\
    SGL-3-H & 2.47e-04 & 2.98e-04 & 1.74e-04 & 2.40e-04 $\pm$ 6.23e-05 & 11264.41 & 3.15 \\
    \hline
    \end{tabular}%
  }
\end{table}

\subsubsection{Empirical discussion on loss weighting}

A common challenge in training Physics-Informed Neural Networks (PINNs) is the imbalance between different loss components, which often leads to gradient pathologies. Theoretically, introducing penalty weights $\lambda_I$ and $\lambda_B$ in the loss function \eqref{loss:PINN_weighted} is intended to balance the gradients of the PDE residual with those of the initial and boundary conditions. However, our empirical investigations using the TS-MGDL framework reveal an interesting counter-intuitive phenomenon where manual weight tuning is largely unnecessary.

To investigate the impact of different weights $\lambda_I$ and $\lambda_B$ on the final solution accuracy, we designed a sensitivity analysis experiment using the Adam optimizer. We employed the same three-grade network configuration for TS-MGDL as detailed in Table \ref{table:1D:Burgers:MGL:structure}. The training schedule was standardized across all trials. In Stage 1, Grade 1, Grade 2, and Grade 3 were each trained for 100,000 epochs, with the initial learning rates set to 1.0e-03, 3.0e-04, and 2.0e-04, respectively. Subsequently, a Stage 2 fine-tuning phase was conducted for 250,000 epochs with a reduced initial learning rate of 1.0e-04. Throughout all training phases, an inverse time decay schedule with a rate of 1.0e-04 was applied. All other settings, including the activation function and sampling points, remained consistent with the previous examples. To ensure statistical robustness, we conducted repeated experiments across three independent random seeds (12300, 12301, and 12302) for three distinct weight configurations: heavy penalties ($\lambda_I=\lambda_B=30$), moderate penalties ($\lambda_I=\lambda_B=10$), and unit weights ($\lambda_I=\lambda_B=1$).

The statistical results of the test metrics are summarized in Table \ref{tab:ts_mgdl_results}. We observe a clear trend where increasing the penalty weights correlates with a degradation in model performance. Specifically, the unit weight configuration ($[1, 1]$) yields the highest accuracy with a mean test metric of 3.91e-05 and the lowest variance (std. dev. of 4.06e-06). In contrast, the heavy penalty configuration ($[30, 30]$) results in an error nearly an order of magnitude higher.

\begin{table}[htbp]
\centering
\setlength{\tabcolsep}{4pt}
\caption{Detailed relative $L_2$ error for TS-MGDL under different loss weights across three random seeds (mean $\pm$ standard deviation).}
\label{tab:ts_mgdl_results}
\begin{tabular*}{\linewidth}{@{\extracolsep{\fill}}lcccc}
\hline
\multirow{2}{*}{Weights ($\lambda_I, \lambda_B$)} & \multicolumn{3}{c}{Random Seed} & \multirow{2}{*}{\begin{tabular}[c]{@{}c@{}}Relative $L_2$ Error \\ (Mean $\pm$ Std. Dev.)\end{tabular}} \\ \cline{2-4}
 & 12300 & 12301 & 12302 & \\ \hline
$[30, 30]$ & 2.23e-04 & 1.49e-04 & 3.07e-04 & 2.26e-04 $\pm$ 0.79e-04 \\
$[10, 10]$ & 1.55e-04 & 4.89e-05 & 7.75e-05 & 9.38e-05 $\pm$ 5.49e-05 \\
$[1, 1]$   & 3.52e-05 & 3.87e-05 & 4.33e-05 & 3.91e-05 $\pm$ 0.41e-05 \\ \hline
\end{tabular*}
\end{table}

This result stands in contrast to standard PINN training, where large weights are frequently necessary to prevent the PDE loss from dominating the boundary loss and  the initial condition loss. The hierarchical training strategy of TS-MGDL appears to facilitate a more natural alignment of gradients. By refining the solution through progressive grades, the model may effectively avoid the severe gradient imbalances seen in single-grade architectures. While the exact theoretical mechanism by which the hierarchical structure mitigates the need for high penalty weights remains an open question, our empirical evidence suggests that the TS-MGDL framework is robust enough to operate effectively with a simplified, unweighted loss function.

\subsection{The 2D Burgers equation}
\begin{figure}[htbp]
  \centering
  \includegraphics[width=0.45\textwidth]{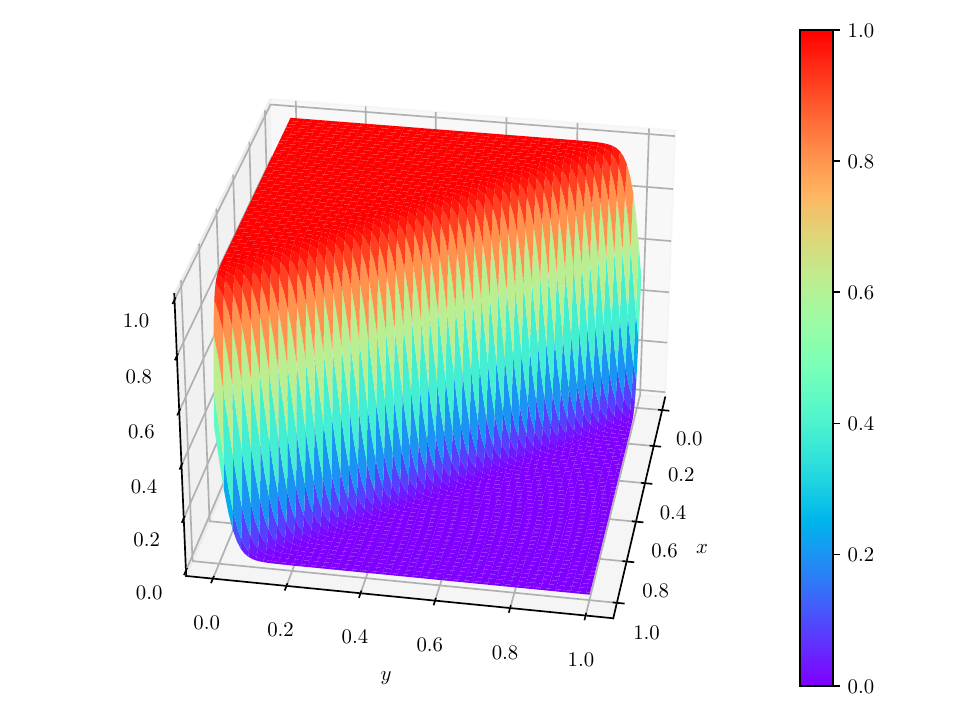}
  \caption{Exact solution of the 2D Burgers equation at $t=1$.}
  \label{fig:2d:exact}
\end{figure}

In this example, we consider the 2D Burgers equation over a square domain $\Omega:=[0, 1]^2$ in the form
\begin{equation}
\label{eq:Burgers:2D}
    \begin{gathered}
 \frac{{\partial}u}{\partial t}+u\frac{{\partial}u}{\partial x}+u\frac{{\partial}u}{\partial y}=\nu \left(\frac{{\partial}^2u}{\partial x^2}+\frac{{\partial}^2u}{\partial y^2}\right), \ \
  (x, y) \in \Omega, \ \ t \in(0, 1], 
 \end{gathered}
\end{equation}
where $u$ represents the velocity, $\nu := \frac{1}{Re}$ is the kinematic viscosity coefficient, and $Re$ is the Reynolds number. In this example, we set Reynolds number $Re := 100$.
The initial condition is given by:
\begin{equation}
    \label{init:Burgers:2D}
    u(0, x, y) = (1+e^{Re(x+y)})^{-1},\ \  (x, y) \in \Omega.
\end{equation} 
The boundary conditions are given by:
\begin{equation}
\label{boundary:Burgers:2D}
    \begin{aligned}
u(t, 0, y)=(1+e^{\frac{R e(y-t)}{2}})^{-1}, \ \ u(t, 1, y) & =(1+e^{\frac{R e(1+y-t)}{2}})^{-1}, y \in[0,1],\ \ t>0,\\
u(t, x, 0)=(1+e^{\frac{R e(x-t)}{2}})^{-1}, \ \ u(t, x, 1) & =(1+e^{\frac{R e(1+x-t)}{2}})^{-1}, x \in[0,1],\ \ t>0.\\
\end{aligned}
\end{equation}
The exact solution of the initial boundary value problem takes the form:
\begin{equation}
    u(x, y, t)=(1+e^{\frac{Re(x+y-t)}{2}})^{-1}, \quad(x, y) \in \Omega,\ \ t \geq 0 .
\end{equation}
We show in Figure \ref{fig:2d:exact} the exact solution $u(x,y,1)$ at $t=1$, which exhibits a significant increase across the line $x+y=1$.

In this example, the training point set is randomly generated using the Hammersley sampling method and consists of three components: $N_f := 20,000$ points located within the interior region $(0, 1] \times (0, 1)^2$, $N_0 := 1,000$ points satisfying the initial condition \eqref{init:Burgers:2D}, and $N_b := 2,000$ points satisfying the boundary conditions \eqref{boundary:Burgers:2D}. The testing point set consists of $15\times15$ equally spaced grid points in the spatial domain $[0, 1]^2$ and a total of 224 equally spaced sampling points in the time interval [0, 1]. Therefore, there are a total of $224 \times 15 \times 15 = 50,400$ points in the testing set.

Table \ref{table:2D:Burgers:SGL:structure} provides an overview of the network structures used in the single-grade learning (SGL) methods to solve the 2D Burgers equation. It shows that as the method number increases, the network structures become progressively deeper, with SGL-3 having the deepest network structure.

Table \ref{table:2D:Burgers:MGL:structure} presents the network structures used in the TS-MGDL method. Notably, Grades 1,  2, and 3 share the identical network structure with SGL-1, SGL-2, and SGL-3, respectively. In Table \ref{table:2D:Burgers:MGL:structure}, the asterisk (*) notation indicates that the layer's parameters were frozen during the first stage of training. It is noteworthy that during the second stage of training, the last eight hidden layers of the Grade 3 neural network are unfrozen and undergo retraining to further refine the model.

\begin{table}[htbp]
		\centering
\begin{tabular}{c|l}
 \hline
   Methods & Network structure  \\
   \hline
   SGL-1 &  [3, 1024, 1024, 1]      \\
   SGL-2 &  [3, 1024, 1024, 1024, 512, 512, 1]     \\
   SGL-3 &  [3, 1024, 1024, 1024, 512, 512, 512, 512, 256, 256, 1]  \\
   \hline
\end{tabular}
  \caption{Network structure of single-grade learning for the 2D Burgers equation.}
 \label{table:2D:Burgers:SGL:structure}
\end{table}

\begin{table}[htbp]
		\centering
\begin{tabular}{c|l}
 \hline
   Grade & Network structure  \\
   \hline
    1 &  [3, 1024, 1024, 1]      \\
    2 &  [3, 1024*, 1024*, 1024, 512, 512, 1]     \\
    3 &  [3, 1024*, 1024*, 1024*, 512*, 512*, 512, 512, 256, 256, 1]  \\
   \hline
\end{tabular}
  \caption{Network structure of multi-grade learning for the 2D Burgers equation.}
 \label{table:2D:Burgers:MGL:structure}
\end{table}

Table \ref{table:2D:Burgers:MG} presents the numerical results of TS-MGDL methods. It shows that the relative $L_2$ error decreases as the grade number increases in the first stage, and the second stage further improves the accuracy. The relative $L_2$ error of the TS-MGDL method is 4.41e-05. 

\begin{table}[htbp]
		\centering
\begin{tabular}{c|c|c|c|c}
 \hline
   Stage \& grade &  Learning rate & Decay rate & Epochs & Relative $L_2$ error  \\
   \hline
   Stage 1, Grade 1 & 1e-3  & 1e-4  & 10,000   & 1.63e-03 \\
   Stage 1, Grade 2 & 3e-4 & 1e-4   & 20,000  & 1.77e-04 \\   
   Stage 1, Grade 3 & 3e-4 & 1e-4   & 40,000  & 1.28e-04\\
   Stage 2          & 3e-4 & 7e-5   & 80,000  & \textbf{4.41e-05}\\
   \hline
\end{tabular}
  \caption{Numerical results of TS-MGDL for the 2D Burgers equation.}
 \label{table:2D:Burgers:MG}
\end{table}

Table \ref{table:2D:Burgers} compares the numerical results of TS-MGDL with those of SGL methods. Hyperparameters used to train
SGL neural networks include a learning rate of 1e-3, a decay rate of 1e-4, and 150,000 epochs for each method. Table \ref{table:2D:Burgers} shows that SGL-2 achieves the relative $L_2$ error 1.39e-04, which is the smallest relative $L_2$ error among the SGL methods. It shows that for single-grade learning, it may become harder to optimize as the network goes deeper, which can lead to a decrease in accuracy. The TS-MGDL method achieves a significantly smaller relative $L_2$ error than any of the SGL methods, demonstrating the effectiveness of multi-grade learning in solving the 2D Burgers equation.

\begin{table}[htbp]
    \centering
    \small
    \setlength{\tabcolsep}{4pt}
    \begin{tabular}{c|c|c|c|c|c}
        \hline
        Method & Learning rate & Epochs & Relative $L_2$ error & Time (min) & Memory (GB) \\
        \hline
        SGL-1 & 1e-3 & 150,000 & 2.00e-04 & 248.80 & 8.15 \\
        SGL-2 & 1e-3 & 150,000 & 1.39e-04 & 640.06 & 14.87 \\
        SGL-3 & 1e-3 & 150,000 & 1.46e-03 & 835.14 & 14.88 \\
        TS-MGDL & -- & 150,000 & \textbf{4.41e-05} & 910.17 & 14.88 \\
        \hline
    \end{tabular}
    \caption{Numerical results for the 2D Burgers equation.}
    \label{table:2D:Burgers}
\end{table}

\begin{figure}[htbp]
  \centering
  \subfloat[SGL-1]{
    \includegraphics[width=0.45\textwidth]{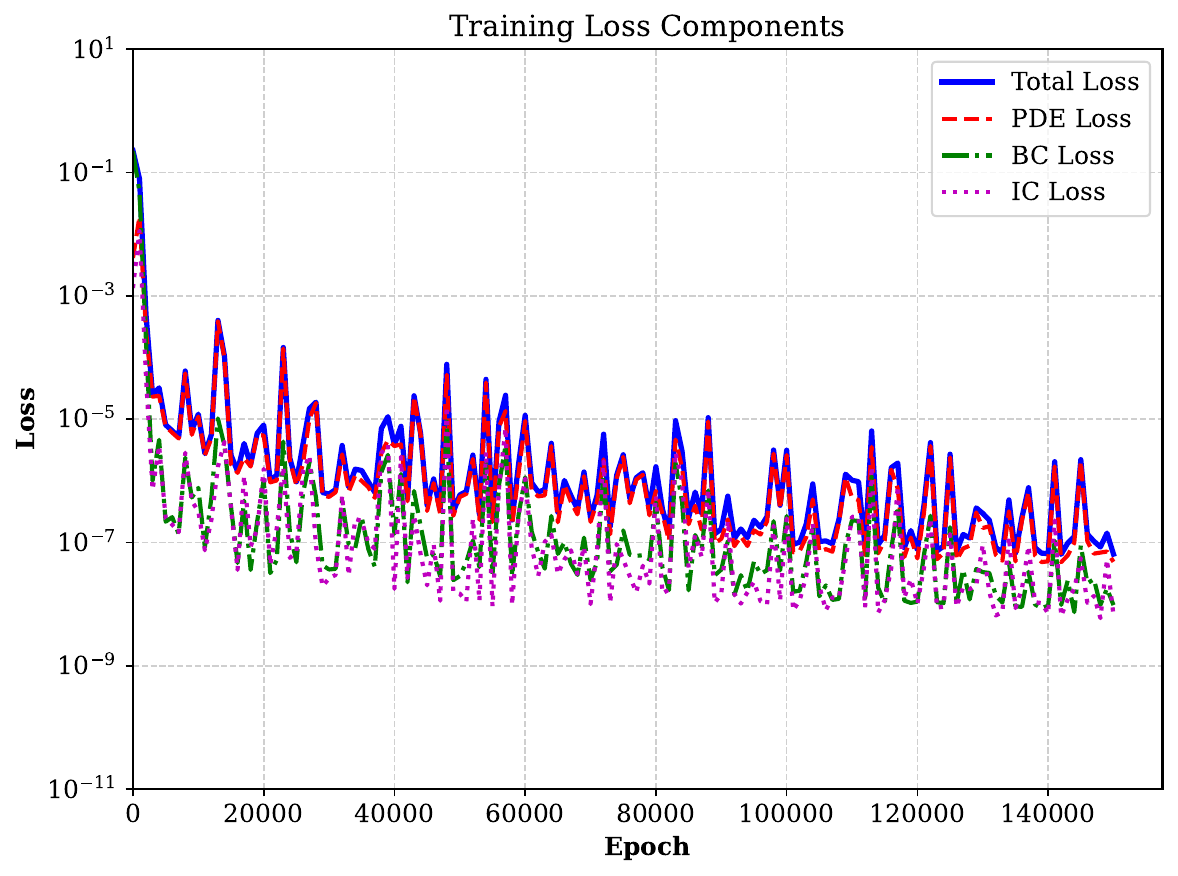}
    \label{fig:loss:sgl-1}
  }
  \subfloat[SGL-2]{
    \includegraphics[width=0.45\textwidth]{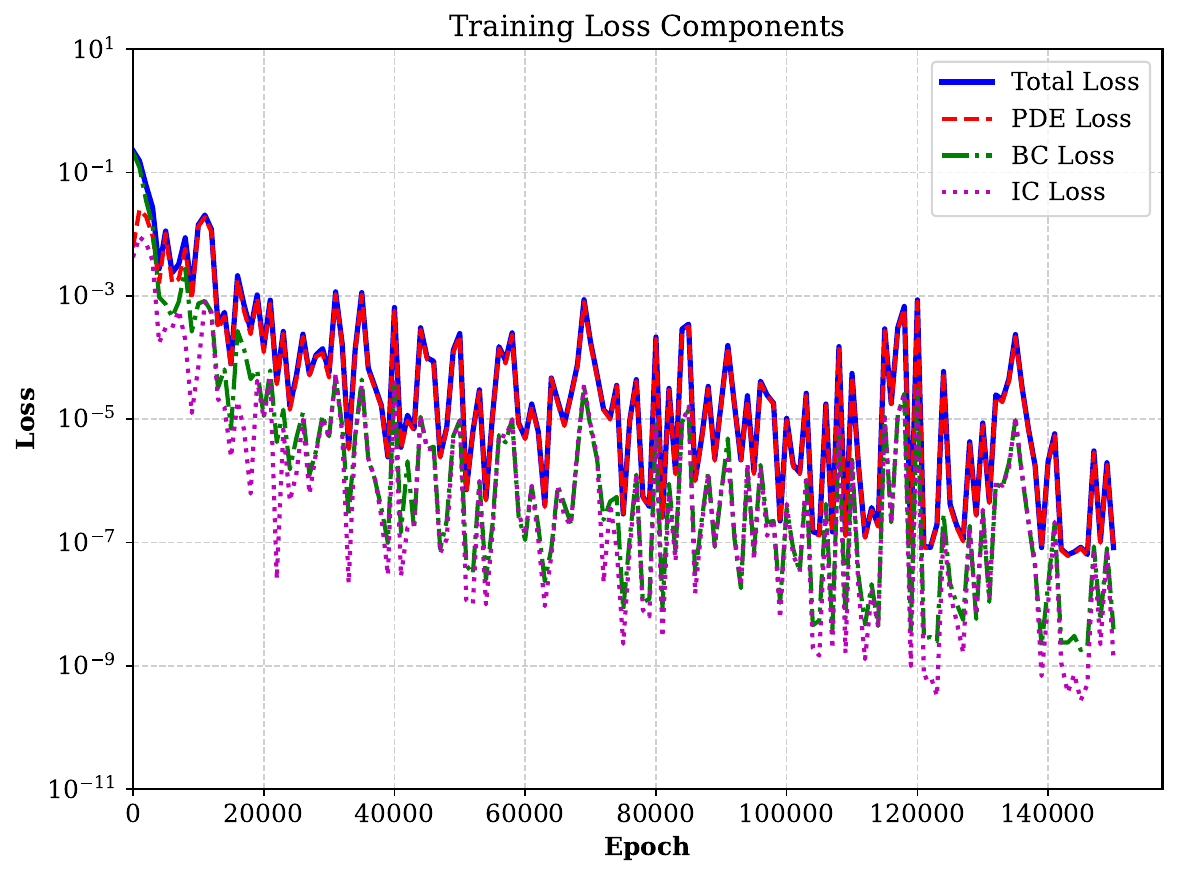}
    \label{fig:loss:sgl-2}
  }\\
    \subfloat[SGL-3]{
    \includegraphics[width=0.45\textwidth]{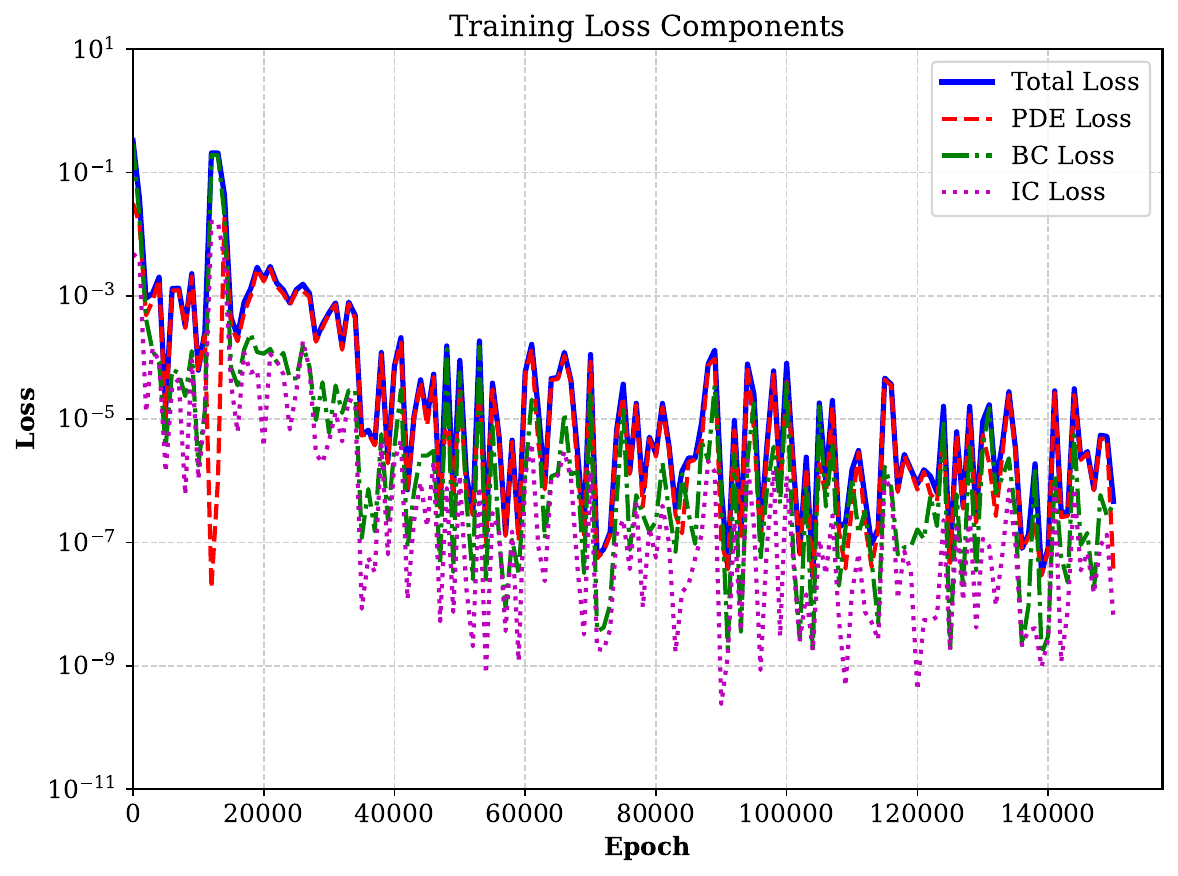}
    \label{fig:loss:sgl-3}
    }
   \subfloat[TS-MGDL]{
    \includegraphics[width=0.45\textwidth]{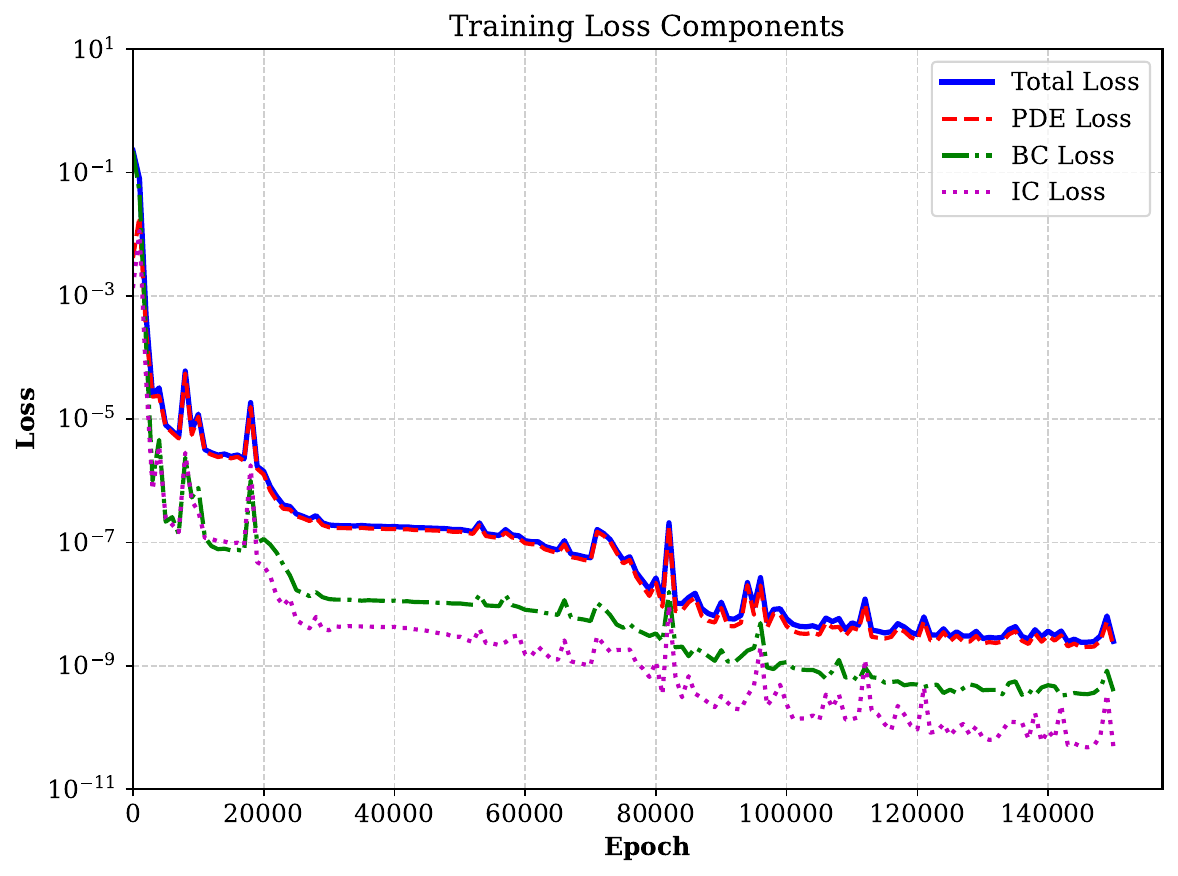}
    \label{fig:loss:mgl}
  }
  \caption{Training loss for SGL-1, SGL-2, SGL-3, and TS-MGDL for the 2D Burgers equation.}
  \label{fig:loss:2d}
\end{figure}

\begin{figure}[htpb]
  \centering
  
  %\subfloat[Predicted solution of SGL-1]{
  %  \includegraphics[width=0.45\textwidth]{2d-sgl-1-predict-1.pdf}
  %  \label{fig:ex2-sgl-1}
  %}
  \subfloat[Absolute error of SGL-1]{
    \includegraphics[width=0.45\textwidth]{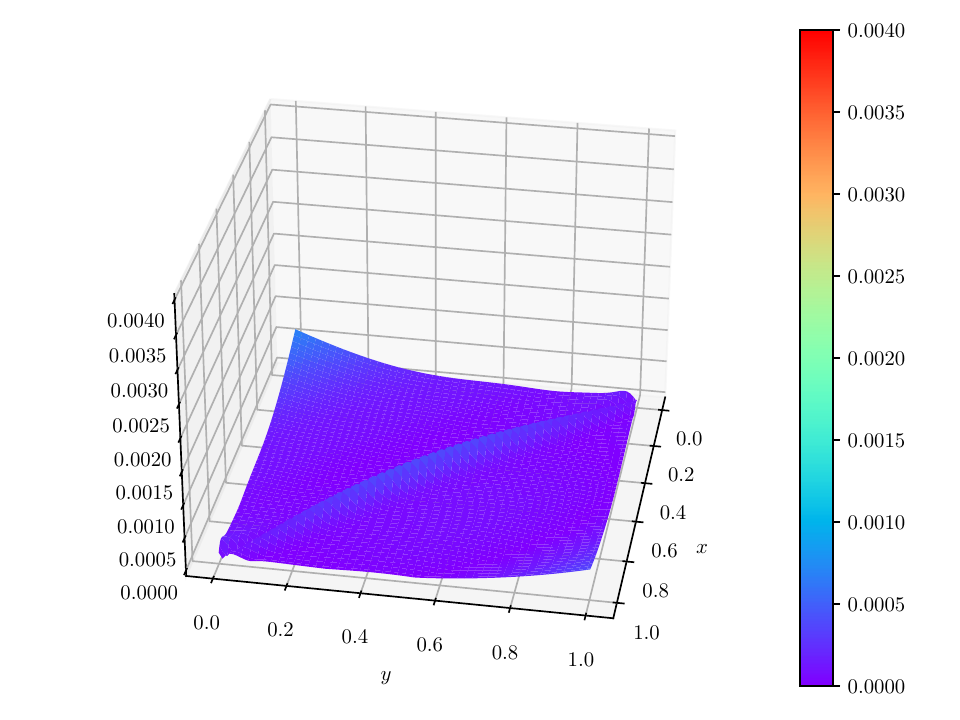}
    \label{fig:ex2-sgl-1:error}
  }
  %  \subfloat[Predicted solution of SGL-2]{
  %  \includegraphics[width=0.45\textwidth]{2d-sgl-2-predict-1.pdf}
  %  \label{fig:ex2-sgl-2}
  %}
  \subfloat[Absolute error of SGL-2]{
    \includegraphics[width=0.45\textwidth]{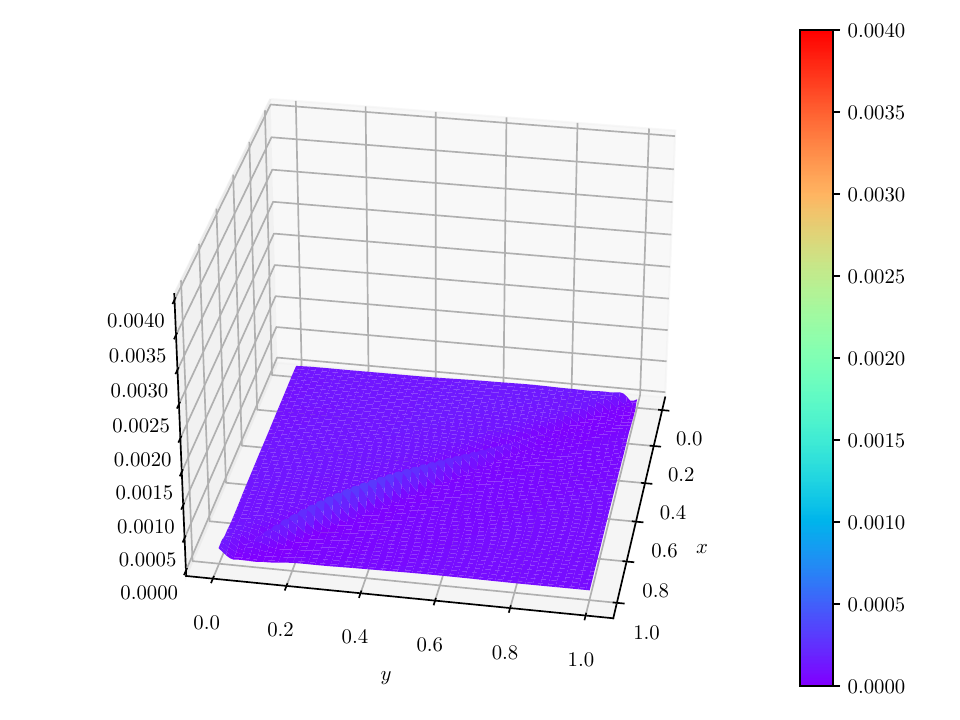}
    \label{fig:ex2-sgl-2:error}
  }\\
  %  \subfloat[Predicted solution of SGL-3]{
  %  \includegraphics[width=0.45\textwidth]{2d-sgl-3-predict-1.pdf}
  %  \label{fig:ex2-sgl-3}
  %}
  \subfloat[Absolute error of SGL-3]{
    \includegraphics[width=0.45\textwidth]{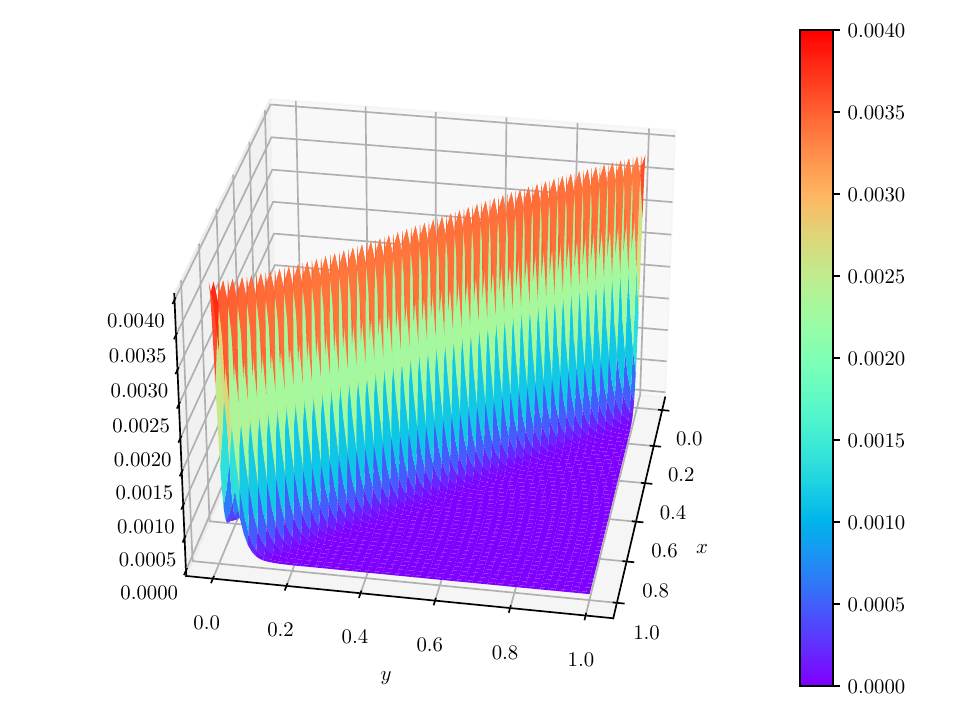}
    \label{fig:ex2-sgl-3:error}
  }
  %  \subfloat[Predicted solution of TS-MGDL]{
  %  \includegraphics[width=0.45\textwidth]{mgl-2d-predict-1.pdf}
  %  \label{fig:ex2-mgl}
  %}
  \subfloat[Absolute error of TS-MGDL]{
    \includegraphics[width=0.45\textwidth]{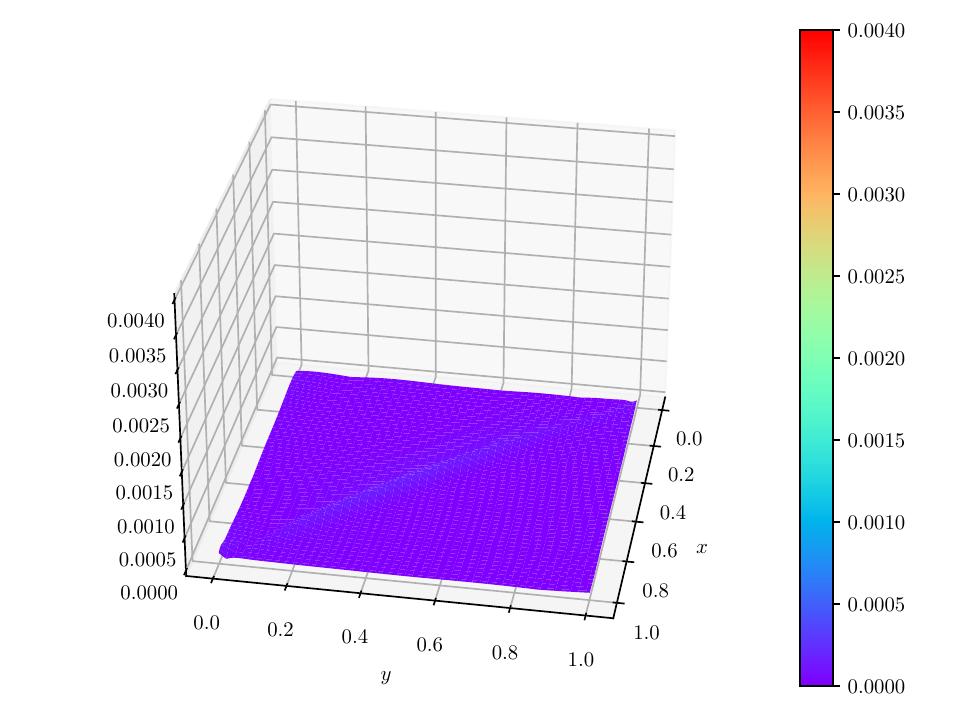}
    \label{fig:ex2-mgl:error}
  }

  \caption{Comparison between SGL methods and the TS-MGDL method for the 2D Burgers equation at $t = 1$.}
  \label{fig:2d}
\end{figure}

 Figure \ref{fig:loss:2d} shows the training loss of the SGL-1, SGL-2, SGL-3, and TS-MGDL methods for solving the 2D Burgers equation. The figures in Figure \ref{fig:2d} depict the approximation solutions and absolute errors of four methods. From the results, it can be observed that the TS-MGDL method exhibits smaller errors compared to the single-grade learning (SGL) methods. The loss curves of the SGL methods show significant oscillations and slow convergence. On the other hand, the loss curve of the TS-MGDL method demonstrates a smoother and faster decrease in loss.  Figure \ref{fig:loss:2d}\subref{fig:loss:mgl} confirms the theoretical results stated in Propositions \ref{thm:loss} and \ref{thm:loss-TS} for the 2D equation. These findings suggest that the TS-MGDL method outperforms the SGL methods in terms of accuracy and convergence speed.

\begin{figure}[htbp]
    \centering
    \includegraphics[width=0.8\textwidth]{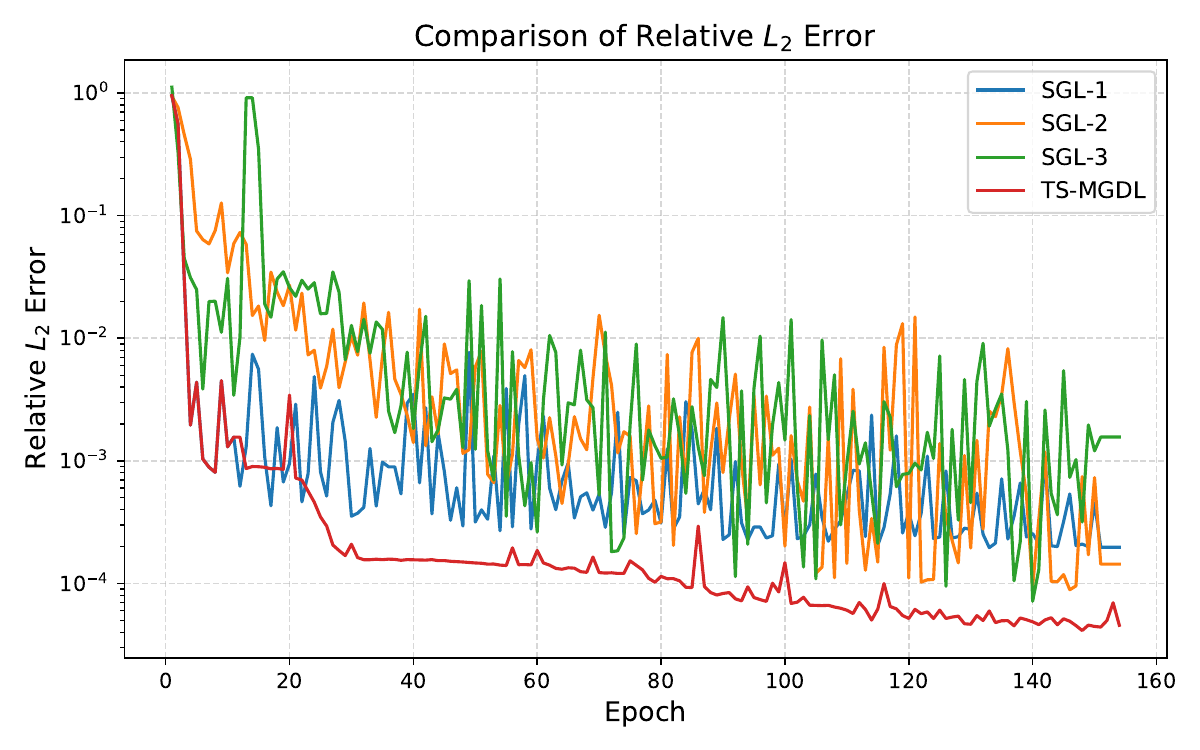}
    \caption{Relative $L_2$ error comparison between TS-MGDL and SGL methods for the 2D Burgers equation.}
    \label{fig:2d-burgers}
\end{figure}

To further assess the predictive performance, Figure \ref{fig:2d-burgers} compares the relative $L_2$ error trajectories computed on the testing dataset. This metric reflects the capability of the neural networks to generalize beyond the training locations. It is evident that the TS-MGDL method significantly outperforms the SGL approaches, maintaining a stable and rapid descent to a relative error of 4.41e-05. The results imply that while single-grade learning (e.g., SGL-3) struggles with optimization and generalization in deeper architectures, the multi-grade learning framework robustly captures the complex solution features, leading to high-fidelity predictions.

\subsection{The 3D Burgers equation}
\begin{figure}[htbp]
  \centering
  \includegraphics[width=0.4\textwidth]{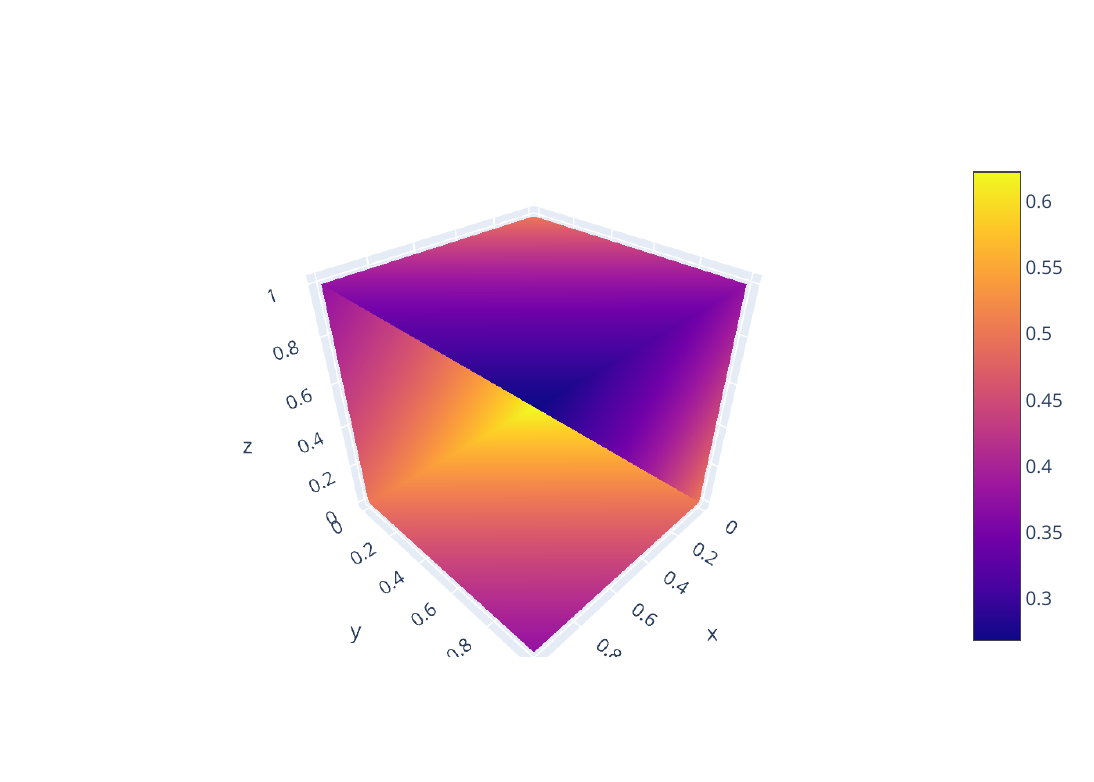}
  \caption{Exact solution of the 3D Burgers equation at $t=1$.}
  \label{fig:3d:exact}
\end{figure}

In this example, we consider solving the initial boundary value problem of the 3D Burgers equation defined in $\Omega := [0, 1]^3$, which can be expressed as follows:
\begin{equation}
\label{eq:Burgers:3D}
\begin{aligned}
&u_t + u u_x + u u_y + u u_z = \nu\left(u_{x x} + u_{y y} + u_{z z}\right) + f, \ \ (t, x, y, z) \in  (0, 1] \times \Omega, \\    
&u(t, x, y, z) = g(t, x, y, z), \quad t \in   (0, 1],  (x, y, z) \in \partial\Omega, \\
&u(0, x, y, z) = h(x, y, z), \quad (x, y, z) \in \Omega,
\end{aligned}
\end{equation}
where $\nu:=\frac{1}{Re}>0$ represents the viscosity coefficient,  and $Re$ is the Reynolds number. 
The exact solution of the equation \eqref{eq:Burgers:3D}  is given by:
\begin{equation}
\label{eq:exact:3D}
    u(t, x, y, z)=(1+\mathrm{e}^{\frac{x+y+z-t}{2 \nu}})^{-1}.
\end{equation}
In the experiment, the function values of $f$, $g$, and $h$ in equation \eqref{eq:Burgers:3D} were computed according to the exact solution \eqref{eq:exact:3D}. In this example, the Reynolds number $Re$ is set to 1. The exact solution of the 3D Burgers equation at $t=1$ is shown in Figure \ref{fig:3d:exact}.

The training point set is generated randomly using the Hammersley sampling method. It is structured as follows: $N_f := 10,000$ points are located within the interior region $(0, 1] \times (0, 1)^3$, $N_b := 2,000$ points satisfy the boundary condition, and $N_0 := 1,000$ points meet the initial condition. As for the testing point set, it consists of a uniform grid comprising $7^3$ points in the spatial domain $[0, 1]^3$, along with 265 equally spaced sampling points in the time interval [0, 1]. This combination results in a total of $265 \times 7^3 = 90,895$ points in the testing set.

Table \ref{table:3D:Burgers:SGL:structure} illustrates the network structures employed by the SGL methods. As the method number increases, the network structures of the SGL methods become progressively deeper. Table \ref{table:3D:Burgers:MGL:structure} presents the network structures of the TS-MGDL method. In this particular example, the TS-MGDL method consists of three grades. The network structure of Grade 1,  2, and 3 is identical to that of SGL-1, SGL-2, and SGL-3, respectively. In Table \ref{table:3D:Burgers:MGL:structure}, the asterisk (*) notation indicates that the parameters of the corresponding layer were kept fixed during the first stage of training.  During the second stage of training, the last eight hidden layers of the Grade 3 neural network are unfrozen and retrained.

\begin{table}[htbp]
		\centering
\begin{tabular}{c|l}
 \hline
   Methods & Network structure  \\
   \hline
   SGL-1 &   [4, 256, 256, 256, 256, 1]    \\
   SGL-2 &  [4, 256, 256, 256, 256, 512, 512, 256, 1]     \\
   SGL-3 &   [4, 256, 256, 256, 256, 512, 512, 256, 512, 512, 512, 256, 1] \\
   \hline
\end{tabular}
  \caption{Network structure of single-grade learning for the 3D Burgers equation.}
 \label{table:3D:Burgers:SGL:structure}
\end{table}

\begin{table}[htbp]
		\centering
\begin{tabular}{c|l}
 \hline
   Grade & Network structure  \\
   \hline
    1 &   [4, 256, 256, 256, 256, 1]    \\
    2 &  [4, 256*, 256*, 256*, 256*, 512, 512, 256, 1]     \\
    3 &   [4, 256*, 256*, 256*, 256*, 512*, 512*, 256*, 512, 512, 512, 256, 1] \\
   \hline
\end{tabular}
  \caption{Network structure of multi-grade learning for the 3D Burgers equation.}
 \label{table:3D:Burgers:MGL:structure}
\end{table}

Table \ref{table:3D:Burgers:MG} shows the numerical results of the TS-MGDL method.  %They demonstrate that the TS-MGDL method outperforms the SGL methods in terms of achieving a smaller relative $L_2$ error within the same number of training epochs. 
In Stage 1, the method progressively refines its approximate solution of the 3D Burgers equation through different grades, resulting in decreasing relative $L_2$ errors. In Stage 2, the TS-MGDL method further enhances its accuracy, attaining an impressive relative $L_2$ error of 8.90e-05.  The TS-MGDL method's ability to refine and adapt its learning through multiple stages and grades enables it to capture intricate dynamics, thereby resulting in enhanced accuracy when solving the 3D Burgers equation.

\begin{table}[htbp]
		\centering
\begin{tabular}{c|c|c|c|c}
 \hline
   Stage \& grade &  Learning rate & Decay rate & Epochs & Relative $L_2$ error  \\
   \hline
   Stage 1, Grade 1 & 4e-4 & 1e-4  & 3,000  & 2.95e-03\\
   Stage 1, Grade 2 & 4e-4 & 1e-4   & 6,000  & 7.39e-04\\   
   Stage 1, Grade 3 & 4e-4 & 1e-4  & 6,000  & 3.92e-04\\
   Stage 2 & 5e-4 & 1e-4  & 25,000  & \textbf{8.90e-05}\\
   \hline
\end{tabular}
  \caption{Numerical results of multi-grade learning for the 3D Burgers equation.}
 \label{table:3D:Burgers:MG}
\end{table}

Table \ref{table:3D:Burgers} presents a comparison of the numerical results obtained from SGL-1, SGL-2, SGL-3, and TS-MGDL methods for solving the 3D Burgers equation. The results demonstrate that the TS-MGDL method outperforms the SGL methods in terms of achieving a smaller relative $L_2$ error within the same number of training epochs. For SGL-1, SGL-2, and SGL-3, the learning rate is set to 4e-4, the decay rate is 1e-4, and the training is conducted for 40,000 epochs. Among the SGL methods, SGL-2 achieves the highest accuracy with a relative $L_2$ error of 1.77e-04. However, the TS-MGDL method exhibits the best overall performance, surpassing all SGL methods, with a relative $L_2$ error of 8.90e-05.

\begin{table}[htbp]
    \centering
    \small
    \setlength{\tabcolsep}{4pt}
    \begin{tabular}{c|c|c|c|c|c}
        \hline
        Method & Learning rate & Epochs & Relative $L_2$ error & Time (min) & Memory (GB) \\
        \hline
        SGL-1 & 4e-4 & 40,000 & 1.06e-03 & 22.56 & 2.57 \\
        SGL-2 & 4e-4 & 40,000 & 1.77e-04 & 61.59 & 4.59 \\
        SGL-3 & 4e-4 & 40,000 & 3.35e-04 & 105.58 & 8.57 \\
        TS-MGDL & -- & 40,000 & \textbf{8.90e-05} & 108.28 & 14.87 \\
        \hline
    \end{tabular}
    \caption{Numerical results for the 3D Burgers equation.}
    \label{table:3D:Burgers}
\end{table}

\begin{figure}[htbp]
  \centering
  \subfloat[SGL-1]{
    \includegraphics[width=0.45\textwidth]{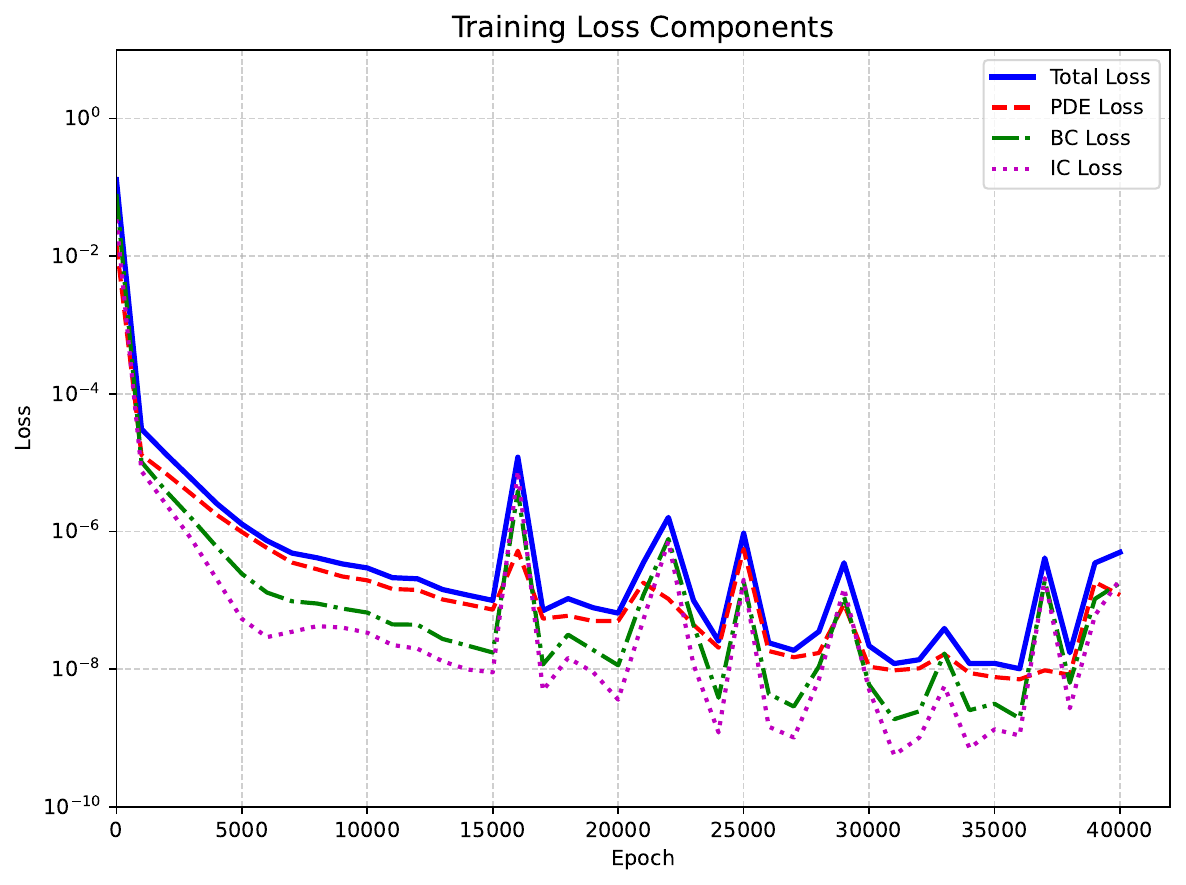}
    \label{fig:loss:sgl-1:3d}
  }
  \subfloat[SGL-2]{
    \includegraphics[width=0.45\textwidth]{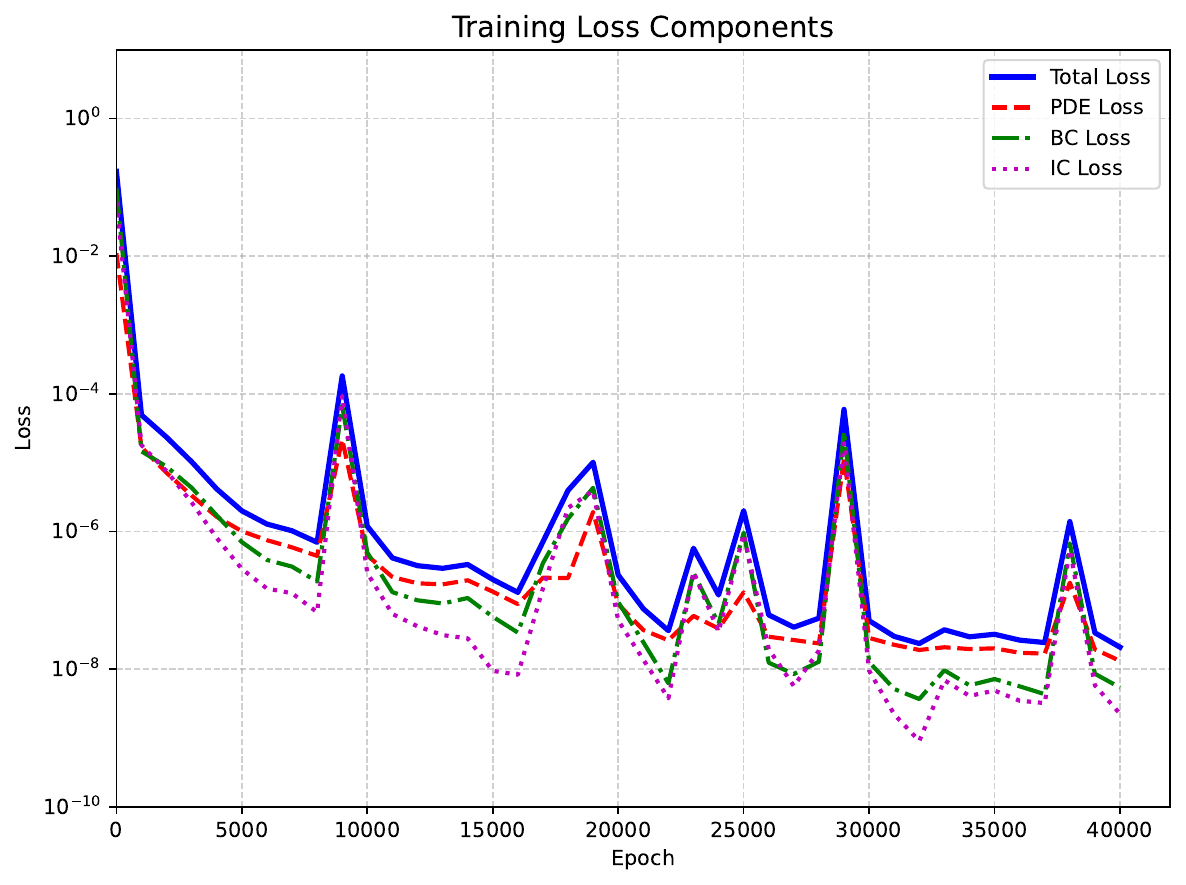}
    \label{fig:loss:sgl-2:3d}
  }\\
    \subfloat[SGL-3]{
    \includegraphics[width=0.45\textwidth]{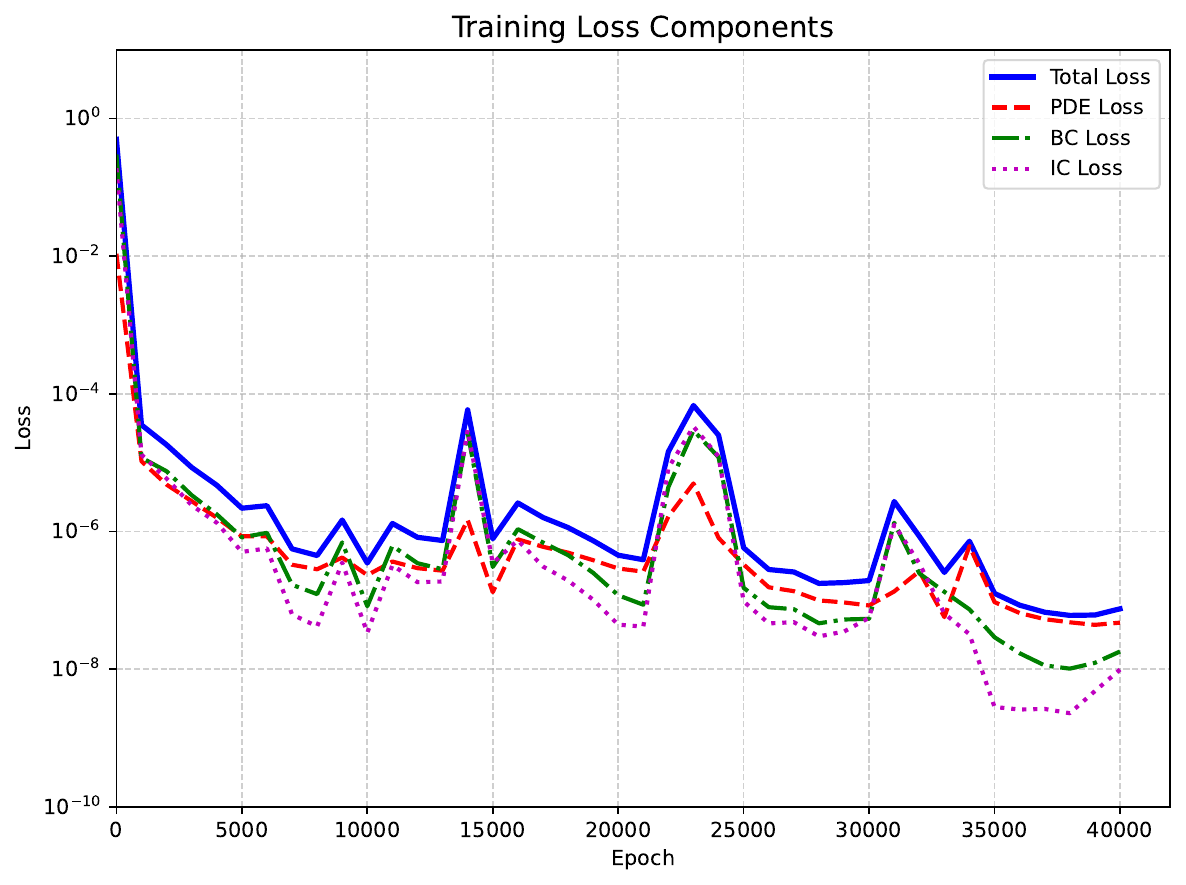}
    \label{fig:loss:sgl-3:3d}
    }
   \subfloat[TS-MGDL]{
    \includegraphics[width=0.45\textwidth]{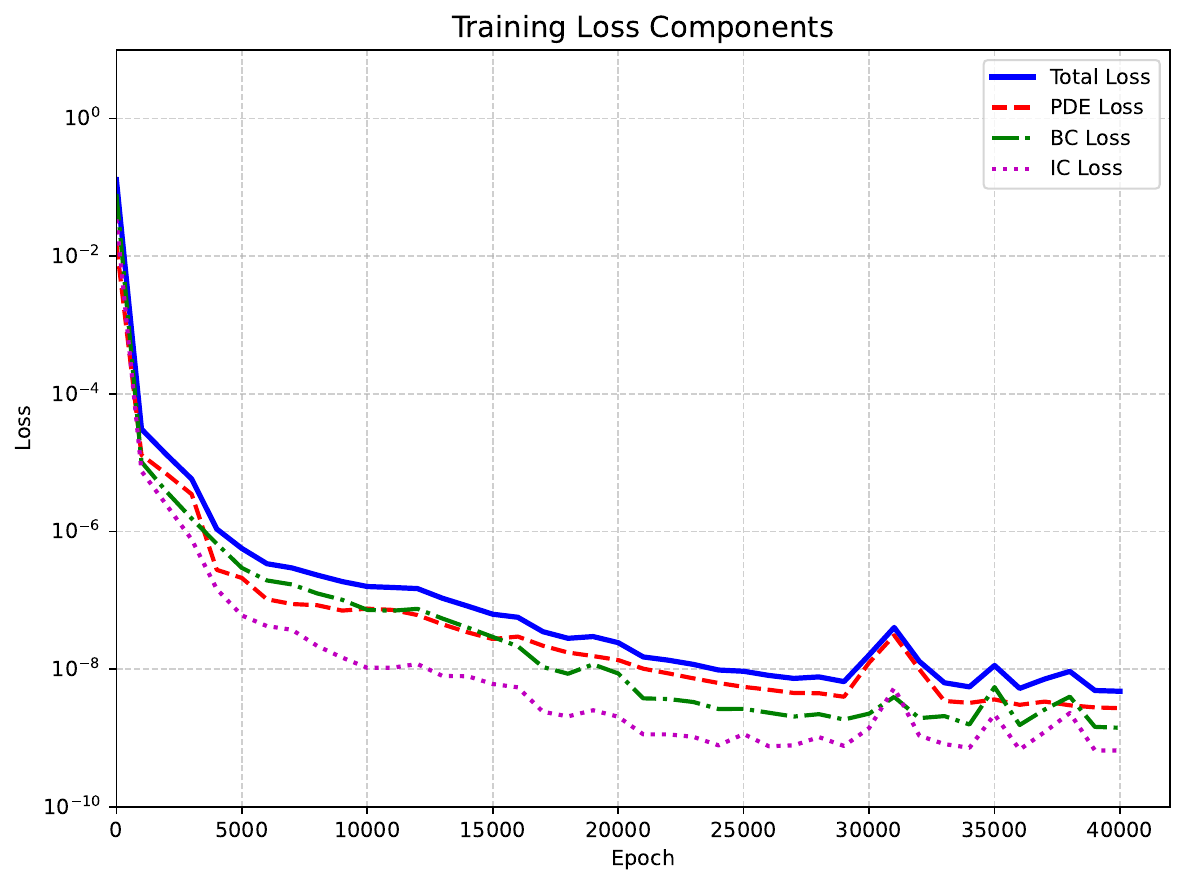}
    \label{fig:loss:mgl:3d}
  }
  \caption{Training loss for SGL-1, SGL-2, SGL-3, and TS-MGDL for the 3D Burgers equation.}
  \label{fig:loss:3d}
\end{figure}

\begin{figure}[htbp]
    \centering
    \includegraphics[width=0.8\textwidth]{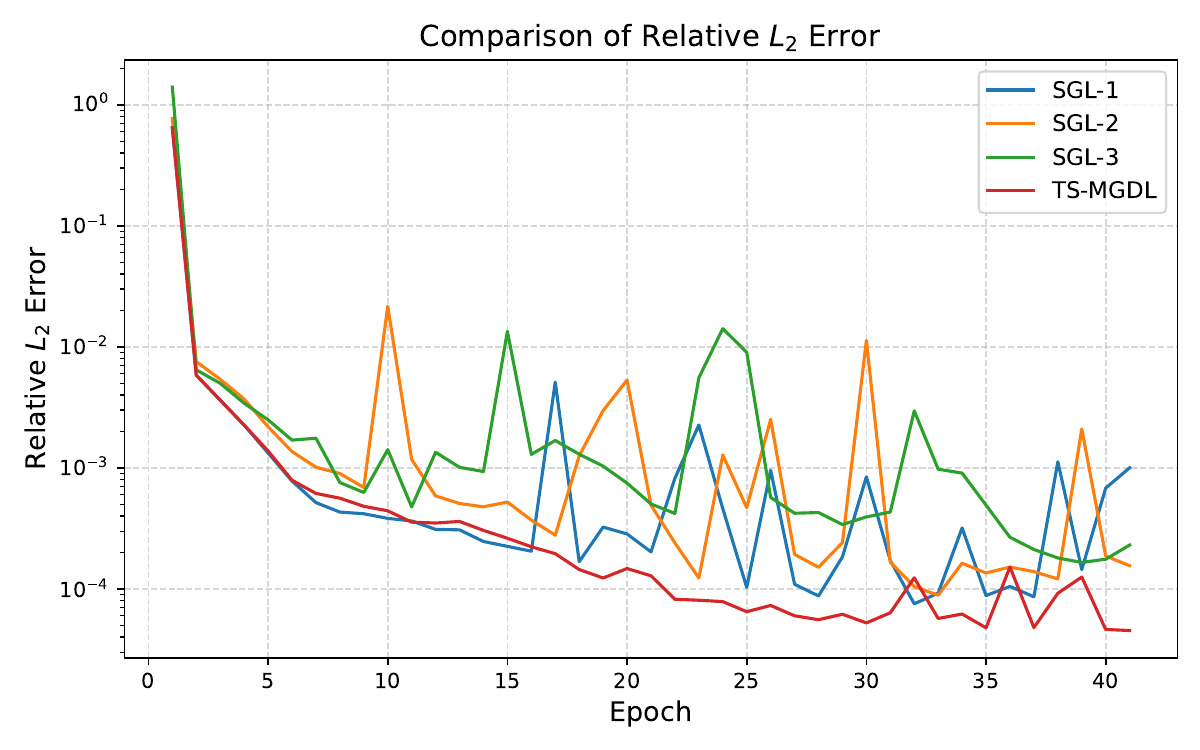}
    \caption{Relative $L_2$ error comparison between TS-MGDL and SGL methods for the 3D Burgers equation.}
    \label{fig:3d-burgers}
\end{figure}

Figure \ref{fig:loss:3d} illustrates the training loss curves for the SGL-1, SGL-2, SGL-3, and TS-MGDL methods. We can see that the TS-MGDL method demonstrates a significantly faster rate of decline in training loss compared to the SGL methods. The TS-MGDL method achieves the lowest overall loss among the evaluated methods. Moreover, the loss curve for TS-MGDL exhibits minimal oscillations, indicating a stable convergence process.  Figure \ref{fig:3d-burgers} presents the evolution of the relative $L_2$ error. It can be observed that the TS-MGDL method achieves a stable and fast convergence to the lowest error level, whereas the SGL methods exhibit varying degrees of oscillation and stagnation, failing to reach comparable accuracy.

\begin{figure}[htbp]
  \centering
  
  %\subfloat[Predicted solution of SGL-1]{
  %  \includegraphics[width=0.45\textwidth]{3d-sgl-1-pred-1.pdf}
  %  \label{fig:ex3-sgl-1}
  %}
  \subfloat[Absolute error of SGL-1]{
    \includegraphics[width=0.3\textwidth]{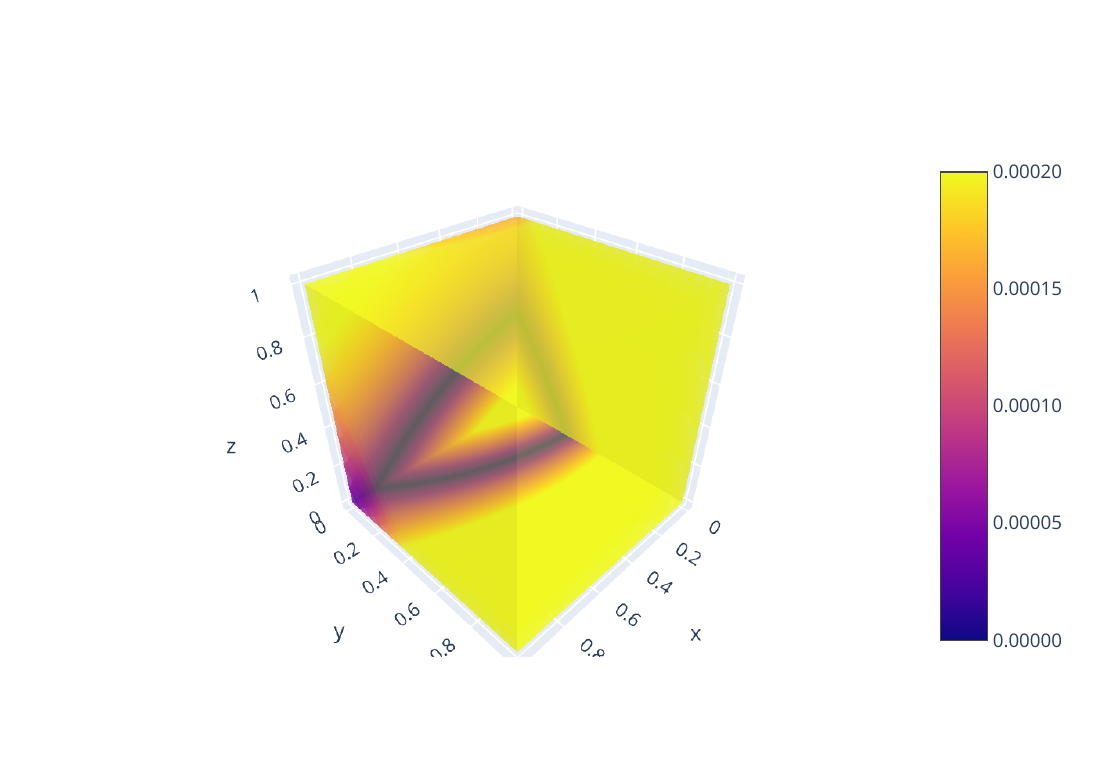}
    \label{fig:ex3-sgl-1:error}
  }
  %  \subfloat[Predicted solution of SGL-2]{
  %  \includegraphics[width=0.45\textwidth]{3d-sgl-2-pred-1.pdf}
  %  \label{fig:ex3-sgl-2}
  %}
  \subfloat[Absolute error of SGL-2]{
    \includegraphics[width=0.3\textwidth]{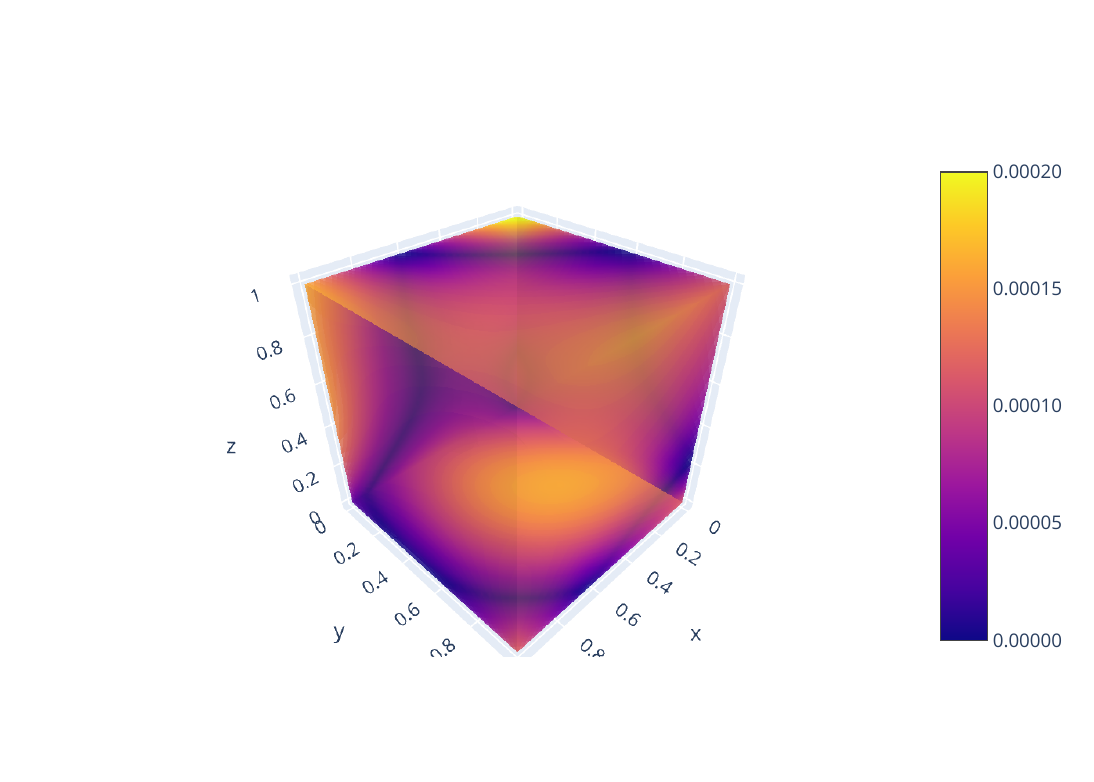}
    \label{fig:ex3-sgl-2:error}
  }\\
  %  \subfloat[Predicted solution of SGL-3]{
  %  \includegraphics[width=0.45\textwidth]{3d-sgl-3-predict-1.pdf}
  %  \label{fig:ex3-sgl-3}
  %}
  \subfloat[Absolute error of SGL-3]{
    \includegraphics[width=0.3\textwidth]{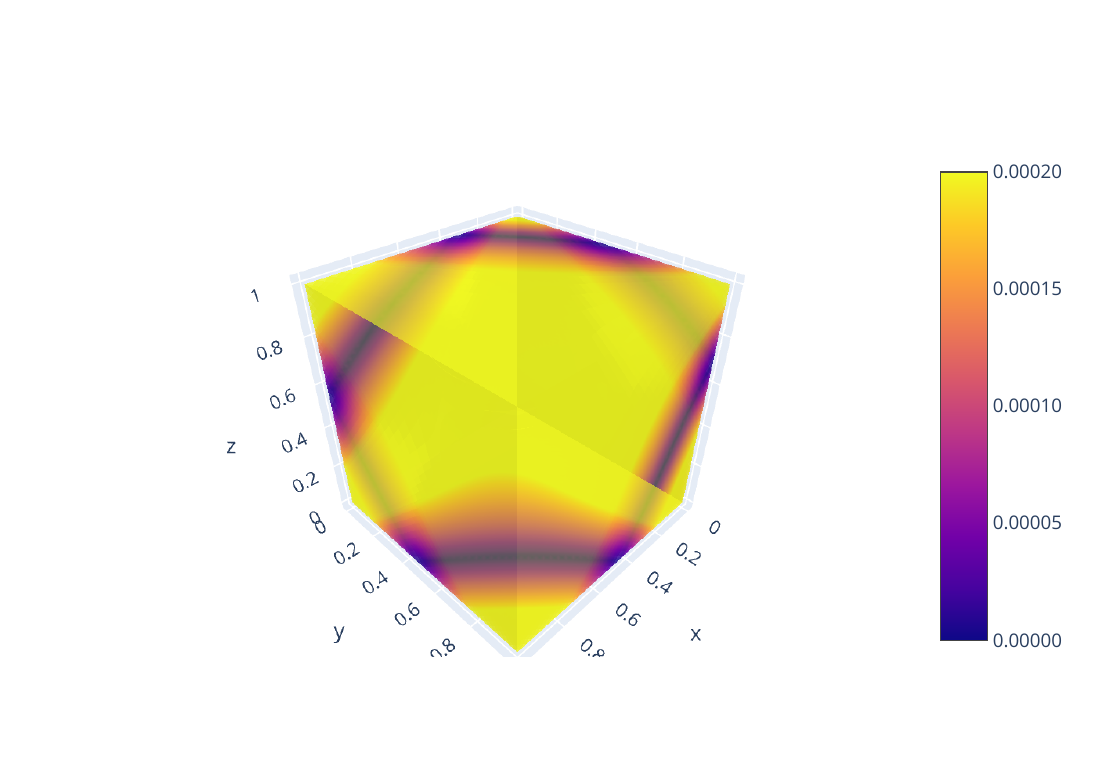}
    \label{fig:ex3-sgl-3:error}
  }
  %  \subfloat[Predicted solution of TS-MGDL]{
  %  \includegraphics[width=0.45\textwidth]{mgl-3d-pred-1.pdf}
  %  \label{fig:ex3-mgl}
  %}
  \subfloat[Absolute error of TS-MGDL]{
    \includegraphics[width=0.3\textwidth]{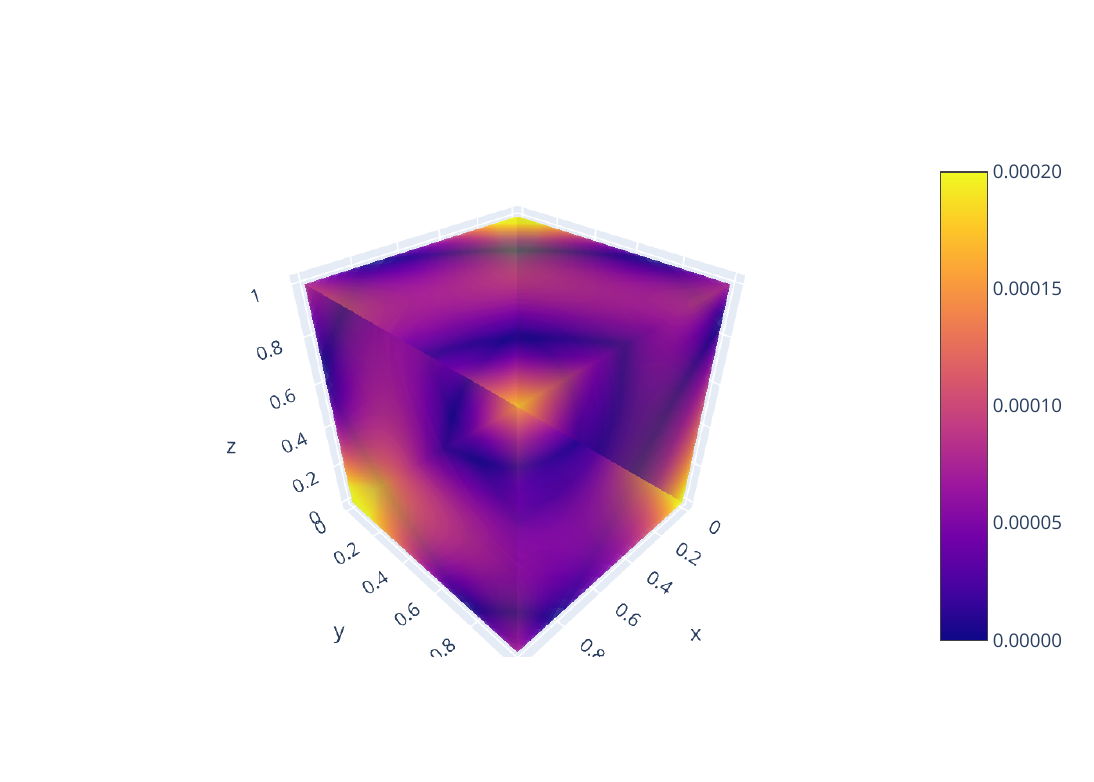}
    \label{fig:ex3-mgl:error}
  }

  \caption{Comparison between SGL methods and TS-MGDL method for the 3D Burgers equation at $t = 1$.}
  \label{fig:3d}
\end{figure}

Figure \ref{fig:3d} shows the computational results and error plots of the SGL-1, SGL-2, SGL-3, and TS-MGDL methods for the 3D Burgers equation at $t = 1$. As can be seen in Figure \ref{fig:3d},  the TS-MGDL method provides a better approximation of the true solution, with smaller absolute errors. Once again, the theoretical results stated in Propositions \ref{thm:loss} and \ref{thm:loss-TS} are verified by Figure \ref{fig:loss:3d} \subref{fig:loss:mgl:3d}.

In summary, the experimental results for the three examples presented above confirm that the proposed TS-MGDL method enables effective
learning of solutions of Burgers equations and outperforms existing single-grade deep learning methods in predictive accuracy and convergence. In particular, the predictive errors of single-grade deep learning are larger than those of the multi-grade deep learning
in 26-60, 4-31, and 3-12 times, for the 1D, 2D, and 3D equations, respectively. The numerical experiments show that the proposed TS-MGDL method can effectively capture the intrinsic multiscale components of the solutions near regions with steep gradients.

\section{Numerical Findings} We discuss in this section the findings from our numerical experiments. The TS-MGDL method, which integrates the concepts of greedy learning, residual learning, pretraining, and fine-tuning, demonstrates significant effectiveness in solving PDEs. The numerical experiments reveal four key advantages:

%\textbf{Improved Global Optima:} The two-stage approach helps in finding better global optima. In the first stage, the TS-MGDL method incrementally increases the learning grades, allowing for efficient local optimization and the discovery of a good initial solution. The second stage involves unfreezing and additional training, which explores a broader parameter space and refines the model. By considering a wider range of possible solutions, the method increases the chance of finding better global optima, leading to improved accuracy.

\textbf{Accelerated Convergence:}  The TS-MGDL method accelerates the convergence of the optimization process. By training on relatively shallow networks in the first stage, the method helps mitigate issues like gradient vanishing and exploding, resulting in stable and effective learning. As the complexity and capacity of the network are gradually increased in the second stage, the method facilitates the convergence of deeper networks, leading to faster optimization and improved convergence rates.

From Figures \ref{fig:loss:1d}, \ref{fig:loss:2d}, and \ref{fig:loss:3d}, we observe that the SGL methods often exhibit noticeable oscillations and slower convergence in their loss curves. In contrast, the loss curve of the TS-MGDL method demonstrates a smoother and more rapid decrease in loss. Additionally, the TS-MGDL method ensures continuous loss reduction across each training stage, maintaining efficient optimization throughout the process.

\textbf{Implicit Regularization:} The TS-MGDL method incorporates implicit regularization throughout its training stages. In the first stage, the algorithm gradually increases the complexity of the network, preventing premature attempts to learn intricate relationships that could lead to overfitting. In the second stage, only the unfrozen layers are allowed to adapt and fine-tune their parameters. This approach enables the model to refine its learned representations in a controlled manner, balancing the risk of overfitting with the need for improved prediction accuracy. Additionally, during the first stage of training in TS-MGDL, there are relatively fewer training epochs for each grade compared to single-grade training. This strategy mirrors the concept of early stopping in neural network training, providing an implicit regularization effect.

The implicit regularization capability of TS-MGDL can be observed by comparing the results of regularized single-grade neural networks. We apply $L_2$ regularization to the SGL method. The $L_2$ regularization term is defined as:
$$
L_R := \frac{\lambda_R}{2} \sum_{i}\sum_{k}\sum_{j} (W_{i}^{k,j})^2,
$$
where $W_{i}$ is the weight matrix in the $i$-th layer and $\lambda_R$ is the regularization parameter. We solve the Burgers equation of 1D, 2D, and 3D using the regularized SGL method and compare their accuracy with that of the MGDL method.

We report the numerical results for the equations of 1D, 2D, and 3D in Tables  \ref{table:1D:Burgers:regularization},  \ref{table:2D:Burgers:SG:Regularization}, and  \ref{table:3D:Burgers:regularization}, respectively. In these tables, we use SGL-X-R-1D (resp. 2D, 3D) to denote the regularized single-grade methods with networks SGL-X, X=1,2,3, for the 1D (resp. 2D, 3D) equation. All these examples show that the TS-MGDL method outperforms significantly the regularized SGL method in preventing overfitting, indicating its superior performance in implicit regularization.

%
%Tables \ref{table:1D:Burgers:regularization},  \ref{table:2D:Burgers:SG:Regularization} and  \ref{table:3D:Burgers:regularization} show that the MGDL method outperforms the regularized SGL method in all cases, indicating its superior performance in implicit regularization.
%Table \ref{table:1D:Burgers:regularization} presents the numerical results of the 1D Burgers equation after regularization. SGL-1-R-1D, SGL-2-R-1D, and SGL-3-R-1D refer to the single-grade methods SGL-1, SGL-2, and SGL-3 for the 1D equation with $L_2$ regularization, respectively. The regularization parameter $\lambda_R$ is set to be 1e-7 for all three methods.  Their relative $L_2$ errors are 2.10e-04, 8.56e-05, and 3.77e-05, respectively. Notably, TS-MGDL achieves the lowest relative $L_2$ error of 9.65e-06 among all the methods, indicating its superior performance in implicit regularization.
%
%
\begin{table}[htpb]
		\centering
\begin{tabular}{c|c|c|c|c}
 \hline
   Methods &  Learning rate &  Decay rate & Epochs & Relative $L_2$ error  \\
   \hline
   SGL-1-R-1D &3e-4 & 1e-4 & 550,000 & 2.10e-04 \\
   SGL-2-R-1D & 3e-4  & 1e-4     & 550,000    & 8.56e-05   \\ 
   SGL-3-R-1D & 3e-4  & 1e-4    & 550,000   & 3.77e-05  \\
   TS-MGDL  & -- & -- &550,000   & \textbf{9.65e-06}
   \\
   \hline
\end{tabular}
  \caption{Comparison of regularized single-grade methods with TS-MGDL for the 1D Burgers equation.}
 \label{table:1D:Burgers:regularization}
\end{table}
%
%
%Table \ref{table:2D:Burgers:SG:Regularization} displays the numerical results of the 2D Burgers equation after regularization. The methods SGL-1-R-2D, SGL-2-R-2D, and SGL-3-R-2D correspond to the single-grade methods SGL-1, SGL-2, and SGL-3 for the 2D equation, respectively, with $L_2$ regularization. The regularization parameter $\lambda_R$ is set to 1e-6 for the single-grade methods. The relative $L_2$ errors for these methods are 3.81e-04, 1.81e-04, and 2.68e-04, respectively. It is worth noting that TS-MGDL achieves the lowest relative $L_2$ error of 4.75e-05 among all the method.
%
%
\begin{table}[htpb]
		\centering
\begin{tabular}{c|c|c|c|c}
 \hline
   Methods &  Learning rate &  Decay rate & Epochs & Relative $L_2$ error  \\
   \hline
   SGL-1-R-2D & 3e-4  & 1e-4   & 150,000    & 3.81e-04  \\
   SGL-2-R-2D &   3e-4  & 1e-4   & 150,000   & 1.81e-04   \\ 
   SGL-3-R-2D &   3e-4  & 1e-4   & 150,000   & 2.68e-04  \\
    TS-MGDL  & -- & --    & 150,000  & \textbf{4.75e-05}\\
   \hline
\end{tabular}
  \caption{Comparison of regularized single-grade methods with TS-MGDL for the 2D Burgers equation.}
 \label{table:2D:Burgers:SG:Regularization}
\end{table}
%
%Table \ref{table:3D:Burgers:regularization} shows the numerical results for the 3D Burgers equation with regularization. The methods include SGL-1-R-3D, SGL-2-R-3D, and SGL-3-R-3D, which are single-grade methods with $L_2$ regularization, and their regularization parameters are set to 1e-5. Once again, TS-MGDL achieves the lowest relative $L_2$ error of 8.90e-05, surpassing the performance of the single-grade methods with regularization.

\begin{table}[htbp]
		\centering
\begin{tabular}{c|c|c|c|c}
 \hline
   Methods &  Learning rate &  Decay rate & Epochs & Relative $L_2$ error  \\
   \hline
   SGL-1-R-3D & 4e-4 &  1e-4 & 40,000   & 1.67e-04  \\
   SGL-2-R-3D & 4e-4 &  1e-4 & 40,000   & 2.67e-04  \\
   SGL-3-R-3D & 4e-4 &  1e-4 & 40,000   & 2.23e-04 \\
   TS-MGDL & --  & --   & 40,000  & \textbf{8.90e-05}\\
   \hline
\end{tabular}
  \caption{Comparison of regularized single-grade methods with TS-MGDL for the 3D Burgers equation.}
 \label{table:3D:Burgers:regularization}
\end{table}

\textbf{Intrinsic Multiscale Representation:} The TS-MGDL method enables an intrinsic multiscale representation of solutions. From the solution representation \eqref{Multiscale-pre} obtained through TS-MGDL, it is evident that different grades of the model capture solution components at varying scales, arranged from large to small. This hierarchical structure, combined with the incorporation of multiple grades, facilitates a more accurate and economically efficient representation of the solution. The lower grades are adept at capturing the global or low-frequency components of the solution. Given that the residuals of the solution typically exhibit high-frequency characteristics, subsequent grades learn the residuals by adding extra hidden layers. This progressive refinement allows the method to effectively handle problems with intrinsic multiscale features.
By combining multiple grades with different levels of refinement, the multi-grade method provides a comprehensive and adaptable approach to tackling such problems. The intrinsic multiscale representation of the TS-MGDL method enhances its ability to handle complex multiscale features, contributing to its effectiveness in solving challenging problems.

\begin{figure}[htbp]
\centering
\subfloat[$u_1^*$ \label{fig:ex1-grade-1}]{\includegraphics[width=0.25\textwidth]{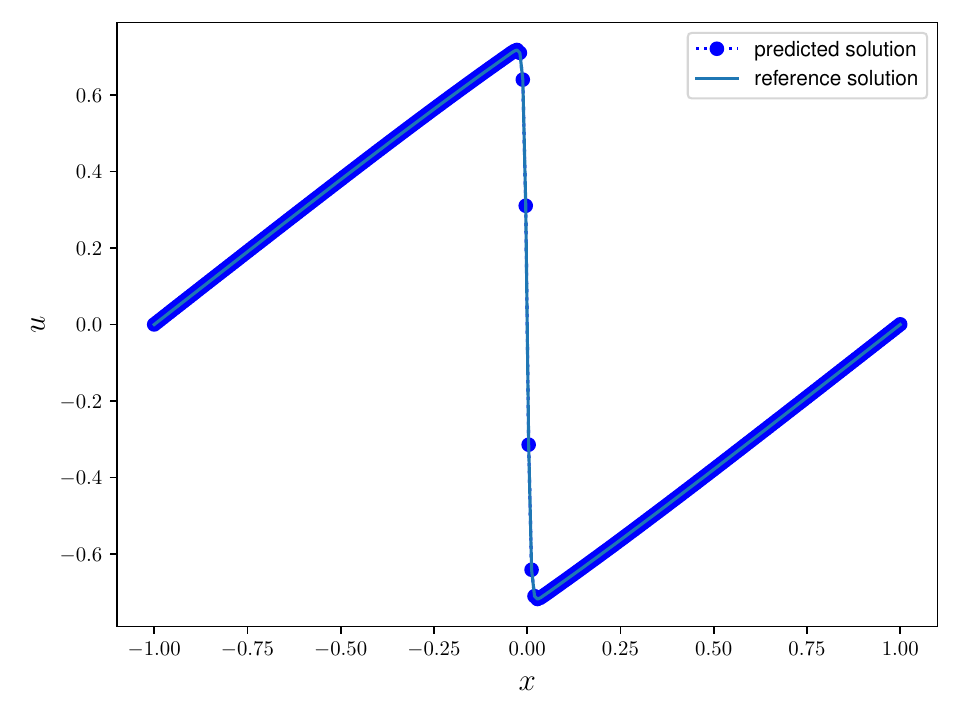}}\hfill
\subfloat[$u_2^*$ \label{fig:ex1-grade-2}]{\includegraphics[width=0.25\textwidth]{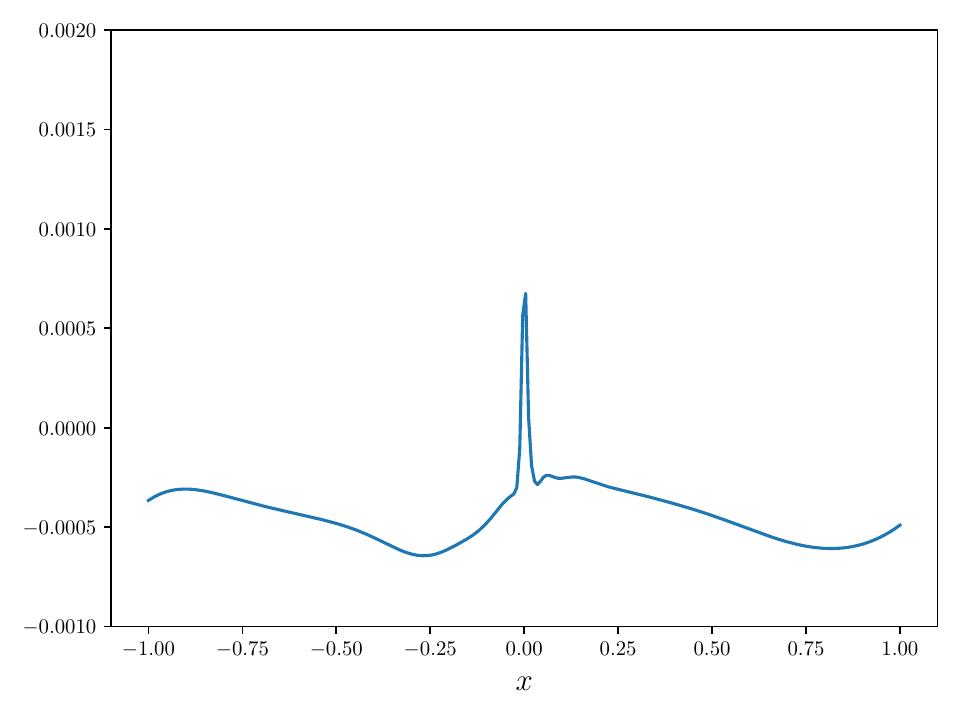}}\hfill
\subfloat[$\widetilde{u}_{3}^{*}$ \label{fig:ex1-stage2:add}]{\includegraphics[width=0.25\textwidth]{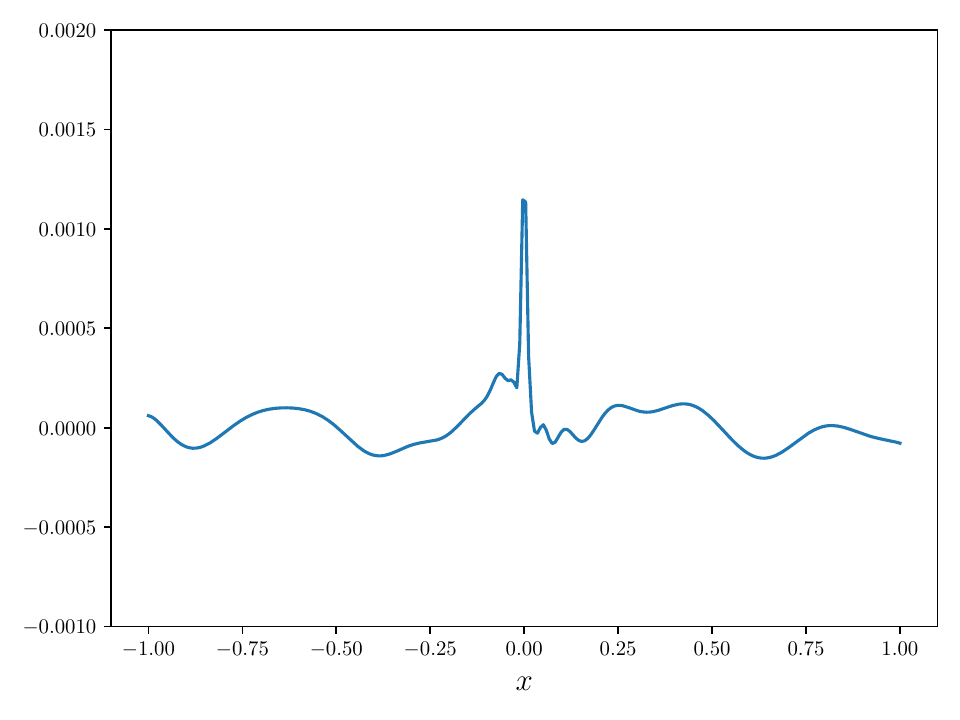}}\hfill
\caption{Numerical results of TS-MGDL for the 1D Burgers equation at $t=1$. (a) $u_1^*$ generated from Grade 1. (b)$u_2^*$ generated from Grade 2. (c) $\widetilde{u}_{3}^{*}$ generated from Stage 2.}
\label{fig:1d:TS-MGDL}
\end{figure}

\begin{figure}[htbp]
\centering
\subfloat[$u_1^*$]{\includegraphics[width=0.3\textwidth]{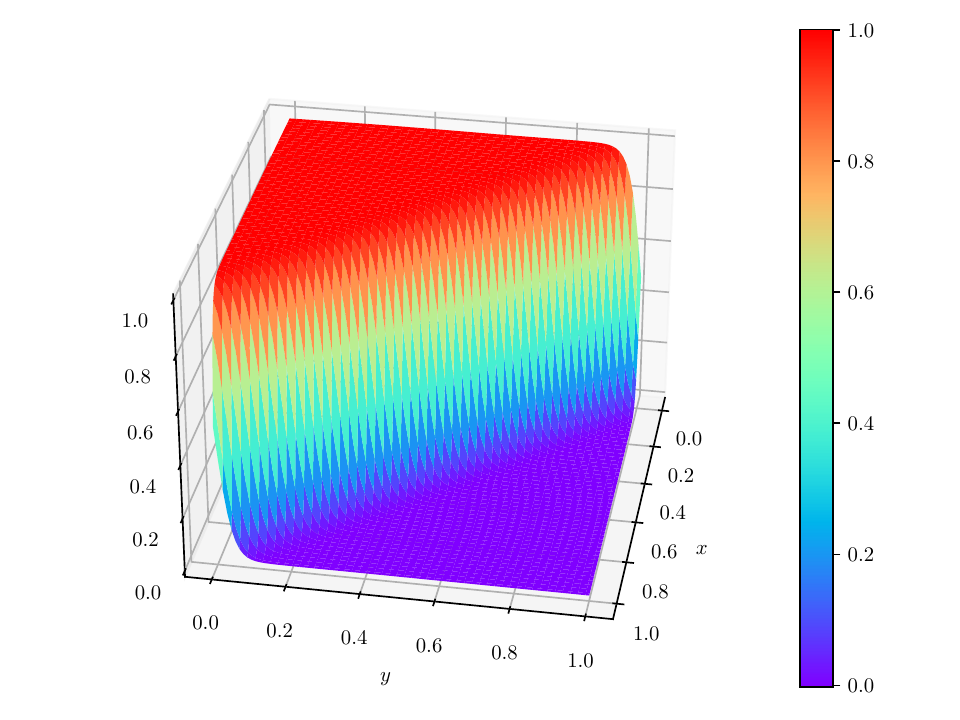}}\hfill
\subfloat[$u_2^*$]{\includegraphics[width=0.3\textwidth]{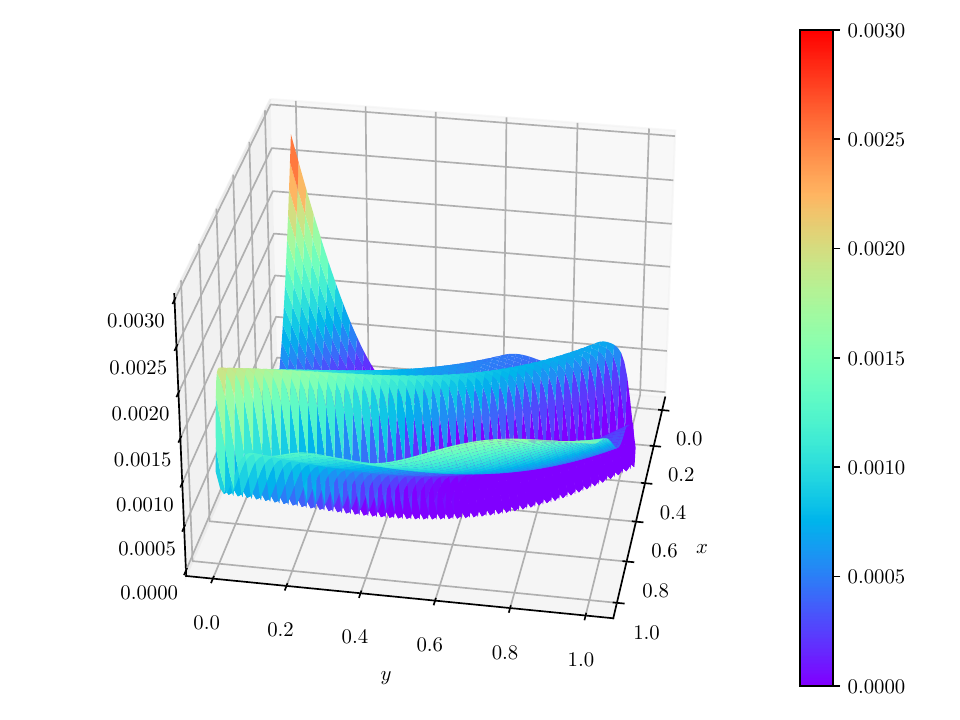}}\hfill
\subfloat[$\widetilde{u}_{3}^{*}$]{\includegraphics[width=0.3\textwidth]{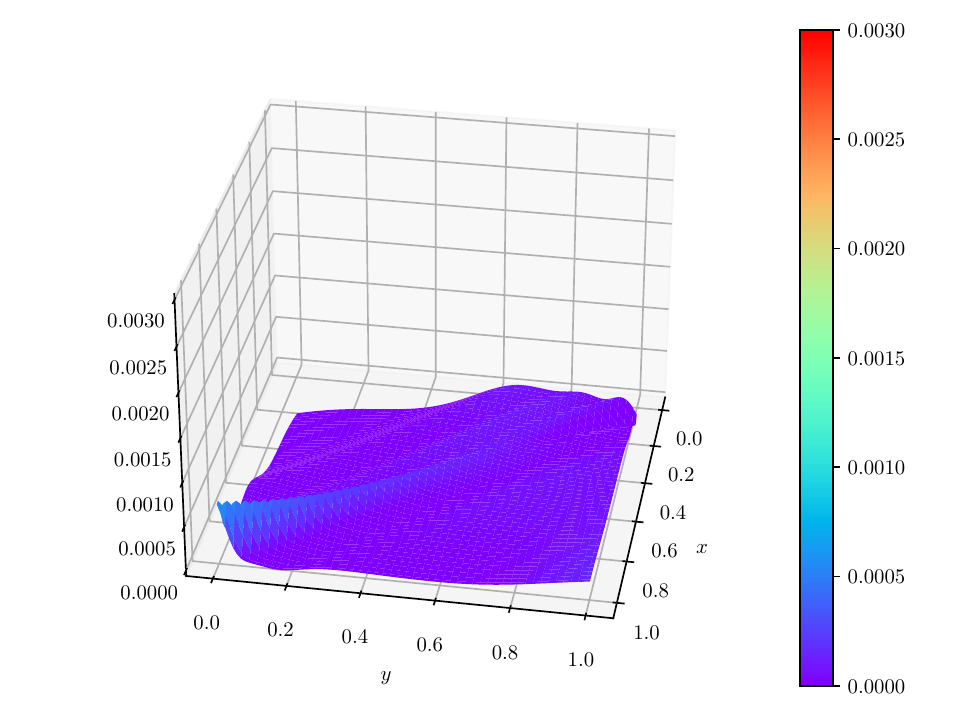}}\hfill
\caption{Numerical results of TS-MGDL for the 2D Burgers equation at $t=1$. (a) $u_1^*$ generated from Grade 1. (b)$u_2^*$ generated from Grade 2. (c) $\widetilde{u}_{3}^{*}$ generated from Stage 2.}
\label{fig:2d:TS-MGDL}
\end{figure}

We validate this point for the 1D and 2D cases in Figures \ref{fig:1d:TS-MGDL} and \ref{fig:2d:TS-MGDL}, respectively. In both cases, it can be observed that the first grade captures the main part of the solution, while the subsequent grades gradually extract the high-frequency and local features of the solution through residual learning.

%, showing that $\widetilde{u}_{3}^{*}$ contains richer information than $u_3^*$.

\textbf{Improved Accuracy:}  The TS-MGDL method enhances the accuracy of the learned solution via a combination of implicit regularization, multiscale representation, residual learning, and fine-tuning. By incorporating implicit regularization, the TS-MGDL method effectively prevents overfitting. Additionally, multiscale representation enables the model to handle solution sharp gradients. Furthermore, the inclusion of residual learning and fine-tuning techniques further contributes to improving the accuracy of the learned solution. All of these contribute to improved solution accuracy.  

Specifically, in the first stage, the network is trained incrementally, one grade at a time, progressively incorporating increasingly detailed information to improve accuracy. We can see from Tables \ref{table:1D:Burgers:MG}, \ref{table:2D:Burgers:MG}, and \ref{table:3D:Burgers:MG} that by capturing residual information, the TS-MGDL method captures intricate details and improves accuracy after adding a new grade. 
The second stage involves unfreezing and retraining, which explores a broader parameter space and refines the solution. By considering a wider range of possible solutions, the method increases the chance of finding better global optima, leading to improved accuracy.
It is confirmed from Tables \ref{table:1D:Burgers}, \ref{table:2D:Burgers}, and \ref{table:3D:Burgers} that after the second stage training, the TS-MGDL method achieves higher approximation accuracy compared to the single-grade method with the same number of epochs.
These results indicate that TS-MGDL exhibits superior accuracy in solving the Burgers equations of 1D, 2D, and 3D compared to single-grade methods with $L_2$ regularization.

\section{Conclusion}
A significant limitation of traditional DNNs is their inability to maintain or improve prediction accuracy as the network depth increases. To address this issue, we introduce a novel two-stage multi-grade deep learning (TS-MGDL) method specifically designed for solving nonlinear partial differential equations (PDEs). Unlike traditional DNNs, which operate as single-grade learning models, our multi-grade approach comprises a hierarchical superposition of neural networks arranged in a stair-like structure, with each level representing a distinct learning grade.

We have demonstrated that each grade and stage of the TS-MGDL method effectively reduces the loss function, a result validated through numerical experiments. Our findings show that the multi-grade approach significantly outperforms its single-grade counterpart by mitigating overfitting and achieving superior prediction accuracy. Consequently, the TS-MGDL method presents an innovative solution to the limitations of traditional single-grade models.

The numerical experiments further confirm that the TS-MGDL method is particularly well-suited for solving nonlinear PDEs with sharp gradients. Additionally, the study illustrates that DNNs can overcome the curse of dimensionality in PDE-solving tasks. Although the proposed TS-MGDL method has been applied primarily to the Burgers equation in this paper, it holds promise for broader applications to other types of nonlinear PDEs. While this work focuses primarily on the computational methodology and its numerical implementation, it also lays the foundation for exploring several critical theoretical issues in future research.

%\section*{Acknowledgments}
%We would like to acknowledge the assistance of volunteers in putting
%together this example manuscript and supplement.

\bibliographystyle{plain}
% \bibliography{Reference}

\end{document}